\documentclass[12pt,a4paper,reqno]{amsart}
\usepackage[T1]{fontenc}
\usepackage[utf8]{inputenc}
\usepackage[english]{babel}
\usepackage{lmodern,amsmath,amssymb,amsthm,mathtools}
\usepackage[margin=28mm,headheight=12pt,headsep=8mm]{geometry}
\usepackage{microtype,enumitem,tikz,pgfplots,mathrsfs,graphicx}
\usepackage[hidelinks,hypertexnames=false]{hyperref}
\allowdisplaybreaks[2]
\numberwithin{equation}{section}
\newtheorem{theorem}{Theorem}[section]
\newtheorem{lemma}[theorem]{Lemma}
\newtheorem{proposition}[theorem]{Proposition}
\newtheorem{corollary}[theorem]{Corollary}
\theoremstyle{definition}
\newtheorem{definition}[theorem]{Definition}
\theoremstyle{remark}

\newcommand{\R}{\mathbb R}
\newcommand{\N}{\mathbb N}
\newcommand{\Ftwo}{\mathbb F_2}
\newcommand{\Pp}{\mathcal P}
\newcommand{\Sp}{\mathcal S^+}
\newcommand{\Sb}{\overline{\mathcal S^+}}
\newcommand{\Cs}{\mathscr C_{\mathrm s}}
\newcommand{\F}{\mathfrak F}
\newcommand{\Lin}{\mathfrak L}
\newcommand{\Red}{\mathfrak G}
\newcommand{\Yred}{\mathcal Y}
\newcommand{\Null}[1]{N[#1]}
\newcommand{\Range}[1]{R[#1]}
\newcommand{\pr}{\mathcal P_\lambda}
\newcommand{\wpr}{\mathcal P_a}
\newcommand{\Reflection}{\mathscr R}
\newcommand{\A}{\mathcal A}
\newcommand{\norm}[1]{\lVert#1\rVert}
\newcommand{\supp}{\operatorname{supp}}
\newcommand{\dist}{\operatorname{dist}}
\newcommand{\ind}{\operatorname{ind}}
\setlist[enumerate]{label=\textup{(\roman*)},itemsep=4pt,topsep=5pt}
\title[Global topology and symmetry breaking]{Global topology and symmetry breaking in an indefinite superlinear Dirichlet problem}
\author{Juan Carlos Sampedro}
\address{Departamento de Matem\'atica Aplicada y Ciencias de la Computaci\'on, Avenida de los Castros 46, Universidad de Cantabria (UC), Santander, 39005, Spain}
\email{juancarlos.sampedro@unican.es}
\thanks{The author acknowledges support from the Spanish Ministry of Science and Innovation through grants PID2025-171414NA-I00 and PID2024-155890NB-I00.}
\date{}
\hypersetup{pdftitle={Global topology and symmetry breaking in an indefinite superlinear Dirichlet problem},pdfauthor={Juan Carlos Sampedro},pdfsubject={Positive solutions, global topology and symmetry breaking}}
\begin{document}
\begin{abstract}
Motivated by the work of Fencl and L\'opez-G\'omez~\cite{FL} on global
bifurcation diagrams, we study the connected components of positive solutions of $-u''=\lambda u+a(x)u^p$ on $(-1,1)$ with $u(-1)=u(1)=0$, where $a$ changes sign and $p\ge2$ is an integer. For generic $\mathcal C^r$ weights, $r\ge2$, with $n$ positive intervals, simple zeros and negative endpoint values, we obtain at least $2^{n-1}$ noncompact components. Apart from the principal component, each is a curve with both ends at $\lambda=-\infty$. The proof distinguishes at least $2^n-1$ ends and yields bounds for branching and compactly supported cohomology. For even weights, an odd Morse index of at least $k$ on the principal symmetric component forces at least $k$ symmetry-breaking points. In particular, symmetry breaking occurs if $a(0)<0$.
\end{abstract}
\maketitle

\section{Introduction}\label{Sec1}

We study how the positive solutions of an indefinite superlinear Dirichlet problem are distributed among connected components as the parameter varies. Fencl and L\'opez-G\'omez~\cite{FL} raised this question in their analytical
and numerical study of global bifurcation diagrams.

Throughout the paper, $p\in\N$ is fixed with $p\ge2$, and $a\in\mathcal C([-1,1])$ is a sign-changing weight. For $\lambda\in\R$, we consider the positive solutions of
\begin{equation}\label{Eq1.1}
	\begin{cases}
		-u''=\lambda u+a(x)u^p & \text{in }(-1,1),\\
		u(-1)=u(1)=0.
	\end{cases}
\end{equation}

Under the hypotheses below, a weight with $n$ positive intervals gives at least $2^n-1$ positive solutions for sufficiently negative $\lambda$ \cite[Th.~1.1]{FLS}. Several solutions at the same value of $\lambda$ may, however, lie on a single curve that turns back as $\lambda$ varies. This leaves open how many connected components are needed to account for this multiplicity. We give conditions on the weight under which the diagram has at least $2^{n-1}$ noncompact components, including the one that issues from the zero solution.

\begin{figure}[!htbp]
\centering
\begin{tikzpicture}[x=5.7cm,y=1.35cm,font=\small,>=stealth]
\fill[black!7] (-.75,-.48) rectangle (-.25,1.05);
\fill[black!7] (.25,-.48) rectangle (.75,1.05);
\draw[->,gray!70] (-1.04,0)--(1.07,0) node[right] {$x$};
\draw[->,gray!70] (0,-.55)--(0,1.18) node[above] {$a(x)$};
\draw[thick] (-1,-.38) .. controls (-.87,-.4) and (-.8,-.28) .. (-.75,0)
 .. controls (-.64,.92) and (-.38,.92) .. (-.25,0)
 .. controls (-.15,-.43) and (.15,-.43) .. (.25,0)
 .. controls (.38,.72) and (.64,.72) .. (.75,0)
 .. controls (.83,-.38) and (.92,-.43) .. (1,-.37);
\foreach \x/\lab in {-1/-1,-.75/\sigma_1,-.25/\tau_1,.25/\sigma_2,.75/\tau_2,1/1}{
 \draw (\x,-.035)--(\x,.035);
 \node[below] at (\x,-.51) {$\lab$};
}
\node at (-.5,1.02) {$I_1^+$};
\node at (.5,1.02) {$I_2^+$};
\node at (-.88,-.14) {$I_0^-$};
\node at (.1,-.18) {$I_1^-$};
\node at (.89,-.14) {$I_2^-$};
\end{tikzpicture}
\caption{A schematic weight with two positive intervals. The shaded intervals are $I_i^+=(\sigma_i,\tau_i)$, where $a>0$. The intervening and exterior intervals have $a<0$.}\label{Fig1}
\end{figure}
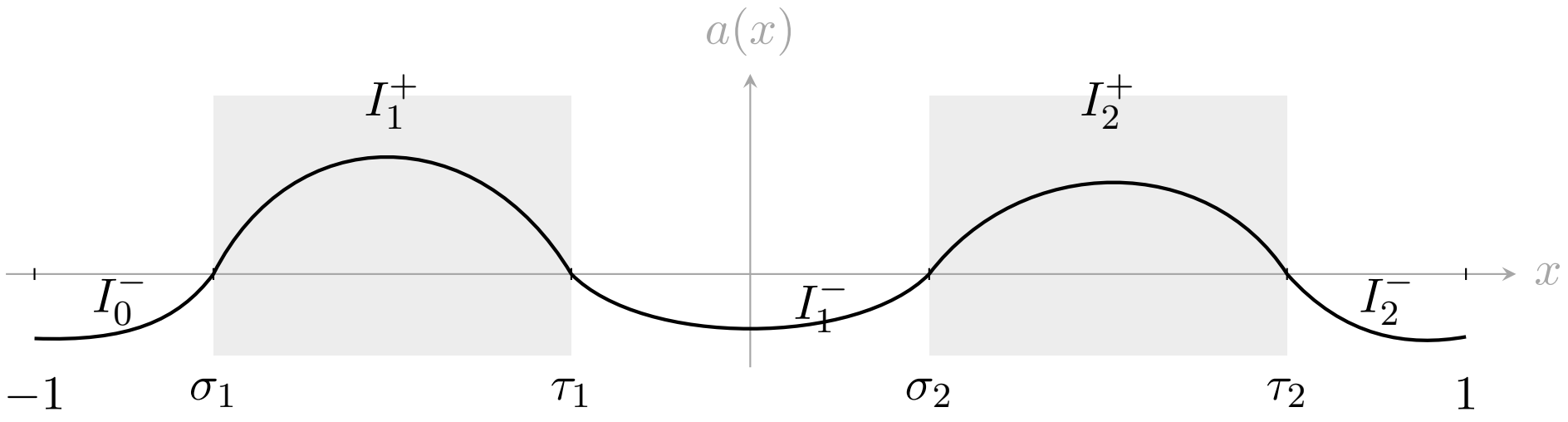

The following hypotheses from \cite{FLS} remain in force throughout the paper, with evenness of the weight assumed only where stated.
\begin{description}[style=unboxed,leftmargin=0pt,labelsep=.5em,font=\normalfont,itemsep=6pt,topsep=6pt]
\item[$(H_a)$] There are $n\ge1$ positive intervals $I_i^+=(\sigma_i,\tau_i)$ with
\[
 -1=\tau_0\le\sigma_1<\tau_1<\sigma_2<\cdots<\sigma_n<\tau_n\le\sigma_{n+1}=1.
\]
The weight is positive on $I_i^+$ for $1\le i\le n$ and negative on every nonempty interval $I_i^-=(\tau_i,\sigma_{i+1})$, $0\le i\le n$.
\item[$(H_g)$] On each $I_i^+$,
\[
 a(x)=\rho_i(x)\dist(x,\partial I_i^+)^{\gamma_i},
\]
where $\gamma_i>0$, $\rho_i\in \mathcal C([\sigma_i,\tau_i])$ is strictly positive, and $\dist$ denotes distance.
\end{description}
At least one negative interval is nonempty, and the fixed weight is assumed only continuous. Hypothesis $(H_g)$ specifies its vanishing at the endpoints of the positive intervals and is used in the a priori and degree estimates below. The arrangement of positive and negative intervals is illustrated in Figure~\ref{Fig1}.

Let $X=\{u\in \mathcal C^2([-1,1]):u(-1)=u(1)=0\}$ and $Y=\mathcal C([-1,1])$. Write $\Sp$ for the set of positive solution pairs $(\lambda,u)\in\R\times X$ of~\eqref{Eq1.1}, and set $\F_a(\lambda,u)=-u''-\lambda u-au^p$.
The closure of $\Sp$ in $\R\times X$ is denoted by $\Sb$. Let $\lambda_1=\pi^2/4$ be the first Dirichlet eigenvalue on $(-1,1)$ and put $\mathfrak{o}=(\lambda_1,0)$. Lemma~\ref{Lem2.2} will show that $\Sb=\Sp\cup\{\mathfrak{o}\}$ and that exactly one positive half-branch issues from $\mathfrak{o}$.

Let $G$ be the nonisolated part of $\Sb$, obtained by deleting its isolated points. We shall prove that $G$ is a locally finite graph, locally consisting of finitely many solution arcs meeting at their endpoints. Throughout the paper, components are understood to be connected, and the one containing $\mathfrak{o}$ is called the \emph{principal component}.

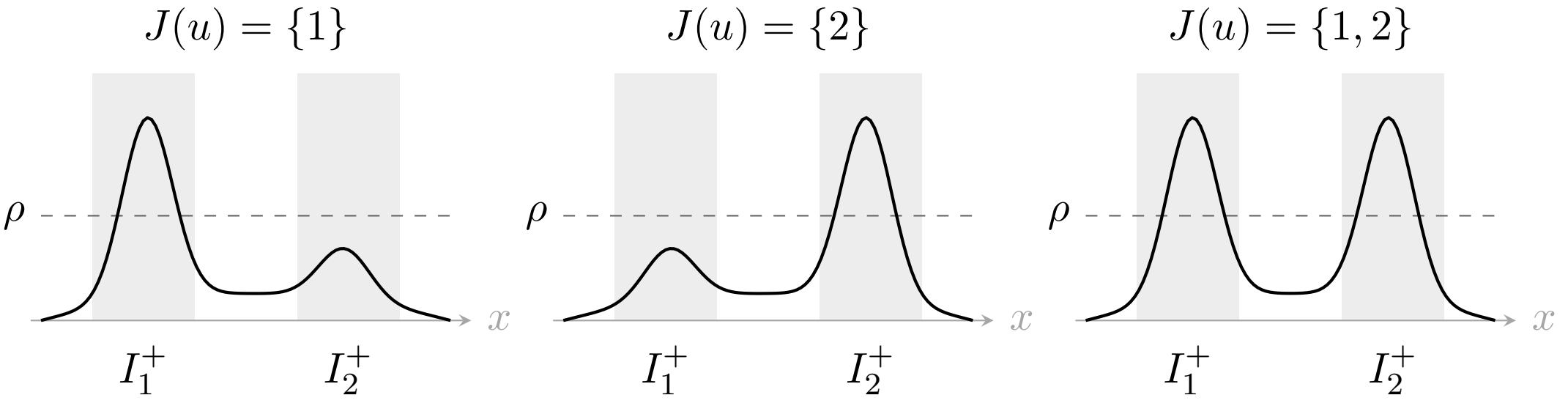
\begin{figure}[h!]
\centering
\begin{tikzpicture}[x=1.85cm,y=1.35cm,font=\small,>=stealth]
\foreach \shift/\aa/\bb/\lab in {0/1.6/.45/\{1\},2.55/.45/1.6/\{2\},5.1/1.6/1.6/{\{1,2\}}}{
\begin{scope}[shift={({\shift},0)}]
\fill[black!7] (-.75,0) rectangle (-.25,1.65);
\fill[black!7] (.25,0) rectangle (.75,1.65);
\draw[->,gray!70] (-1.05,0)--(1.1,0) node[right] {$x$};
\draw[dashed,black!60] (-1,.7)--(1,.7);
\node[left] at (-1,.7) {$\rho$};
\draw[thick,domain=-1:1,samples=101,variable=\x]
 plot ({\x},{(1-\x*\x)*(.18+\aa*exp(-((\x+.5)/.18)^2)+\bb*exp(-((\x-.5)/.18)^2))});
\node[above] at (0,1.7) {$J(u)=\lab$};
\node[below] at (-.5,-.06) {$I_1^+$};
\node[below] at (.5,-.06) {$I_2^+$};
\end{scope}}
\end{tikzpicture}
\caption{Reading the label $J(u)$ in~\eqref{Eq1.2}. An interval belongs to $J(u)$ precisely when the maximum of $u$ on it exceeds the dashed threshold $\rho$. The three schematic positive profiles illustrate all nonempty labels for two intervals.}\label{Fig2}
\end{figure}

By \cite{FLS}, there are $\rho>0$ and $\lambda_*<0$ such that, for $\lambda<\lambda_*$, every positive solution determines a nonempty index set
\begin{equation}\label{Eq1.2}
 J(u)=\{i\in\{1,\ldots,n\}:\norm u_{L^\infty(I_i^+)}>\rho\},
\end{equation}
and none of these interval norms equals $\rho$. The set $J(u)$ labels a solution by the intervals on which its maximum exceeds $\rho$, as illustrated in Figure~\ref{Fig2}. Since the threshold is excluded, this label is constant on every connected family with $\lambda<\lambda_*$. The first result shows that each label persists along a continuation to $\lambda=-\infty$. Two such curves determine different \emph{ends} if, for all sufficiently
negative $c$, their tails lie in different connected components of the part of $G$ with $\lambda<c$.

\begin{theorem}[Distinct ends]\label{Th1.1}
Under $(H_a)$ and $(H_g)$, for every nonempty $J\subset\{1,\ldots,n\}$ there is a continuous injective curve $r:[0,\infty)\to G$, $r(t)=(\lambda(t),u(t))$, such that
\[
 J(u(t))=J,\qquad
 \lim_{t\to\infty}\lambda(t)=-\infty,\qquad
 \lim_{t\to\infty}\norm{u(t)}_\infty=\infty.
\]
Curves carrying different labels determine different ends. In particular, $G$ has at least $2^n-1$ ends.
\end{theorem}

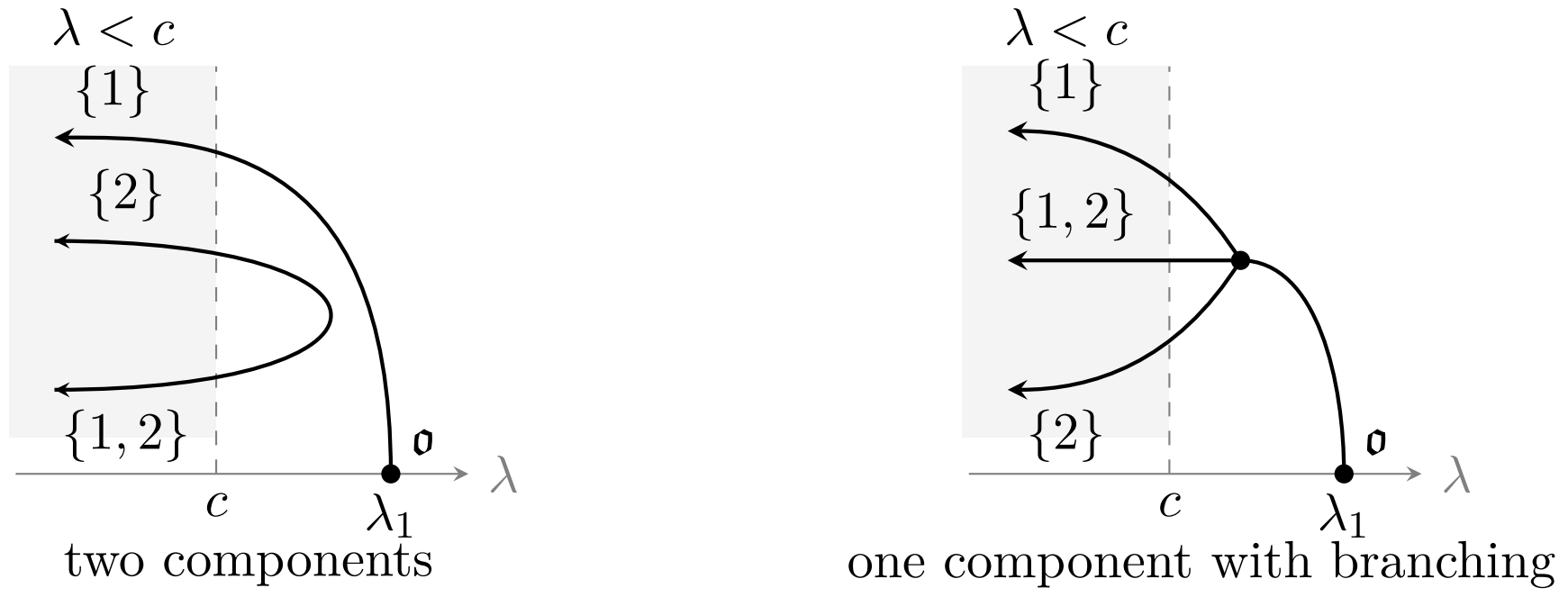
\begin{figure}[!htbp]
\centering
\begin{tikzpicture}[x=.95cm,y=.95cm,>=stealth,font=\small]
\begin{scope}[xshift=0cm]
\fill[gray!9] (-2.65,.08) rectangle (-1.05,2.95);
\draw[dashed,gray] (-1.05,-.2)--(-1.05,2.95);
\node[above] at (-1.85,2.95) {$\lambda<c$};
\draw[->,gray] (-2.6,-.2)--(.9,-.2) node[right] {$\lambda$};
\node[below] at (-1.05,-.2) {$c$};
\draw[thick,->] (.3,-.2)..controls (.3,2.4) and (-1.2,2.4)..(-2.3,2.4);
\fill (.3,-.2) circle (2pt) node[above right] {$\mathfrak{o}$};
\node[below] at (.3,-.2) {$\lambda_1$};
\node[above] at (-1.85,2.4) {$\{1\}$};
\draw[thick] (-2.3,.45)..controls (.55,.45) and (.55,1.6)..(-2.3,1.6);
\draw[->] (-1.95,.45)--(-2.3,.45);
\draw[->] (-1.95,1.6)--(-2.3,1.6);
\node[above] at (-1.75,1.6) {$\{2\}$};
\node[below] at (-1.75,.45) {$\{1,2\}$};
\node at (-.8,-.9) {two components};

\end{scope}
\begin{scope}[xshift=7cm]
\fill[gray!9] (-2.65,.08) rectangle (-1.05,2.95);
\draw[dashed,gray] (-1.05,-.2)--(-1.05,2.95);
\node[above] at (-1.85,2.95) {$\lambda<c$};
\draw[->,gray] (-2.6,-.2)--(.9,-.2) node[right] {$\lambda$};
\node[below] at (-1.05,-.2) {$c$};
\draw[thick] (.3,-.2)..controls (.3,.5) and (.1,1.45)..(-.5,1.45);
\fill (.3,-.2) circle (2pt) node[above right] {$\mathfrak{o}$};
\node[below] at (.3,-.2) {$\lambda_1$};
\fill (-.5,1.45) circle (2pt);
\draw[thick,->] (-.5,1.45)..controls (-1,2.25) and (-1.65,2.45)..(-2.3,2.45);
\draw[thick,->] (-.5,1.45)--(-2.3,1.45);
\draw[thick,->] (-.5,1.45)..controls (-1,.65) and (-1.65,.45)..(-2.3,.45);
\node[above] at (-1.85,2.45) {$\{1\}$};
\node[above] at (-1.8,1.45) {$\{1,2\}$};
\node[below] at (-1.85,.45) {$\{2\}$};
\node at (-.8,-.9) {one component with branching};

\end{scope}
\end{tikzpicture}
\caption{Two abstract graphs with the same three labeled ends. The cut $\lambda=c$ separates three tails continuing to $\lambda=-\infty$, shaded in each panel. The full graphs have two components on the left and one on the right.}\label{Fig3}
\end{figure}

For two positive intervals, the three labels are $\{1\}$, $\{2\}$ and $\{1,2\}$. Their ends may lie on separate components or meet through a branching point, as in Figure~\ref{Fig3}. Restricting to $\lambda<c$ separates three tails in either model, although the full graphs have different numbers of components. This relation between ends, components and branching is quantified in Section~\ref{Sec4}.

The labeled ends also determine classes in $H_c^1(G,\Ftwo)$, the first cohomology group with compact support and coefficients in $\Ftwo=\{0,1\}$, with addition modulo two. On a graph, this group records finite edge assignments modulo those obtained by changing values at finitely many vertices. Compact support retains information about separation at infinity, so a graph can have nonzero cohomology of this kind even in the absence of cycles. For example, $H_c^1(\R,\Ftwo)\simeq\Ftwo$ reflects the two ends of a line. Section~\ref{Subsec4.2} develops this description through vertex and edge functions and uses the interval labels to obtain quantitative constraints on the solution set. We introduce the graph-theoretic and cohomological tools in some detail, with readers working in differential equations in mind.

\begin{theorem}[Cohomological bounds]\label{Th1.2}
Under $(H_a)$ and $(H_g)$, each nonempty $J\subset\{1,\ldots,n\}$ determines a canonical class $\alpha_J\in H_c^1(G,\Ftwo)$, and
\[
 \sum_{J\ne\varnothing}\alpha_J=0,\qquad
 \dim H_c^1(G,\Ftwo)\ge2^{n-1}-1.
\]
If $G$ is connected, this is the only nontrivial linear relation among the $\alpha_J$. In particular,
\[
 \dim\operatorname{span}\{\alpha_J:J\ne\varnothing\}=2^n-2,
 \qquad \dim H_c^1(G,\Ftwo)\ge2^n-2.
\]
\end{theorem}

To obtain a lower bound on the number of components, we vary the weight in the following class. For an integer $r\ge2$, let $\A_n$ consist of all $a\in\mathcal C^r([-1,1])$ with $a(-1)<0$, $a(1)<0$ and exactly $2n$ simple zeros in $(-1,1)$, where a zero $x$ is simple if $a'(x)\ne0$. This class is open in $\mathcal C^r([-1,1])$, and every weight in it satisfies $(H_a)$ and $(H_g)$. Here generic means on a residual subset of $\A_n$, one containing a countable intersection of open dense subsets.

At a positive solution, the derivative with respect to $u$ is the Dirichlet linearization $\Lin_{\lambda,u}:=D_u\F_a(\lambda,u):X\to Y$, given by
\begin{equation}\label{Eq1.3}
 \Lin_{\lambda,u}v=-v''-\lambda v-pa(x)u^{p-1}v,\qquad v\in X.
\end{equation}
A \emph{turning point} is a strict local extremum of the parameter along a local solution curve. Let $K$ denote the number of noncompact components of $G$.

\begin{theorem}[Generic components]\label{Th1.3}
For a residual subset $\mathcal R_n\subset\A_n$, the set $\Sp$ is a one-dimensional analytic manifold with
\begin{equation}\label{Eq1.4}
 K\ge2^{n-1}.
\end{equation}
The principal component, including $\mathfrak{o}$, is homeomorphic to $[0,\infty)$. Every other noncompact component is homeomorphic to $\R$, has $\lambda\to-\infty$ at both ends and contains a turning point where $\Lin_{\lambda,u}$ is singular.

If $K=2^{n-1}$, each nonempty $J\subset\{1,\ldots,n\}$ labels exactly one negative end.
\end{theorem}

In general, the exact number of components remains undetermined, and the theorem permits compact components and several turning points on a curve, possibly of higher even order.

For even weights, we also study bifurcation from the family of symmetric solutions. Define $\Reflection u(x)=u(-x)$ and let $\Cs$ be the component of
\[
 \{(\lambda,u)\in\Sp:\Reflection u=u\}\cup\{\mathfrak{o}\}
\]
containing $\mathfrak{o}$. An even positive solution is a \emph{symmetry-breaking point} if it is a limit of positive solutions $u$ with $\Reflection u\ne u$. On $\Cs$, this is a \emph{secondary bifurcation} along the nontrivial symmetric family arising from the primary bifurcation at $\mathfrak{o}$.

At an even solution, the Dirichlet linearization preserves odd functions. Its number of negative eigenvalues on that subspace is the \emph{odd Morse index}, denoted by $m_-(\lambda,u)$.

\begin{theorem}[Secondary bifurcation]\label{Th1.4}
Assume $(H_a)$ and $(H_g)$, with $a$ even. If $\Cs$ contains a solution with $m_-(\lambda,u)\ge k$, where $k\ge1$ is an integer, then it contains at least $k$ distinct symmetry-breaking points, each with at least two even and two nonsymmetric incident half-arcs. In particular, $a(0)<0$ forces symmetry breaking on $\Cs$.
\end{theorem}

The degree computations in \cite{FZ,FLS}, together with the a priori estimates of \cite{AG}, provide the starting point for our analysis. Exact multiplicity for a different family, parametrized by the size of the negative part of the weight, is obtained in \cite{FT}.

Related questions for solutions with one interior zero of logistic equations are studied by L\'opez-G\'omez and Rabinowitz~\cite{LR15,LR20}. Spatial heterogeneity produces multiplicity in this nodal setting~\cite{LR16}, including for strictly positive weights~\cite{CLT2022}. Numerical studies by L\'opez-G\'omez, Molina-Meyer and Rabinowitz~\cite{LGMR2017} and Cubillos, L\'opez-G\'omez and Tellini~\cite{CLT2023} show how branches reconnect and separate into different components as the weight varies.

For positive solutions, Cubillos, L\'opez-G\'omez and Tellini~\cite{CLT2024} consider nonnegative piecewise constant weights with vanishing intervals, formulate an exponential multiplicity conjecture and compute global diagrams with secondary branches. For a weight with two positive intervals separated by a central zero interval, L\'opez-G\'omez, Mu\~noz-Hern\'andez and Zanolin~\cite[Ths.~2.1 and~4.1]{LMZ2025} prove uniqueness of the symmetric positive solution and existence of a reflected pair of asymmetric positive solutions for sufficiently negative $\lambda$.

Our continuation arguments build on classical bifurcation theory~\cite{CR,Rabinowitz,Dancer}, unilateral global bifurcation~\cite{LG2001}, and the analytic graph framework of L\'opez-G\'omez and Sampedro~\cite[Sec.~7.1, Th.~7.4]{LGB}, applied there to the same polynomial equation. Compactness, finite parameter fibers and the local branch at $\mathfrak{o}$ yield the solution graph, whose parts with $\lambda$ bounded below are compact. We follow each degree class to an individual end, quantify the components or branching needed to accommodate the ends, and construct the cohomology classes carried by their interval labels. The passage from connected continua to paths requires justification, as illustrated by~\cite{Toland}.

Generic regularity under perturbations has precedents in~\cite{Rynne,BHJY}. Here it converts the end count into a lower bound on the number of noncompact components. For even weights, the lower bound for the odd Morse index supplied by $J(u)$ yields a count of distinct symmetry-breaking points on $\Cs$ and the criterion $a(0)<0$, complementing~\cite{Tanaka,Mandel,PW}.

Sections~\ref{Sec2}--\ref{Sec4} construct the graph and prove the topological bounds. Sections~\ref{Sec5} and~\ref{Sec6} treat generic components and symmetry breaking, followed by examples in Section~\ref{Sec7}.

\section{Analytic structure of the positive solution set}\label{Sec2}
Continuation along curves requires compactness on bounded parameter intervals and a local description of how solution branches meet. We first establish compactness, then combine analytic reduction with Cauchy uniqueness to prove finiteness of the solutions at each fixed parameter.

\subsection{Functional setting}
We endow $X$ and $Y$ with $\norm u_X=\sum_{j=0}^2\norm{u^{(j)}}_\infty$ and $\norm f_Y=\norm f_\infty$, respectively. Here $\norm{\cdot}_\infty$ denotes the uniform norm on $[-1,1]$. For functions on this interval, we write $(v,w)_2=\int_{-1}^1vw\,dx$ whenever the integral is defined.

For bounded linear operators $T:X\to Y$, we follow \cite{LGB}, using $\Null{T}$ and $\Range{T}$ for the kernel and range, $GL(X,Y)$ for the topological isomorphisms, and $\Phi_0(X,Y)$ for the Fredholm operators of index zero. A Fredholm operator has closed range, finite-dimensional kernel and finite-codimensional range, and its index is $\ind T:=\dim\Null{T}-\operatorname{codim}_Y\Range{T}$. Thus $T\in\Phi_0(X,Y)$ satisfies
\[
 \dim\Null{T}=\operatorname{codim}_Y\Range{T}<\infty.
\]
A function in $X$ is \emph{strongly positive} if it belongs to the open set
\[
 \Pp=\{u\in X:u>0\text{ in }(-1,1),\ u'(-1)>0,\ u'(1)<0\}.
\]
Define the nonlinear operator and its positive solution set by
\[
 \F_a:\R\times X\to Y,\qquad
 \F_a(\lambda,u):=-u''-\lambda u-au^p,\qquad
 \Sp=\F_a^{-1}(0)\cap(\R\times\Pp).
\]
Every positive solution satisfies the boundary-slope conditions by Cauchy uniqueness, so this definition agrees with that in the introduction. The derivative of $\F_a$ with respect to $u$ is the Dirichlet linearization introduced in~\eqref{Eq1.3},
\[
 \Lin_{\lambda,u}=D_u\F_a(\lambda,u):X\to Y,\qquad
 \Lin_{\lambda,u}v=-v''-\lambda v-pa(x)u^{p-1}v.
\]

\begin{lemma}\label{Lem2.1}
The map $\F_a:\R\times X\to Y$ is real analytic. At each positive solution, $\Lin_{\lambda,u}\in\Phi_0(X,Y)$ and $\dim\Null{\Lin_{\lambda,u}}\le1$.
\end{lemma}
\begin{proof}
Pointwise multiplication is continuous from $X^p$ to $Y$, so $u\mapsto au^p$ is a continuous homogeneous polynomial and $\F_a$ is jointly real analytic. The operator $-\partial_x^2:X\to Y$ is an isomorphism, while multiplication by $\lambda+pa u^{p-1}$ is compact from $X$ to $Y$ through the compact inclusion $X\hookrightarrow Y$. Hence $\Lin_{\lambda,u}$ is Fredholm of index zero. The linear map
\[
 \ell:\Null{\Lin_{\lambda,u}}\longrightarrow\R,\qquad \ell(v):=v'(-1),
\]
is injective because a kernel element with $v(-1)=v'(-1)=0$ vanishes identically by Cauchy uniqueness. The kernel of $\Lin_{\lambda,u}$ therefore has dimension at most one.
\end{proof}

For the degree arguments, we use a compact formulation that selects the positive solutions. The inverse Dirichlet operator is
\[
 \mathcal T:Y\longrightarrow X,\qquad
 (\mathcal Tf)(x):=\int_{-1}^1 k_D(x,y)f(y)\,dy,
\]
where the Green kernel on $[-1,1]$ is
\[
 k_D(x,y):=\frac12
 \begin{cases}
 (1+x)(1-y),&x\le y,\\
 (1+y)(1-x),&y\le x.
 \end{cases}
\]
Differentiation shows that $\mathcal Tf$ solves $-u''=f$, and the kernel vanishes at $x=\pm1$, giving the zero Dirichlet data. If $j:X\hookrightarrow Y$ is the inclusion, then
\[
 T:=j\mathcal T:Y\longrightarrow Y
\]
is compact. Writing $v^+(x)=\max\{v(x),0\}$, define
\[
 \Phi_\lambda:Y\longrightarrow Y,\qquad
 \Phi_\lambda(v):=T\bigl(\lambda v^++a(v^+)^p\bigr).
\]
Every fixed point lies in $X$. On a component of $\{v<0\}$ it would be affine with zero values at both endpoints, which rules out such a component. A nonzero fixed point is therefore nonnegative, and Cauchy uniqueness makes it strongly positive. Conversely, every positive solution of~\eqref{Eq1.1} is a fixed point of $\Phi_\lambda$. Thus
\[
 \F_a(\lambda,u)=0,\quad u\in\Pp
 \quad\Longleftrightarrow\quad u=\Phi_\lambda u,\quad u\ne0.
\]
For a bounded open set $\Omega\subset Y$ with no fixed point on its boundary, we denote the Leray--Schauder degree by $\deg_{\mathrm{LS}}(I-\Phi_\lambda,\Omega)$. We use the normalization $\deg_{\mathrm{LS}}(I,\Omega)=1$ when $0\in\Omega$. A nonzero degree guarantees a fixed point in $\Omega$, and the degree remains unchanged as $\lambda$ varies over a bounded interval provided no fixed point meets $\partial\Omega$.

The following a priori and degree estimates will be used throughout the argument.
\begin{enumerate}[label=\textup{(E\arabic*)}]
\item\label{E1} Positive solutions are uniformly bounded in $L^\infty$ on compact parameter intervals, and their parameters have a finite upper bound.
\item\label{E2} For $\lambda<\lambda_*$, the set $J(u)$ in~\eqref{Eq1.2} is nonempty and no interval norm equals $\rho$.
\item\label{E3} For $\lambda<\lambda_*$, choose $R_\lambda>\rho$ larger than every solution norm at that parameter. For each nonempty $J\subset\{1,\ldots,n\}$, the region
\[
 \Lambda^J_\lambda=\left\{v\in Y:
 \begin{array}{l}
 \norm v_\infty<R_\lambda,\\
 \norm v_{L^\infty(I_i^+)}>\rho\ \text{if }i\in J,\\
 \norm v_{L^\infty(I_i^+)}<\rho\ \text{if }i\notin J
 \end{array}\right\}
\]
contains no fixed point on its boundary and satisfies
\begin{equation}\label{Eq2.1}
 d_J:=\deg_{\mathrm{LS}}(I-\Phi_\lambda,\Lambda^J_\lambda)=(-1)^{|J|}.
\end{equation}
\item\label{E4} As $\lambda\to-\infty$, solutions tend uniformly to zero on compact subintervals of negative intervals, uniformly over all positive solutions.
\end{enumerate}
The uniform bound on compact parameter intervals in \ref{E1} is the one recalled for this Dirichlet problem in \cite[Sec.~1]{FLS}. The bound used in the degree homotopies at fixed $\lambda<0$ is given separately in \cite[Lem.~2.3]{FLS}.

The upper parameter bound in \ref{E1} can be obtained directly by testing on a positive interval. Let $\varphi_i>0$ be a first Dirichlet eigenfunction on $I_i^+$ and write $\nu_i=\pi^2/(\tau_i-\sigma_i)^2$. Integration by parts gives
\[
 (\lambda-\nu_i)\int_{I_i^+}u\varphi_i\,dx
 +\int_{I_i^+}a u^p\varphi_i\,dx
 =u(\tau_i)\varphi_i'(\tau_i)-u(\sigma_i)\varphi_i'(\sigma_i)\le0.
\]
The second integral is positive, hence $\lambda<\nu_i$ for every $i$ and $\lambda<\min_i\nu_i$.

The threshold exclusion in \ref{E2}, the degree formula in \ref{E3}, and the decay in \ref{E4} follow from \cite[Lem.~2.5, Prop.~3.3 and Lem.~2.4]{FLS}, respectively. For $\lambda<0$, the global maximum lies in a positive interval and satisfies the lower bound in \cite[Lem.~2.2]{FLS}. Decreasing $\lambda_*$ until this bound exceeds $\rho$ ensures that $J(u)$ is nonempty. The decay estimate originates in \cite[Th.~3.1]{FL}.

\subsection{Compactness and finiteness at fixed parameter}\label{Subsec2.2}
We consider solutions with $b\le\lambda\le c$, referring to $[b,c]\times X$ as a parameter strip. The compactness argument must account for the possible convergence of positive solutions to zero.

\begin{lemma}\label{Lem2.2}
One has $\Sb=\Sp\cup\{\mathfrak{o}\}$, and $\Sb\cap([b,c]\times X)$ is compact for every bounded closed interval $[b,c]$. Near $\mathfrak{o}$ there is exactly one positive analytic half-branch.
\end{lemma}
\begin{proof}
Consider a sequence $\{(\lambda_j,u_j)\}_{j\ge1}\subset\Sp$ with $\lambda_j\in[b,c]$, and put $L=\max\{|b|,|c|\}$. By \ref{E1}, there is $M>0$ such that $\norm{u_j}_\infty\le M$. The equation gives
\[
 \norm{u_j''}_\infty
 \le LM+\norm a_\infty M^p=:A.
\]
Rolle's theorem provides $\xi_j\in(-1,1)$ with $u_j'(\xi_j)=0$, so
\[
 |u_j'(x)|=\left|\int_{\xi_j}^x u_j''(y)\,dy\right|\le2A,
 \qquad \norm{u_j}_X\le M+3A.
\]
After passing to a subsequence and retaining the same notation, compactness of $[b,c]$ and the inclusion $X\hookrightarrow\mathcal C^1([-1,1])$ give $\lambda_j\to\lambda$ and $u_j\to u$ in $\mathcal C^1$. Since
\[
 \lim_{j\to\infty}u_j''
 =\lim_{j\to\infty}(-\lambda_j u_j-au_j^p)=-\lambda u-au^p
 \quad\text{in }Y,
\]
the limit belongs to $X$, the convergence holds in $X$, and $u$ solves the equation with $u\ge0$. If $u\not\equiv0$, an interior zero or a zero boundary slope would give zero Cauchy data and hence $u\equiv0$. Thus every nonzero limit is strongly positive.

To identify a zero limit, suppose that $\{u_j\}_{j\ge1}$ converges to zero and define the normalized sequence $\{v_j\}_{j\ge1}$ by $v_j=u_j/\varepsilon_j$, where $\varepsilon_j=\norm{u_j}_\infty$. Then
\[
 \begin{cases}
 -v_j''=\lambda_jv_j+a\varepsilon_j^{p-1}v_j^p,\\
 v_j(-1)=v_j(1)=0,
 \end{cases}
 \qquad v_j\ge0,\quad\norm{v_j}_\infty=1.
\]
Because $p\ge2$, the nonlinear term tends to zero in $Y$. The preceding compactness argument applies to $\{v_j\}_{j\ge1}$, and a subsequence converges in $X$ to a function $v\ge0$ with $\norm v_\infty=1$ and $-v''=\lambda v$. Let $\varphi_1>0$ be the first Dirichlet eigenfunction normalized by $(\varphi_1,\varphi_1)_2=1$. Integration by parts yields
\[
 (\lambda-\lambda_1)(v,\varphi_1)_2=0,
 \qquad (v,\varphi_1)_2>0.
\]
Consequently $\lambda=\lambda_1$ and $v=\varphi_1/\norm{\varphi_1}_\infty$. This shows that $\mathfrak{o}$ is the only possible limit with zero solution.

To construct the analytic branch through $\mathfrak{o}$, write $u=sz$ with $(z,\varphi_1)_2=1$ and introduce the polynomial map
\[
 \begin{gathered}
 \mathcal H:\R\times\R\times X\longrightarrow Y\times\R,\\
 \mathcal H(s,\lambda,z):=
 \bigl(-z''-\lambda z-s^{p-1}az^p,\ (z,\varphi_1)_2-1\bigr).
 \end{gathered}
\]
Analyticity follows from Lemma~\ref{Lem2.1}, and $\mathcal H(0,\lambda_1,\varphi_1)=0$. Its derivative in $(\lambda,z)$ at this point is
\[
 \begin{gathered}
 \mathcal B_0:\R\times X\longrightarrow Y\times\R,\\
 \mathcal B_0(\mu,w):=
 \bigl((-\partial_x^2-\lambda_1)w-\mu\varphi_1,
       (w,\varphi_1)_2\bigr).
 \end{gathered}
\]
For given $(f,\eta)\in Y\times\R$, the equation $\mathcal B_0(\mu,w)=(f,\eta)$ is solved by
\[
 \begin{aligned}
 \mu&=-(f,\varphi_1)_2,\qquad w=\eta\varphi_1+w_\perp,\\
 (-\partial_x^2-\lambda_1)w_\perp
 &=f-(f,\varphi_1)_2\varphi_1,\qquad
 (w_\perp,\varphi_1)_2=0.
 \end{aligned}
\]
The Fredholm alternative gives a unique $w_\perp\in X$, depending continuously on $f$. Hence $\mathcal B_0$ is an isomorphism. The analytic implicit function theorem now gives unique analytic functions $\lambda(s)$ and $z(s)$ near $s=0$, with
\[
 \lambda(0)=\lambda_1,\qquad z(0)=\varphi_1,\qquad
 \mathcal H(s,\lambda(s),z(s))=0.
\]
Since $z(0)=\varphi_1\in\Pp$ and $\Pp$ is open in $X$, we may assume $z(s)\in\Pp$ for $|s|<\delta$, with $\delta>0$. Multiplication of the normalized equation by $s$ gives the analytic solution curve
\[
 \gamma:(-\delta,\delta)\longrightarrow\R\times X,\qquad
 \gamma(s):=\bigl(\lambda(s),\ s z(s)\bigr),\qquad \gamma(0)=\mathfrak{o}.
\]
The normalization gives $(s z(s),\varphi_1)_2=s$, so $\gamma$ is injective and $s$ is a coordinate along the curve. Its tangent at zero is $(\lambda'(0),\varphi_1)\ne0$, and its positive solutions are exactly those with $s>0$. In particular, $\mathfrak{o}\in\Sb$, which proves $\Sb=\Sp\cup\{\mathfrak{o}\}$. The sequential compactness established above then gives compactness of the solution set over every bounded closed parameter interval.

Finally, suppose that a sequence $\{(\lambda_j,u_j)\}_{j\ge1}\subset\Sp$ outside this curve converges to $\mathfrak{o}$. With $\varepsilon_j$ and $v_j$ as above, pass to a subsequence such that $\lim_{j\to\infty}v_j=\varphi_1/\norm{\varphi_1}_\infty$ in $X$. Put $s_j=(u_j,\varphi_1)_2>0$ and $z_j=u_j/s_j$, so that
\[
 \begin{aligned}
 \lim_{j\to\infty}\frac{s_j}{\varepsilon_j}
 &=\lim_{j\to\infty}(v_j,\varphi_1)_2
 =\frac1{\norm{\varphi_1}_\infty},\\
 \lim_{j\to\infty}z_j
 &=\lim_{j\to\infty}\frac{v_j}{(v_j,\varphi_1)_2}
 =\varphi_1\quad\text{in }X.
 \end{aligned}
\]
Since $\mathcal H(s_j,\lambda_j,z_j)=0$, local uniqueness in the implicit function theorem gives $\lambda_j=\lambda(s_j)$ and $z_j=z(s_j)$ for large $j$. This contradicts the choice of the sequence and proves uniqueness of the positive half-branch.
\end{proof}

To study solutions at a fixed parameter and local branches at singular points, we use the analytic Lyapunov--Schmidt reduction of \cite{LGB}. Its projections have an explicit form in this self-adjoint setting. At a positive solution $(\lambda_0,u_0)$ where $\Lin_0:=\Lin_{\lambda_0,u_0}$ is singular, choose $\phi\in\Null{\Lin_0}$ with $(\phi,\phi)_2=1$. The kernel and range projections are
\[
 \begin{aligned}
 P&:X\longrightarrow X,& P(v)&:=(v,\phi)_2\phi,\\
 Q&:Y\longrightarrow Y,& Q(f)&:=f-(f,\phi)_2\phi.
 \end{aligned}
\]
The Fredholm alternative gives
\[
 \begin{aligned}
 X_1:=\Null{P}&=\{v\in X:(v,\phi)_2=0\},\qquad
 X=\Null{\Lin_0}\oplus X_1,\\
 Y_0:=\Range{\Lin_0}&=\{f\in Y:(f,\phi)_2=0\},\qquad
 Y=Y_0\oplus\operatorname{span}\{\phi\}.
 \end{aligned}
\]
Thus $(P,Q)$ is a pair of $\Lin_0$-projections in the notation of \cite{LGB}, with $P(X)=\Null{\Lin_0}$ and $Q(Y)=Y_0$. To solve for the component orthogonal to the kernel, define
\[
 \F_{\mathrm r}:\R^2\times X_1\longrightarrow Y_0,\qquad
 \F_{\mathrm r}(\lambda,t,w):=Q\F_a(\lambda,u_0+t\phi+w).
\]
Its derivative in $w$ at $(\lambda_0,0,0)$ is $Q\Lin_0|_{X_1}\in GL(X_1,Y_0)$. The analytic implicit function theorem therefore solves $\F_{\mathrm r}(\lambda,t,w)=0$ uniquely for small $w\in X_1$, giving an analytic map
\[
 \Yred:\mathcal U\longrightarrow X_1,\qquad
 \mathcal U:=I\times(-\delta,\delta)\subset\R^2,
\]
where $I$ is an open interval containing $\lambda_0$ and $\delta>0$. It satisfies $\Yred(\lambda_0,0)=0$ and $\partial_t\Yred(\lambda_0,0)=0$. The original equation now reduces to the scalar equation $\Red(\lambda,t)=0$, where
\begin{equation}\label{Eq2.2}
\Red:\mathcal U\longrightarrow\R,\quad
\Red(\lambda,t):=\bigl(\F_a(\lambda,u_0+t\phi+\Yred(\lambda,t)),\phi\bigr)_2.
\end{equation}

After restricting $\mathcal U$ and a neighborhood $\mathcal W$ of $(\lambda_0,u_0)$ in $\R\times X$, write
\[
 \mathcal S_{\mathrm{loc}}:=\F_a^{-1}(0)\cap\mathcal W,\qquad
 Z_{\mathrm{loc}}:=\{(\lambda,t)\in\mathcal U:\Red(\lambda,t)=0\}.
\]
The maps
\begin{equation}\label{Eq2.3}
 \begin{aligned}
 \Psi&:\mathcal S_{\mathrm{loc}}\longrightarrow Z_{\mathrm{loc}},&
 \Psi(\lambda,u)&:=(\lambda,(u-u_0,\phi)_2),\\
 \Psi^{-1}&:Z_{\mathrm{loc}}\longrightarrow\mathcal S_{\mathrm{loc}},&
 \Psi^{-1}(\lambda,t)&:=(\lambda,u_0+t\phi+\Yred(\lambda,t))
 \end{aligned}
\end{equation}
are inverse analytic maps, with $t$ recording the kernel coordinate of a solution. The reduction therefore preserves both the parameter and the local topology of the solution set.

\begin{lemma}\label{Lem2.3}
For each fixed $\lambda\in\R$, problem~\eqref{Eq1.1} has finitely many positive solutions.
\end{lemma}
\begin{proof}
Fix $\bar\lambda\in\R$. By Lemma~\ref{Lem2.2}, the fiber $\mathcal K=\{u:(\bar\lambda,u)\in\Sb\}$ is compact, as is its set $D$ of accumulation points. Suppose that $D$ contains a positive solution, and choose $u\in D$ maximizing the initial slope $u'(-1)$. This maximum is positive, so $u\in\Pp$.

The linearization at $u$ must be singular, since an invertible linearization would isolate $u$ in the fiber. Let $\phi$ span its kernel, with $(\phi,\phi)_2=1$. Under the local correspondence~\eqref{Eq2.3}, distinct solutions converging to $u$ give distinct zeros of $t\mapsto\Red(\bar\lambda,t)$ converging to zero. Analyticity forces this function to vanish identically near zero. Hence
\[
 U(t)=u+t\phi+\Yred(\bar\lambda,t),\qquad |t|<\delta,
\]
is an arc of positive solutions at $\bar\lambda$. Since $(U(t)-u,\phi)_2=t$, every point of the arc is an accumulation point of the fiber. For the initial slope functional $\ell(v)=v'(-1)$, differentiation at zero gives
\[
 \left.\frac{d}{dt}\ell(U(t))\right|_{t=0}
 =\phi'(-1)\ne0.
\]
Here $\partial_t\Yred(\bar\lambda,0)=0$, and Cauchy uniqueness gives the inequality because $\phi(-1)=0$ and $\phi$ is nonzero. The arc thus contains a point of $D$ with larger initial slope than $u$, a contradiction.

The only remaining possible accumulation point is $u=0$ at $\bar\lambda=\lambda_1$. On the local branch of Lemma~\ref{Lem2.2}, a sequence in this fiber converging to zero would give distinct $s_j\downarrow0$ after passage to a subsequence, with $\lambda(s_j)=\lambda_1$. By analyticity $\lambda(s)=\lambda_1$ near zero. The positive half-branch would then consist of positive accumulation points of the fiber, which we have just excluded. Thus $D$ is empty and the compact fiber is finite.
\end{proof}

\subsection{The solution graph}\label{Subsec2.3}
We now assemble the local solution curves into the graph suggested by Figure~\ref{Fig3}. We shall show that every nontrivial nonisolated solution has an even number of incident half-arcs, whereas $\mathfrak{o}$ has a single incident half-arc. To describe this structure, we recall the standard notions needed below for readers unfamiliar with graph theory, referring to~\cite[Sec.~1]{Diestel} for locally finite graphs and their ends.

A \emph{graph} consists of vertices joined by edges. Formally, we attach copies of $[0,1]$ to a discrete vertex set $V$ by identifying their endpoints with vertices, and use the topology obtained by this gluing. The interval interiors are the \emph{open edges}. Writing $E$ for their collection gives
\[
 \Gamma=V\sqcup\bigsqcup_{e\in E}e,\qquad
 \text{each }e\text{ is homeomorphic to }(0,1).
\]
For a specified graph $H$, we write $V(H)$ and $E(H)$ for its vertices and open edges. Different edges may join the same vertices, and an edge with coinciding endpoints is a \emph{loop}. In the solution graph, an edge represents a portion of a solution curve, with the choice of endpoints determining its subdivision into pieces.

A \emph{simple edge path} follows successive edges without repeating a vertex. A \emph{simple cycle} is a closed edge path with no repeated vertex apart from its coinciding initial and final vertices, allowing a loop as a cycle with one edge. A \emph{tree} is a connected graph with no simple cycle. In a tree, any two vertices are joined by exactly one simple edge path. A \emph{rooted tree} is a tree with one specified vertex, called its root. These notions apply equally to finite and infinite graphs.

Removing a point $z$ from a sufficiently small connected neighborhood leaves its incident \emph{half-arcs}. Their number is the \emph{valence}, denoted by $\deg_G(z)$, and is defined independently of whether $z$ is chosen as a vertex. At a vertex, a loop is counted twice. A point of valence one is an endpoint, also called a \emph{leaf}, and one of valence at least three is a \emph{branching point}. Interior points of edges have valence two, so all leaves and branching points must be vertices. The graph is \emph{locally finite} when every vertex has finite valence. We denote its set of branching points by
\[
 \mathcal B:=\{v\in G:\deg_G(v)\ge3\}.
\]
The valence is intrinsic to the solution space.

Inserting a vertex inside an edge is a \emph{subdivision}. We use locally finite subdivisions, adding only finitely many vertices in each compact subset. Marking the integers on a half-line is one example. The compactness restriction excludes accumulating vertex sets such as $\{1/j:j\ge1\}$ in $[0,1]$, which would be incompatible with the usual interval topology.

Recall that $G$ is the nonisolated part of $\Sb$, obtained by deleting its isolated points. To study finite pieces of $G$, we use the parameter projection and its upper-level sets,
\[
 \pr:\Sb\longrightarrow\R,\qquad \pr(\lambda,u):=\lambda,
 \qquad G_{\ge c}:=G\cap\pr^{-1}([c,\infty)).
\]
We shall choose the vertices once and for all so that $G_{\ge c}$ is a finite graph after adding the finitely many points on the cut $\lambda=c$. These sets exhaust $G$ as $c$ tends to $-\infty$.

\begin{proposition}\label{Prop2.4}
The space $G$ is a countable locally finite graph. Its vertices and edges can be chosen so that every compact subset meets only finitely many of them. Every vertex other than $\mathfrak{o}$ has even valence, while $\mathfrak{o}$ has valence one. Each set $G_{\ge c}$ is compact, and the principal component is noncompact.
\end{proposition}
\begin{proof}
At a positive solution $z=(\lambda_0,u_0)$, invertibility of $\Lin_{\lambda_0,u_0}$ gives a local analytic graph over $\lambda$ and hence two incident half-arcs. When the linearization is singular, the correspondence~\eqref{Eq2.3} instead identifies the nearby solutions with
\[
 Z_z:=\{(\lambda,t)\in\mathcal U_z:\Red(\lambda,t)=0\}.
\]
Here $\mathcal U_z$ is a small neighborhood of $q_z=(\lambda_0,0)$ in the reduced plane. The analytic function $\Red$ cannot vanish identically, since fixing $\lambda=\lambda_0$ would contradict finiteness at $\lambda_0$ in Lemma~\ref{Lem2.3}. The local structure theorem for real analytic plane curves~\cite[Cor.~2.7, Lem.~6.3]{BM} therefore makes $Z_z$ either an isolated point or finitely many half-arcs meeting only at $q_z$. In the latter case, we count the half-arcs by intersecting $Z_z$ with a small circle around $q_z$. This intersection is called the \emph{link},
\[
 \operatorname{Lk}_z:=Z_z\cap
 \{(\lambda,t):(\lambda-\lambda_0)^2+t^2=\varepsilon^2\}.
\]
For sufficiently small $\varepsilon>0$, each half-arc meets the circle once, even when two half-arcs have the same tangent. The local analytic arc description~\cite[Lem.~6.3]{BM} makes the link independent of this small radius up to homeomorphism. Figure~\ref{Fig4} illustrates the count.

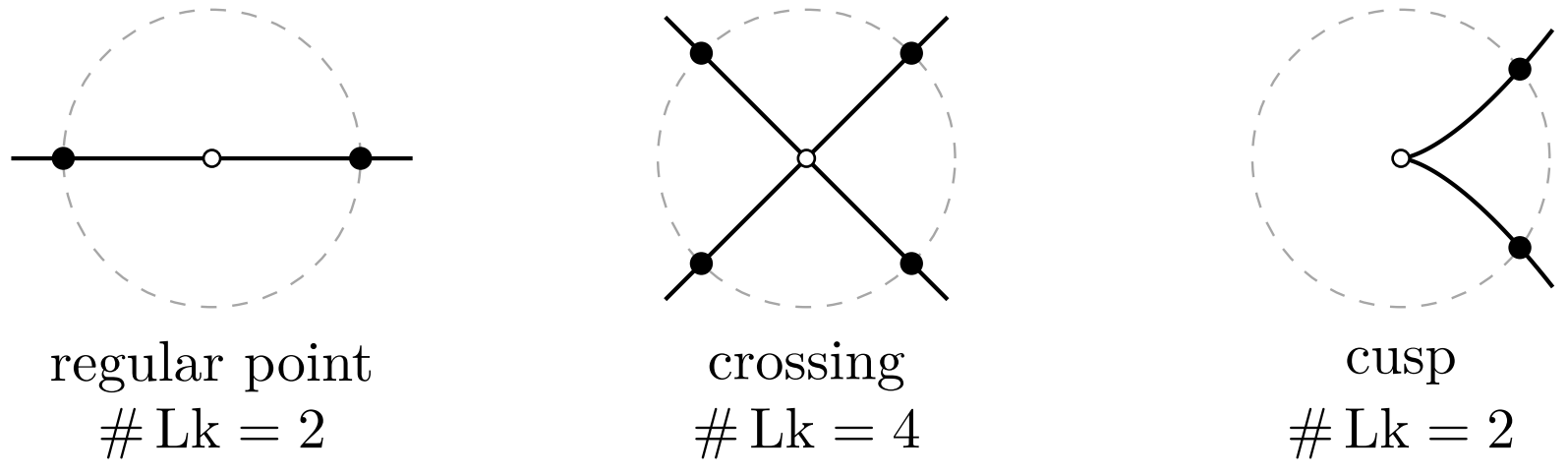
\begin{figure}[!htbp]
\centering
\begin{tikzpicture}[x=1cm,y=1cm,font=\small]
\begin{scope}
\draw[dashed,gray!70,line width=.45pt] (0,0) circle (1);
\draw[thick] (-1.35,0)--(1.35,0);
\draw[fill=white,line width=.5pt] (0,0) circle (1.6pt);
\fill (-1,0) circle (2.2pt);
\fill (1,0) circle (2.2pt);
\node at (0,-1.4) {regular point};
\node at (0,-1.85) {$\#\operatorname{Lk}=2$};
\end{scope}
\begin{scope}[xshift=4cm]
\draw[dashed,gray!70,line width=.45pt] (0,0) circle (1);
\draw[thick] (-.95,-.95)--(.95,.95);
\draw[thick] (-.95,.95)--(.95,-.95);
\draw[fill=white,line width=.5pt] (0,0) circle (1.6pt);
\foreach \x/\y in {.707107/.707107,-.707107/.707107,.707107/-.707107,-.707107/-.707107}
 {\fill (\x,\y) circle (2.2pt);}
\node at (0,-1.4) {crossing};
\node at (0,-1.85) {$\#\operatorname{Lk}=4$};
\end{scope}
\begin{scope}[xshift=8cm]
\draw[dashed,gray!70,line width=.45pt] (0,0) circle (1);
\draw[thick,domain=-1.13:1.13,samples=81,variable=\s]
 plot ({.8*\s*\s},{.6*\s*\s*\s});
\draw[fill=white,line width=.5pt] (0,0) circle (1.6pt);
\fill (.8,.6) circle (2.2pt);
\fill (.8,-.6) circle (2.2pt);
\node at (0,-1.4) {cusp};
\node at (0,-1.85) {$\#\operatorname{Lk}=2$};
\end{scope}
\end{tikzpicture}
\caption{The open circle marks the center of each local model in the reduced plane. The filled points on the dashed circle form its link, one for each incident half-arc. The cusp has two half-arcs with the same tangent, so its link still has two points.}\label{Fig4}
\end{figure}

Sullivan's theorem supplies the parity of this finite set through its \emph{Euler characteristic}. For a finite graph $H$ this invariant is
\[
 \chi(H)=|V(H)|-|E(H)|,
\]
which is unchanged by subdivision, since inserting a vertex also adds one edge. For a finite set $L$ there are no edges, so $\chi(L)=\#L$. Sullivan's \emph{local Euler condition}~\cite{Sullivan}, also proved in~\cite[Th.~4.4]{FuMcCrory}, states that the link of a real analytic set has even Euler characteristic. We need only the case of a finite link, where it says that the number of points is even.

Strong positivity is open in $X$, so every nearby solution represented by $Z_z$ remains positive and contributes to the link count. It follows that
\[
 \deg_G(z)=\#\operatorname{Lk}_z
          =\chi(\operatorname{Lk}_z)\equiv0\pmod2.
\]
Here $\deg_G(z)$ counts incident half-arcs and will be the valence in the graph constructed below. At $\mathfrak{o}$, Lemma~\ref{Lem2.2} gives exactly one positive half-branch, hence valence one.

To assemble the local arcs over bounded parameter intervals, we also need to control their parameter coordinate. For an analytic parametrization $\gamma(s)$ of a half-arc issuing from a positive solution $z$, Lemma~\ref{Lem2.3} excludes a constant parameter. Hence there are $m\ge1$ and $\beta\ne0$ such that
\[
 \pr(\gamma(s))=\lambda_0+\beta s^m+O(s^{m+1}).
\]
The parameter is consequently strictly monotone for sufficiently small $s>0$. The analytic parametrization from Lemma~\ref{Lem2.2} gives the same conclusion at $\mathfrak{o}$, since a constant $\lambda(s)$ would again contradict finiteness at $\lambda_1$.

Each local chart can therefore be chosen so that its center is the only possible isolated point, branching point or local extremum of the parameter. For every bounded closed interval $[b,c]$, compactness gives a finite covering of $\Sb\cap([b,c]\times X)$ by such charts, so only finitely many of these exceptional points occur there. Cutting at them and at the finitely many solutions on the strip boundaries leaves finitely many interval pieces, with strictly monotone parameter on their interiors, together with finitely many isolated solutions. The same local description prevents isolated solutions from accumulating at any point of $\Sb$, proving that $G$ is closed. Combined with the upper parameter bound and Lemma~\ref{Lem2.2}, this also proves that every $G_{\ge c}$ is compact.

We can now specify the vertices. In addition to the branching points in $\mathcal B$, let $\mathcal T_\lambda$ denote the valence-two points at which $\pr$ has a strict local extremum. Such a point lies on a local curve, although that curve need not be smooth at the point. Set
\begin{equation}\label{Eq2.4}
 V:=\{\mathfrak{o}\}\cup\mathcal B\cup\mathcal T_\lambda
       \cup\bigcup_{j\ge1}\bigl(G\cap\pr^{-1}(\{-j\})\bigr).
\end{equation}
The open edges are the connected components of $G\setminus V$. On each of them the parameter is strictly monotone. The integer cuts keep an edge inside a bounded parameter interval, while the preceding finite-arc description shows that its closure is a compact arc with endpoints in $V$. Every compact subset $Q\subset G$ lies in some $G_{\ge c}$, where only finitely many points of $\mathcal B\cup\mathcal T_\lambda$ and finitely many integer levels occur. Each level contains finitely many solutions by Lemma~\ref{Lem2.3}, hence
\begin{equation}\label{Eq2.5}
 \#(V\cap Q)<\infty\qquad\text{for every compact }Q\subset G.
\end{equation}
The finite-arc description in a slightly larger parameter strip also shows that $Q$ meets only finitely many edges. Gluing the closed arcs at their endpoints therefore recovers the original topology of $G$. With this fixed choice, the finite subgraphs $G_{\ge-j}$ exhaust $G$, proving countability and local finiteness without further choices of vertices.

The valence calculation now forces the principal component to be noncompact. Otherwise it would be a finite graph $C$ with $\mathfrak{o}$ as its only odd-valence vertex. Counting the two endpoints of every edge, including a loop twice, would give
\[
0\equiv2|E(C)|=\sum_{v\in V(C)}\deg_C(v)
\equiv\deg_C(\mathfrak{o})=1\pmod2.
\]
This contradiction proves that the principal component is noncompact.
\end{proof}

From now on, $V$ is the vertex set in~\eqref{Eq2.4}, and $E$ is the collection of connected components of $G\setminus V$, its open edges. We allow the locally finite subdivisions specified above when a cut is needed. The negative integers serve only to divide unbounded curves into compact edges. Any other such choice gives the same solution space, with the same components, ends and valences. Two locally finite decompositions have a common refinement obtained by taking the union of their vertex sets, which still satisfies~\eqref{Eq2.5}.

For a cut to separate each local solution arc into an upper and a lower part, its level must avoid branching points, isolated solutions and turns of the parameter. We use the following terminology.

\begin{definition}\label{Def2.5}
A real number $c$ is an \emph{exceptional level} if the slice $\lambda=c$ contains an isolated point of $\Sb$, the point $\mathfrak{o}$, a branching point of $G$, or a valence-two point where $\pr$ has a strict local extremum. All other levels are called \emph{nonexceptional}.
\end{definition}

Only finitely many exceptional levels occur in each bounded parameter interval, so nonexceptional levels exist arbitrarily far to the left. At each solution on such a level, the graph is locally a single arc and $\pr$ is strictly monotone on it, separating the arc into a part above and a part below $c$. 

To describe unbounded continuations within a noncompact component, we use \emph{rays}, that is, embeddings $r:[0,\infty)\to G$. Here an embedding is a continuous injective map whose inverse on its image is continuous, and we also call the image a ray. For $T\ge0$, its remaining part $r([T,\infty))$ is called a \emph{tail}.

To express escape from compact sets, we require the ray to be \emph{proper}, meaning that $r^{-1}(Q)$ is compact for every compact subset $Q\subset G$. Since this inverse image is already closed, the requirement is precisely that it be bounded, or equivalently,
\[
 \text{for every compact }Q\subset G\text{ there is }T_Q\ge0
 \text{ such that }r([T_Q,\infty))\cap Q=\varnothing.
\]
Writing $r(t)=(\lambda(t),u(t))$, we can read this condition directly from the parameter. Properness makes $r^{-1}(G_{\ge c})$ bounded for every $c\in\R$, so eventually $\lambda(t)<c$. Conversely, every nonempty compact $Q\subset G$ is contained in some $G_{\ge c}$, because $\pr$ attains a minimum on $Q$. Escape to $\lambda=-\infty$ therefore makes $r^{-1}(Q)$ closed and bounded. We obtain
\begin{equation}\label{Eq2.6}
 r\text{ is proper}\quad\Longleftrightarrow\quad
 \lim_{t\to\infty}\lambda(t)=-\infty.
\end{equation}
A proper ray may still have turning points, provided its parameter tends to $-\infty$.

To distinguish the ways in which proper rays escape, we compare their tails in
\[
 G_{<c}:=G\setminus G_{\ge c}=\{(\lambda,u)\in G:\lambda<c\}.
\]
By~\eqref{Eq2.6}, each proper ray has a tail in $G_{<c}$ for every $c<0$. For each fixed $c$, all sufficiently late tails lie in one component of $G_{<c}$, because they are connected and nested. We set
\begin{equation}\label{Eq2.7}
 r_1\sim r_2\quad\Longleftrightarrow\quad
 \begin{gathered}
 \text{for every }c<0,\text{ sufficiently late tails of }r_1\text{ and }r_2\\
 \text{lie in the same connected component of }G_{<c}.
 \end{gathered}
\end{equation}
This defines an equivalence relation on proper rays, whose classes are called the \emph{ends} of $G$.

This is the usual notion of an end defined by removal of compact sets~\cite[Sec.~1]{Diestel}. Indeed, for a compact $Q\subset G$, let $C_Q(r)$ denote the component of $G\setminus Q$ containing any sufficiently late tail of $r$. If $r_1\sim r_2$, choose $c<0$ with $Q\subset G_{\ge c}$. The component of $G_{<c}$ containing both tails lies in a single component of $G\setminus Q$, so $C_Q(r_1)=C_Q(r_2)$. Conversely, if this equality holds for every compact $Q$, taking $Q=G_{\ge c}$ gives~\eqref{Eq2.7}. Compact subsets meet only finitely many edges, so finite subgraphs suffice for these comparisons.

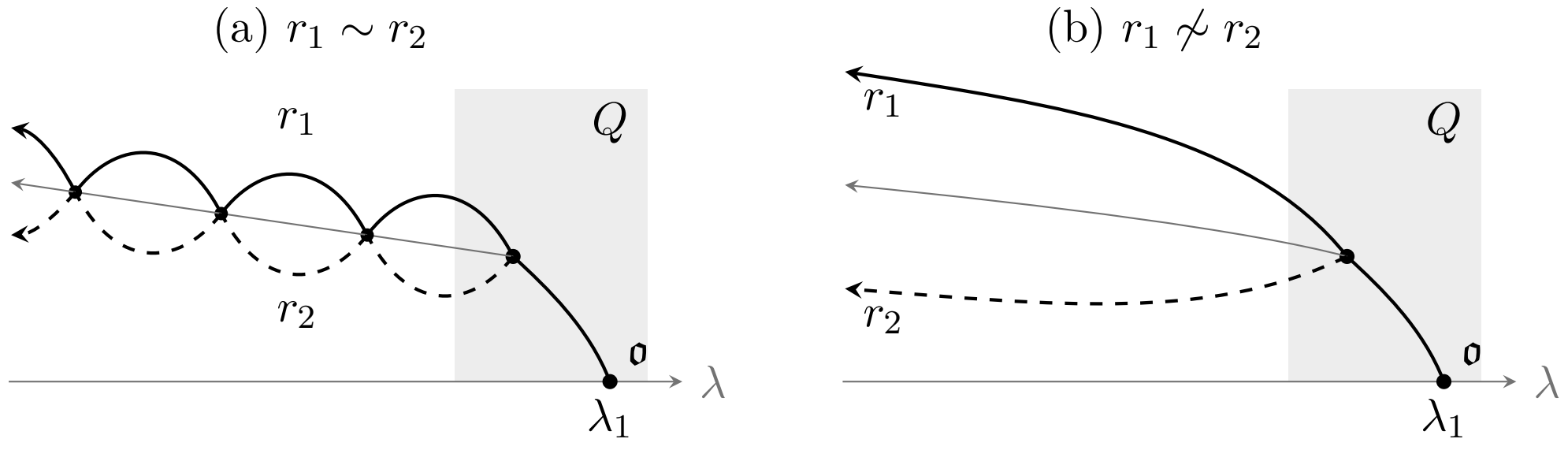
\begin{figure}[!htbp]
\centering
\begin{tikzpicture}[x=.98cm,y=1cm,font=\small,>=stealth]
\foreach \shift/\panel in {0/a,7/b}{
\begin{scope}[xshift=\shift cm]
\node at (2.7,2.95) {\if a\panel (a) $r_1\sim r_2$\else (b) $r_1\not\sim r_2$\fi};
\fill[black!7] (3.85,0) rectangle (5.5,2.45);
\node at (5.18,2.17) {$Q$};
\draw[->,black!55] (.03,0)--(5.8,0) node[right] {$\lambda$};
\draw[thick] (5.18,0) .. controls (5.0,.45) and (4.65,.78) .. (4.35,1.05);
\fill (5.18,0) circle (1.8pt) node[above right] {$\mathfrak{o}$};
\node[below] at (5.18,0) {$\lambda_1$};
\fill (4.35,1.05) circle (1.8pt);
\if a\panel
\foreach \a/\h/\b/\k in {4.35/1.05/3.1/1.23,3.1/1.23/1.85/1.41,1.85/1.41/.6/1.59}{
\draw[thick] (\a,\h) .. controls (\a-.35,\h+.67) and (\b+.35,\k+.45) .. (\b,\k);
\draw[black!55] (\a,\h)--(\b,\k);
\draw[thick,dashed] (\a,\h) .. controls (\a-.35,\h-.45) and (\b+.35,\k-.67) .. (\b,\k);
\fill (\b,\k) circle (1.6pt);
}
\draw[thick,->] (.6,1.59) .. controls (.42,1.94) and (.23,2.10) .. (.05,2.13);
\draw[black!55,->] (.6,1.59)--(.05,1.67);
\draw[thick,dashed,->] (.6,1.59) .. controls (.42,1.37) and (.23,1.24) .. (.05,1.23);
\node at (2.50,2.18) {$r_1$};
\node at (2.50,.56) {$r_2$};
\else
\draw[thick,->] (4.35,1.05) .. controls (3.52,2.08) and (1.9,2.32) .. (.05,2.60) node[below right] {$r_1$};
\draw[black!55,->] (4.35,1.05) .. controls (3.5,1.27) and (1.8,1.48) .. (.05,1.65);
\draw[thick,dashed,->] (4.35,1.05) .. controls (3.2,.50) and (1.85,.64) .. (.05,.78) node[below right] {$r_2$};
\fi
\end{scope}}
\end{tikzpicture}
\caption{Two rays in schematic solution graphs. In (a) the pattern of three arcs between successive vertices repeats towards $\lambda=-\infty$, so $r_1$ and $r_2$ represent the same end. In (b) their last common vertex lies in $Q$, and they represent distinct ends.}\label{Fig5}
\end{figure}

In Figure~\ref{Fig5}(a), the rays meet again arbitrarily far to the left, so their tails belong to the same component outside any compact set. In (b), removing the shaded part separates their tails.

By~\eqref{Eq2.6}, all ends of $G$ are \emph{negative ends}, and we denote their number by $e_-(G)$. A component homeomorphic to $[0,\infty)$ has one end, whereas a component homeomorphic to $\R$ has two. Ends describe escape to infinity, so compact components and the finite leaf $\mathfrak{o}$ contribute none.

\subsection{Local indices and isolation}
The degree arguments below concern isolated solutions at fixed $\lambda$ and isolated collections of components in a strip. For a positive solution $(\lambda,u)$, choose an open ball $B_\varepsilon(u)\subset Y$ whose closure contains no other fixed point of $\Phi_\lambda$, and define its \emph{local index} by
\begin{equation}\label{Eq2.8}
 \iota_\lambda(u)=\deg_{\mathrm{LS}}(I-\Phi_\lambda,B_\varepsilon(u)).
\end{equation}
Lemma~\ref{Lem2.3} ensures that such a ball exists, and excision makes the index independent of sufficiently small $\varepsilon$.

Consider $\mathcal C^1_0=\{v\in \mathcal C^1([-1,1]):v(-1)=v(1)=0\}$ with its usual norm. Since $T$ maps bounded subsets of $Y$ into bounded subsets of $X$, its realization in $\mathcal C^1_0$ is compact. The reduction-of-space property of Leray--Schauder degree therefore allows the same local index to be computed in $\mathcal C^1_0$. This space is convenient because strong positivity, including the boundary-slope conditions, is open in its norm. Near a positive solution the positive-part symbols can thus be omitted. If $\Lin_{\lambda,u}$ is invertible, the local index is $\pm1$.

The next lemma constructs isolating domains that may vary with $\lambda$ and proves constancy of their degrees by homotopy invariance on fixed domains. All boundaries in a closed parameter strip are relative to that strip.

\begin{lemma}\label{Lem2.6}
Let $b<c$ and let $C$ be a nonempty union of connected components of
\[
 S_{[b,c]}=\Sb\cap([b,c]\times X)
\]
such that $\mathfrak{o}\notin C$. There exists a bounded relatively open set $\mathcal U\subset[b,c]\times \mathcal C^1_0$ whose closure consists of strongly positive functions and whose fixed points $(\lambda,u)$ of $\Phi_\lambda$ are exactly $C$, all lying in $\mathcal U$. In particular, its relative boundary contains no fixed points and
\[
 \deg_{\mathrm{LS}}(I-\Phi_\lambda,\mathcal U_\lambda)
\]
is independent of $\lambda\in[b,c]$, where $\mathcal U_\lambda=\{u:(\lambda,u)\in\mathcal U\}$ and an empty section has degree zero.
\end{lemma}
\begin{proof}
Writing $S=S_{[b,c]}$ and $Z=\mathcal C^1_0$, Proposition~\ref{Prop2.4} shows that $S$ is compact with finitely many components. Thus $C$ and $S\setminus C$ are disjoint compact sets in $[b,c]\times Z$. Let $\Pp_1\subset Z$ be the open cone of strongly positive functions. In the product metric, set
\[
 D=(S\setminus C)\cup\bigl([b,c]\times(Z\setminus\Pp_1)\bigr).
\]
The set $D$ is closed and disjoint from the compact set $C$, so $\dist(C,D)>0$. Choose $0<\varepsilon<\dist(C,D)/3$ and define
\[
 \mathcal U=\{z\in[b,c]\times Z:\dist(z,C)<\varepsilon\}.
\]
This set is bounded, and its relative closure lies in $[b,c]\times\Pp_1$. Since every nonzero fixed point of $\Phi_\lambda$ is a positive solution,
\[
 \{(\lambda,u)\in\overline{\mathcal U}:u=\Phi_\lambda u\}
 =C\subset\mathcal U.
\]
For the degree assertion, write $C_\lambda=\{u:(\lambda,u)\in C\}$ and fix $\lambda_0\in[b,c]$. If $C_{\lambda_0}$ is nonempty, its finiteness allows us to choose a finite union $V\subset Z$ of small open balls about its points. After shrinking the balls and choosing $\eta>0$, compactness of $C$ gives, for $|\lambda-\lambda_0|<\eta$ within the strip,
\[
 \overline V\subset\mathcal U_\lambda,\qquad
 C_\lambda\subset V,\qquad
 \{u\in\partial V:u=\Phi_\lambda u\}=\varnothing.
\]
Indeed, a sequence of section zeros outside $V$ with parameters tending to $\lambda_0$ would converge to a point of $C_{\lambda_0}\subset V$. Excision and ordinary homotopy invariance on the fixed domain $V$ now give
\[
 \deg_{\mathrm{LS}}(I-\Phi_\lambda,\mathcal U_\lambda)
 =\deg_{\mathrm{LS}}(I-\Phi_\lambda,V)
 =\deg_{\mathrm{LS}}(I-\Phi_{\lambda_0},V).
\]
If $C_{\lambda_0}$ is empty, compactness gives empty zero sets in all sufficiently nearby sections, whose degrees are zero. The section degree is therefore locally constant and hence constant on $[b,c]$, as in the moving-domain property of \cite[Th.~2.2]{LGS}.
\end{proof}

\section{Degree and proper ends}\label{Sec3}
For each nonempty index set $J$, the degree $d_J\ne0$ cannot be accounted for entirely by compact components. The cancellation lemma below is the compact-component degree balance appearing in global continuation, as in \cite[Th.~3.2]{LGS}, applied within one interval-norm class. Combining it with the locally finite graph structure will yield a proper ray on which $J(u)=J$.

For $c<\lambda_*$ and a nonempty subset $J\subset\{1,\ldots,n\}$, define
\[
 M_J(c)=\{(\lambda,u)\in\Sp:\lambda\le c,\ J(u)=J\}.
\]
Continuity of the interval norms and exclusion of the threshold $\rho$ make each $M_J(c)$ relatively open and closed in $\Sp\cap\{\lambda\le c\}$. Every connected subset below the cut therefore lies in one such class.

\begin{lemma}\label{Lem3.1}
Let $c<\lambda_*$, let $J$ be nonempty and let $C$ be a compact component of $M_J(c)$. Then
\begin{equation}\label{Eq3.1}
 \sum_{(c,u)\in C}\iota_c(u)=0.
\end{equation}
\end{lemma}
\begin{proof}
For $b<\min_C\pr$, the set $C$ is a component of $S_{[b,c]}$ because $J(u)$ is constant on every connected subset of that strip, so any such subset meeting $C$ lies in the same component of $M_J(c)$. Apply Lemma~\ref{Lem2.6} and write
\[
 D_C(\lambda)=\deg_{\mathrm{LS}}(I-\Phi_\lambda,\mathcal U_\lambda).
\]
Since the section at $b$ contains no zeros, constancy, excision and additivity give
\[
 0=D_C(b)=D_C(c)
   =\sum_{(c,u)\in C}\iota_c(u).
\]
The sum is finite by Lemma~\ref{Lem2.3} and remains valid at degenerate solutions because it uses local indices rather than intersection counts.
\end{proof}

\begin{theorem}[Ends distinguished by interval norms]\label{Th3.2}
For each nonempty $J\subset\{1,\ldots,n\}$ there exists a proper ray in $G$ on which $J(u)=J$ and $\lambda\to-\infty$. Rays with different values of $J(u)$ represent different ends, and $\norm u_\infty\to\infty$ along every such ray.
\end{theorem}
\begin{proof}
Choose a nonexceptional level $c<\lambda_*$, and let $C_1,\ldots,C_m$ be the components of $M_J(c)$ meeting $\pr=c$. Finiteness of the fiber gives
\[
 m\le\#\{(c,u):J(u)=J\}<\infty.
\]
If all $C_\ell$ were compact, additivity and Lemma~\ref{Lem3.1} would imply
\[
 0\ne d_J
 =\sum_{J(u)=J}\iota_c(u)
 =\sum_{\ell=1}^m\ \sum_{(c,u)\in C_\ell}\iota_c(u)=0,
\]
a contradiction. Hence some $H=C_\ell$ is noncompact. Since every solution at the chosen level lies on a local arc, any component meeting that level contains nonisolated solutions. The local graph description and the closedness of $M_J(c)$ make $H$ a closed connected subgraph of $G$, after adding the finitely many cut points at $\pr=c$ as vertices.

Use the vertices and edges inherited from $G$ and choose $v_0\in V(H)$. Form an auxiliary rooted tree whose vertices are the finite simple edge paths in $H$ starting at $v_0$. Its root is the path of length zero. A path
\[
 P=(v_0,e_1,v_1,\ldots,e_m,v_m)
\]
is joined to each path obtained by adding one edge and one new vertex. Every nonroot path has a unique predecessor, obtained by deleting its last edge and vertex, so the auxiliary graph is a tree with finitely many successors at each vertex. Its path lengths are unbounded, since otherwise the tree would be finite and connectedness would make $H$ itself finite and compact. Starting at the root, we can therefore choose successive extensions that admit arbitrarily long further extensions. This gives a sequence $\{v_j\}_{j\ge0}$ of distinct vertices joined by edges, and hence an embedding $r:[0,\infty)\to H$ with $r(j)=v_j$.

Every compact subset of $G$ meets finitely many vertices and edges, and the ray traverses each edge only once, so it is proper. By~\eqref{Eq2.6},
\[
 \lim_{t\to\infty}\pr(r(t))=-\infty.
\]
Its image lies in $M_J(c)$, so $J(u)=J$ along the ray.

The compact set $G_{\ge c}$ separates solutions with different index sets through the decomposition
\[
 G\setminus G_{\ge c}
 =\bigsqcup_{J\ne\varnothing}
   \bigl(G\cap\{\pr<c,\ J(u)=J\}\bigr),
\]
where each term is relatively open and closed. Rays with different values of $J(u)$ consequently have tails in different components of this complement, so they cannot represent the same end.

Finally, at a maximum point $x_u$ and for $\lambda<0$, the equation gives $0\le\lambda\norm u_\infty+a(x_u)\norm u_\infty^p$. Consequently,
\[
 \norm u_\infty^{p-1}\ge\frac{-\lambda}{\max a^+},
 \qquad \lim_{\lambda\to-\infty}\frac{-\lambda}{\max a^+}=+\infty.
\]
Thus the uniform norm diverges along every negative end, completing the proof.
\end{proof}

This proves Theorem~\ref{Th1.1}.

\section{Components, branching and cohomology}\label{Sec4}
The $2^n-1$ ends obtained above impose a lower bound on the combined number of components and branching points. Their labels yield a second set of constraints through compactly supported cohomology, which we develop below using functions on the vertices and edges.

\subsection{Counting components and branching}\label{Subsec4.1}
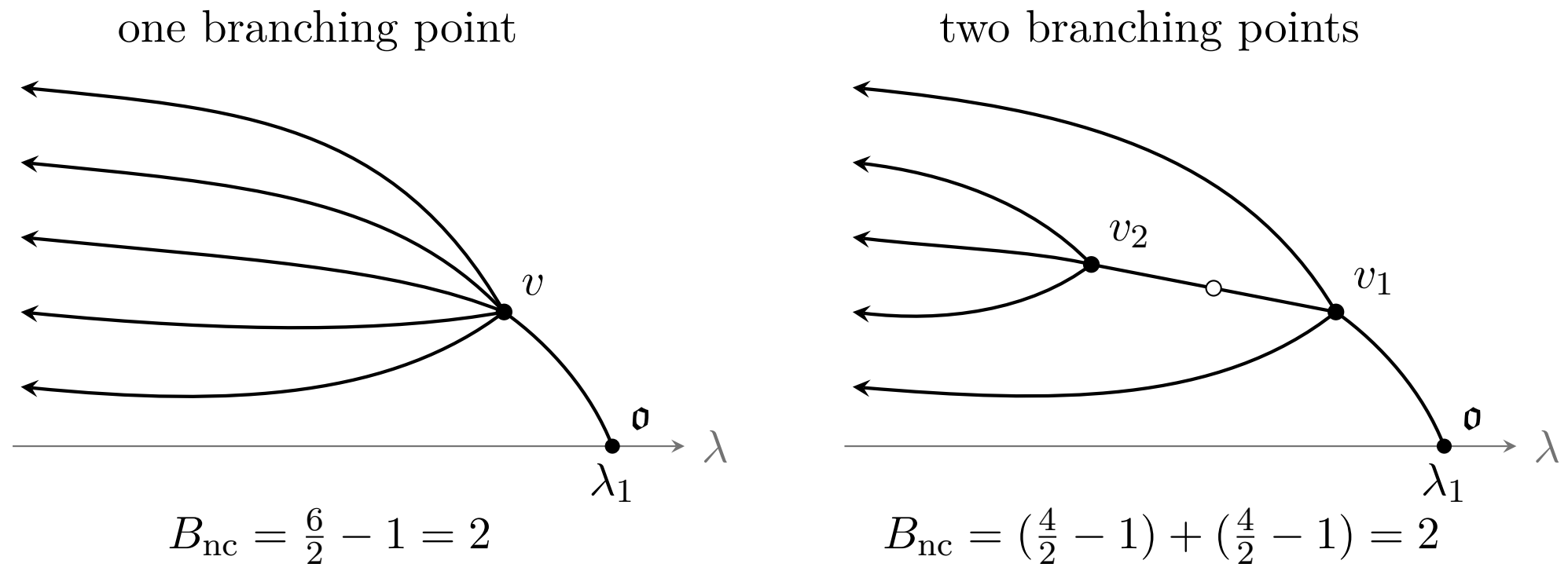
\begin{figure}[!htbp]
\centering
\begin{tikzpicture}[x=.98cm,y=1cm,font=\small,>=stealth]
\foreach \shift/\panel in {0/a,7/b}{
\begin{scope}[xshift=\shift cm]
\node at (2.65,3.5) {\if a\panel one branching point\else two branching points\fi};
\draw[->,black!55] (.03,0)--(5.8,0) node[right] {$\lambda$};
\draw[thick] (5.18,0) .. controls (5.0,.45) and (4.65,.85) .. (4.25,1.13);
\fill (5.18,0) circle (1.8pt) node[above right] {$\mathfrak{o}$};
\node[below] at (5.18,0) {$\lambda_1$};
\fill (4.25,1.13) circle (2pt);
\if a\panel
\foreach \y/\bend in {3.02/2.75,2.39/2.12,1.76/1.51,1.13/.96,.50/.42}{
\draw[thick,->] (4.25,1.13) .. controls (3.30,\bend) and (1.95,\y-.18) .. (.10,\y);
}
\node[above right] at (4.25,1.13) {$v$};
\node at (2.75,-.76) {$B_{\mathrm{nc}}=\frac{6}{2}-1=2$};
\else
\draw[thick,->] (4.25,1.13) .. controls (3.42,2.50) and (1.95,2.84) .. (.1,3.02);
\draw[thick,->] (4.25,1.13) .. controls (3.24,.34) and (1.80,.36) .. (.1,.50);
\draw[thick] (4.25,1.13)--(2.15,1.53);
\draw[fill=white] (3.20,1.33) circle (1.8pt);
\foreach \y/\bend in {2.39/2.06,1.76/1.65,1.13/1.14}{
\draw[thick,->] (2.15,1.53) .. controls (1.63,\bend) and (.88,\y-.09) .. (.10,\y);
}
\fill (2.15,1.53) circle (2pt);
\node[above right] at (4.25,1.13) {$v_1$};
\node[above right] at (2.15,1.53) {$v_2$};
\node at (2.75,-.76) {$B_{\mathrm{nc}}=(\frac{4}{2}-1)+(\frac{4}{2}-1)=2$};
\fi
\end{scope}}
\end{tikzpicture}
\caption{Two ways for the principal component to have five negative ends. The vertex $v$ has valence six, while $v_1$ and $v_2$ each have valence four. The open circle marks a point of valence two. Both graphs have total branching contribution two.}\label{Fig6}
\end{figure}

Let $G_{\mathrm{nc}}$ be the union of the $K$ noncompact components of $G$, and let $e_-(G)$ count their negative ends. Since every component of a locally finite graph is open, this union is open and closed in $G$. The branching points form the set $\mathcal B$ defined above, and the valence calculation gives
\[
 \mathcal B=\{v\in G:\deg_G(v)\ge4\}.
\]
At such a point, the valence is $2m$ with $m\ge2$, and $m-1$ counts branching beyond the two half-arcs of an ordinary curve. The total branching contribution on the noncompact components is therefore
\begin{equation}\label{Eq4.1}
 B_{\mathrm{nc}}=\sum_{v\in\mathcal B\cap G_{\mathrm{nc}}}
 \left(\frac{\deg_G(v)}2-1\right).
\end{equation}
An added valence-two vertex contributes $2/2-1=0$, making this sum invariant under subdivision. A vertex of valence four contributes one and a vertex of valence six contributes two, so the two arrangements in Figure~\ref{Fig6} have the same total branching contribution.

For a graph $H$, let $b_0(H)$ be its number of connected components and $b_1(H)$ its \emph{first Betti number}, the number of independent cycles. Here cycles are finite edge sets with an even number of incidences at every vertex, counting a loop twice. They form a vector space over $\Ftwo$ under addition that cancels shared edges, and $b_1(H)$ is its dimension.

To compute $b_1$ for a finite connected graph, retain all vertices and successively delete edges belonging to simple cycles. Each deletion preserves connectedness because the rest of the cycle still joins the endpoints of the removed edge. When no cycle remains, the resulting graph is a tree containing every original vertex, called a \emph{spanning tree}.

\begin{figure}[!htbp]
\centering
\begin{tikzpicture}[x=1cm,y=1cm,font=\small]
\begin{scope}
\node at (1.2,1.5) {tree};
\draw[line width=.9pt] (0,0)--(1.2,.5)--(2.4,0);
\fill (0,0) circle (2pt);
\fill (1.2,.5) circle (2pt);
\fill (2.4,0) circle (2pt);
\node at (1.2,-1.5) {$b_1=0$};
\end{scope}
\begin{scope}[xshift=4.3cm]
\node at (1.2,1.5) {circle};
\draw[dashed,line width=.9pt] (0,0) .. controls (.4,1.15) and (2,1.15) .. (2.4,0);
\draw[line width=.9pt] (0,0) .. controls (.4,-1.15) and (2,-1.15) .. (2.4,0);
\fill (0,0) circle (2pt);
\fill (2.4,0) circle (2pt);
\node at (1.2,-1.5) {$b_1=1$};
\end{scope}
\begin{scope}[xshift=8.6cm]
\node at (1.2,1.5) {three edges};
\draw[dashed,line width=.9pt] (0,0) .. controls (.4,1.15) and (2,1.15) .. (2.4,0);
\draw[line width=.9pt] (0,0)--(2.4,0);
\draw[dashed,line width=.9pt] (0,0) .. controls (.4,-1.15) and (2,-1.15) .. (2.4,0);
\fill (0,0) circle (2pt);
\fill (2.4,0) circle (2pt);
\node[above] at (1.2,.86) {$e_1$};
\node[above] at (1.2,0) {$e_2$};
\node[below] at (1.2,-.86) {$e_3$};
\node at (1.2,-1.5) {$b_1=2$};
\end{scope}
\end{tikzpicture}
\caption{The solid edges form a spanning tree in each graph. Each dashed edge closes one independent cycle when added to that tree. In the right-hand graph, any two of the three edges form a simple closed path, but only two of these cycles are independent.}\label{Fig7}
\end{figure}
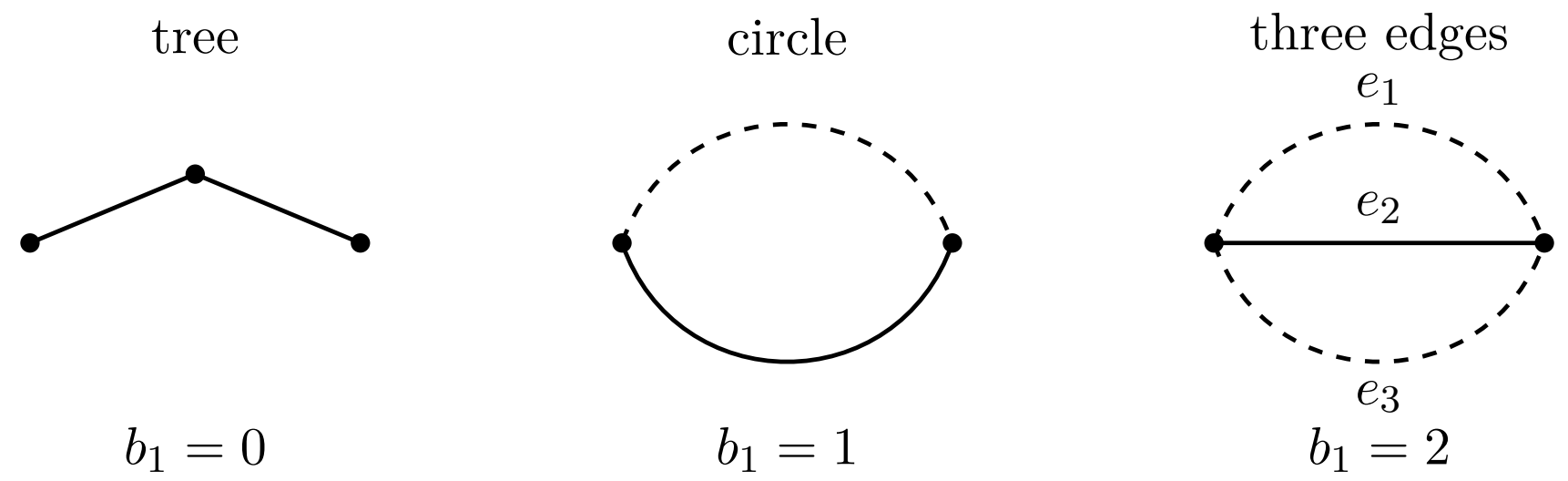

Since a finite tree with $v$ vertices has $v-1$ edges, the spanning tree contains $|V(H)|-1$ edges, leaving
\[
 |E(H)|-(|V(H)|-1)
\]
deleted edges. Restoring one of them creates a cycle with the unique tree path between its endpoints. Each of these cycles contains its own additional edge, which proves independence, and together they generate every cycle by cancellation of the additional edges. Their number is therefore $b_1(H)$. For a graph with several components, the union of one spanning tree in each component is a \emph{spanning forest} with $|V(H)|-b_0(H)$ edges, giving
\begin{equation}\label{Eq4.2}
 \begin{aligned}
 b_1(H)&=|E(H)|-|V(H)|+b_0(H),\\
 \chi(H)&=|V(H)|-|E(H)|=b_0(H)-b_1(H),
 \end{aligned}
\end{equation}
where $\chi(H)$ is the Euler characteristic. Subdivision adds one vertex and one edge, leaving both expressions in~\eqref{Eq4.2} unchanged. Figure~\ref{Fig7} illustrates the spanning-tree count in three connected graphs.

The same definition of $b_1$ applies to $G_{\mathrm{nc}}$, although the quantities $K$, $B_{\mathrm{nc}}$ and $b_1$ may now be infinite. When $K+B_{\mathrm{nc}}<\infty$, the next lemma reduces $G_{\mathrm{nc}}$ to a finite graph by adjoining one point at each end, allowing us to use~\eqref{Eq4.2}.

\begin{lemma}\label{Lem4.1}
If $K+B_{\mathrm{nc}}<\infty$, then $G_{\mathrm{nc}}$ is homeomorphic to a finite graph with finitely many leaves removed, each representing a negative end.
\end{lemma}
\begin{proof}
There are finitely many noncompact components, and each branching point contributes at least one to $B_{\mathrm{nc}}$. Hence the set
\[
 S:=\{\mathfrak{o}\}\cup(\mathcal B\cap G_{\mathrm{nc}}),
 \qquad |S|\le1+B_{\mathrm{nc}},
\]
is finite, while every point outside $S$ has valence two and a neighborhood homeomorphic to an open interval.

For each $v\in S$, take a short initial segment of each incident half-arc. Their union with $v$ is an open neighborhood $U_v$, and we choose these neighborhoods with pairwise disjoint compact closures. Each outer endpoint is an ordinary valence-two point. Set
\[
 M:=G_{\mathrm{nc}}\setminus\bigcup_{v\in S}U_v,
 \qquad A:=\bigcup_{v\in S}\partial U_v,
 \qquad |A|=\sum_{v\in S}\deg_G(v)<\infty.
\]
Points of $M\setminus A$ have interval neighborhoods, while points of $A$ have half-interval neighborhoods. These give $M$ the structure of a one-dimensional manifold with boundary $A$, whose local coordinates take values in an open interval or in $[0,\varepsilon)$ at a boundary point. In Figure~\ref{Fig8}, the solid curves form $M$ and the filled points mark its boundary $A$.

The components of $M$ meeting $A$ number at most $|A|$. A component of $M$ disjoint from $A$ is an entire component of $G_{\mathrm{nc}}$, since a path to a removed neighborhood must first cross $A$. There are at most $K$ components of this second kind, and therefore
\[
 b_0(M)\le |A|+K<\infty.
\]
Choose a positive length for each edge and use distance along the edges as a local coordinate, with either orientation. Changes of coordinate have the form $t=\pm s+\beta$, giving $M$ a smooth structure to which the classification of connected one-dimensional manifolds~\cite[App., Th., p.~55]{Milnor} applies. Each component is a compact interval, a half-line, a line or a circle. A circle has no boundary and would be an entire compact component of $G_{\mathrm{nc}}$, which is impossible. The remaining possibilities are therefore
\[
 [0,1],\qquad [0,\infty),\qquad\R.
\]

\begin{figure}[!htbp]
\centering
\begin{tikzpicture}[x=1cm,y=.79351cm,font=\small,>=stealth,
    retained/.style={line width=.95pt},
    removed/.style={densely dashed,line width=.75pt},
    axis/.style={->,line width=.45pt,black!60}]
\draw[axis] (0,0)--(11.58,0) node[below left,text=black] {$\lambda$};
\draw[axis] (0,0)--(0,3.98) node[left,text=black] {$\norm{u}_\infty$};
\coordinate (o)  at (10.38,0);
\coordinate (ao) at (9.85,.46);
\coordinate (v)  at (6.93,1.40);
\coordinate (se) at (7.43,1.10);
\coordinate (sw) at (6.43,1.10);
\coordinate (nw) at (6.43,1.70);
\coordinate (ne) at (7.43,1.70);
\draw[removed] (o) .. controls (10.22,.18) and (10.05,.33) .. (ao);
\draw[retained] (ao) .. controls (9.00,1.00) and (8.23,.62) .. (se);
\foreach \tip in {se,sw,nw,ne} {\draw[removed] (v)--(\tip);}
\draw[retained] (nw) .. controls (5.83,2.06) and (5.52,2.40) .. (5.65,2.83)
  .. controls (5.87,3.58) and (7.78,3.75) .. (8.12,3.04)
  .. controls (8.35,2.56) and (8.08,2.09) .. (ne);
\draw[retained,->] (sw) .. controls (5.63,.62) and (4.68,.66) .. (3.68,.85)
  .. controls (2.48,1.08) and (1.27,1.40) .. (.38,1.65);
\draw[retained,<->] (.38,3.53) .. controls (1.44,3.21) and (2.50,3.05) .. (3.34,2.83)
  .. controls (4.77,2.46) and (4.57,1.88) .. (3.40,1.88)
  .. controls (2.40,1.88) and (1.29,2.18) .. (.38,2.43);
\foreach \point in {ao,se,sw,nw,ne} {\fill (\point) circle (2pt);}
\node[above right,inner sep=2pt] at (o) {$\mathfrak{o}$};
\node[below,inner sep=2pt] at (o) {$\lambda_1$};
\node[below=2pt,inner sep=2.5pt] at (v) {$v$};
\node at (6.93,2.56) {$[0,1]$};
\node at (8.68,1.19) {$[0,1]$};
\node at (4.43,.32) {$[0,\infty)$};
\node at (2.13,2.53) {$\R$};
\draw[retained] (8.94,3.39)--(9.48,3.39) node[right=2pt] {$M$};
\draw[removed] (8.94,2.85)--(9.48,2.85) node[right=2pt] {$U_{\mathfrak{o}},\ U_v$};
\fill (9.21,2.31) circle (2pt);
\node[right=2pt] at (9.48,2.31) {$A=\partial M$};
\end{tikzpicture}
\caption{The solid curves form $M$ after the dashed neighborhoods of $\mathfrak{o}$ and $v$ are removed. The five filled dots belong to $M$ and form its boundary $A$, one on $\partial U_{\mathfrak{o}}$ and four on $\partial U_v$.}\label{Fig8}
\end{figure}
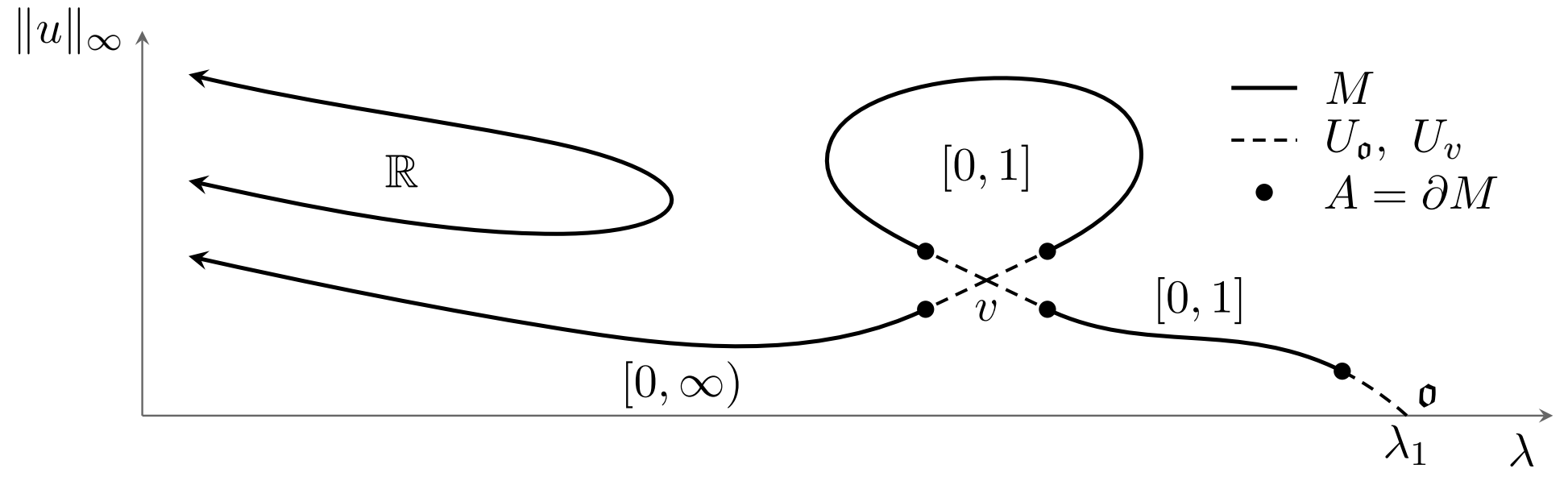

The half-line and line become compact intervals after adjoining one or two endpoints at infinity, respectively.

Compactifying these pieces and gluing them back to the closures of the $U_v$ along $A$ produces a finite graph $\widehat G_{\mathrm{nc}}$. Each added point has valence one, and writing $D_\infty$ for their set gives
\[
 G_{\mathrm{nc}}\cong\widehat G_{\mathrm{nc}}\setminus D_\infty.
\]
To identify $D_\infty$ with the ends, retain the neighborhoods $\overline U_v$, all compact interval pieces and a bounded central part of every remaining piece. Their union is a finite compact core. Outside it there are finitely many disjoint half-lines, one for each point of $D_\infty$. Their tails remain separated outside this core, and every proper ray eventually follows one of them. They therefore represent exactly the ends of $G_{\mathrm{nc}}$, all negative by~\eqref{Eq2.6}. Consequently,
\[
 |D_\infty|=e_-(G).
\]
Adjoining a leaf preserves the number of connected components. Every cycle lies in the finite core, since neither an exterior half-line nor an edge ending at a leaf can belong to a cycle. Thus the compactified graph has the same $K$ components and the same first Betti number as $G_{\mathrm{nc}}$.
\end{proof}

\begin{theorem}[Components and branching]\label{Th4.2}
One always has $K+B_{\mathrm{nc}}\ge2^{n-1}$. If $K+B_{\mathrm{nc}}<\infty$, then
\begin{equation}\label{Eq4.3}
 K+B_{\mathrm{nc}}-b_1(G_{\mathrm{nc}})
 =\frac{e_-(G)+1}{2}\ge2^{n-1}.
\end{equation}
Thus a connected graph satisfies $B_{\mathrm{nc}}\ge2^{n-1}-1$, strengthened to $B_{\mathrm{nc}}\ge2^{n-1}-1+b_1(G)$ when its branching is finite.
\end{theorem}
\begin{proof}
The inequality is immediate if $K+B_{\mathrm{nc}}=\infty$. Otherwise the preceding lemma gives a finite graph $\widehat G_{\mathrm{nc}}$ with $K$ components and the same first Betti number as $G_{\mathrm{nc}}$. Its leaves consist of $\mathfrak{o}$ and one added point for each negative end, for a total of $1+e_-(G)$.

We count the quantity $\deg(v)-2$ over its vertices in two ways. A leaf contributes $-1$, a valence-two vertex contributes zero, and a branching vertex of valence $2m$ contributes $2(m-1)$. The latter is twice its contribution to $B_{\mathrm{nc}}$. Thus, writing $\widehat G=\widehat G_{\mathrm{nc}}$ within this proof,
\[
 \begin{aligned}
 \sum_{v\in V(\widehat G)}(\deg_{\widehat G}(v)-2)
 &=\sum_{\deg(v)=1}(-1)
   +\sum_{\deg(v)=2}0
   +\sum_{\deg(v)\ge4}2\left(\frac{\deg(v)}2-1\right)\\
 &=-(1+e_-(G))+2B_{\mathrm{nc}}.
 \end{aligned}
\]
The added endpoints do not change the old branching valences, so the last sum is exactly the one in~\eqref{Eq4.1}.

On the other hand, every edge has two endpoint incidences, including a loop counted twice at its vertex. Therefore
\[
 \sum_{v\in V(\widehat G)}\deg_{\widehat G}(v)=2|E(\widehat G)|.
\]
Subtracting two for every vertex and applying~\eqref{Eq4.2} gives
\[
 \begin{aligned}
 \sum_{v\in V(\widehat G)}(\deg_{\widehat G}(v)-2)
 &=2|E(\widehat G)|-2|V(\widehat G)|\\
 &=-2\chi(\widehat G)
 =-2K+2b_1(G_{\mathrm{nc}}).
 \end{aligned}
\]
Equating the two counts and rearranging yields
\[
 K+B_{\mathrm{nc}}-b_1(G_{\mathrm{nc}})
 =\frac{1+e_-(G)}2
 \ge\frac{1+(2^n-1)}2=2^{n-1},
\]
where the last inequality uses the end bound. Since $b_1(G_{\mathrm{nc}})\ge0$, this also gives $K+B_{\mathrm{nc}}\ge2^{n-1}$. In the connected case the principal component is noncompact, so substituting $K=1$ and $G=G_{\mathrm{nc}}$ yields the stated bounds. The extra $1$ in the identity accounts for the finite leaf $\mathfrak{o}$.
\end{proof}

\subsection{Interval labels and cohomology}\label{Subsec4.2}
To use the interval labels $J(u)$, we encode them by vertex functions and study their differences across edges. Compactly supported cohomology records which finitely supported edge functions arise in this way from finitely supported vertex functions. We develop this description in detail for readers unfamiliar with cohomology, using coefficients in $\Ftwo=\{0,1\}$ throughout and usually omitting them from the notation.

The chosen vertices $V$ and edges $E$ form a \emph{cellular structure} on $G$, with vertices as zero-dimensional cells and open edges as one-dimensional cells. Assign to each vertex a value $f(v)\in\Ftwo$ and, for an edge $e$ with endpoints $v_0,v_1$, define
\[
 (\delta f)(e)=f(v_0)+f(v_1).
\]
Since addition is modulo two, $\delta f$ is one exactly when the endpoint values differ. For a loop the two contributions cancel, and orientations of the edges are unnecessary.

Vertex functions are called \emph{$0$-cochains}, edge functions are \emph{$1$-cochains}, and $\delta$ is the \emph{coboundary} operator. Compact support means that only finitely many values are nonzero, because compact subsets of $G$ meet finitely many vertices and edges. In notation,
\[
 \begin{aligned}
 C^0(G)&=\prod_{v\in V}\Ftwo,&
 C^1(G)&=\prod_{e\in E}\Ftwo,\\
 C_c^0(G)&=\bigoplus_{v\in V}\Ftwo,&
 C_c^1(G)&=\bigoplus_{e\in E}\Ftwo.
 \end{aligned}
\]
The products allow infinitely many nonzero values, whereas the direct sums require finite support, denoted by the subscript $c$. The coboundary is defined on all vertex functions, $\delta:C^0(G)\to C^1(G)$, and maps $C_c^0(G)$ into $C_c^1(G)$ because only finitely many edges meet each vertex.

The equation $\delta f=0$ requires equal values at the endpoints of every edge. This equality propagates along finite edge paths and hence throughout each connected component, giving
\begin{equation}\label{Eq4.4}
 \Null{\delta}=\{f\in C^0(G):f|_{V\cap C}\text{ is constant for every component }C\}.
\end{equation}
The constants may differ between components. For a finitely supported $f$, the constant on every noncompact component is zero because that component has infinitely many vertices.

The \emph{compactly supported cohomology} of the graph in degree one is
\begin{equation}\label{Eq4.5}
 H_c^1(G)=C_c^1(G)/\delta C_c^0(G).
\end{equation}
For graphs this is the cellular description of compactly supported cohomology~\cite[Def.~A.1, Rem.~A.3]{MW}, with no further condition because there are no cells of dimension two. Concretely,
\[
 [\eta]=[\eta']\quad\Longleftrightarrow\quad
 \eta-\eta'=\delta b\text{ for some }b\in C_c^0(G).
\]
Thus $[\eta]=0$ precisely when $\eta$ is the coboundary of a finitely supported vertex function, and the dimension of $H_c^1(G)$ counts independent obstructions to such a representation. The following subdivision argument shows that this space is intrinsic to the graph.

Let $\widetilde G$ denote the same graph with a finer locally finite cellular structure. Each original edge $e$ is divided into finitely many edges $e_1,\ldots,e_m$, since its closure is compact. Define
\[
 S:C_c^1(\widetilde G)\longrightarrow C_c^1(G),\qquad
 (S\eta)(e):=\sum_{i=1}^m\eta(e_i).
\]
For a vertex function $f$ with finitely supported coboundary, let $v_0,\ldots,v_m$ be the successive vertices along $e$. The intermediate values cancel in pairs,
\[
 (S\delta f)(e)
 =\sum_{i=1}^m\bigl(f(v_{i-1})+f(v_i)\bigr)
 =f(v_0)+f(v_m)
 =(\delta(f|_V))(e).
\]
Thus $S$ takes finite coboundaries to finite coboundaries. It is surjective, since an edge value can be placed on any one of its subedges. Its kernel consists of finitely supported edge functions whose sum along each original edge is zero. Along each such edge, successive partial sums give vertex values that vanish at both original endpoints and have the prescribed coboundary. Only finitely many edges are involved, so these values define a finitely supported function on the new vertices. Finally, every finitely supported function on $V$ extends by zero to the new vertices. Together these observations show that $S$ induces an isomorphism
\[
 \overline S:H_c^1(\widetilde G)\longrightarrow H_c^1(G),\qquad
 \overline S([\eta]):=[S\eta].
\]
Adding edge values is compatible with further subdivision. Applying this construction to a common refinement therefore identifies the cohomology spaces of any two locally finite cellular structures on $G$.

For a concrete example, represent a line by vertices $\{v_j\}_{j\in\mathbb Z}$ and edges $\{e_j\}_{j\in\mathbb Z}$, where $e_j$ joins $v_j$ to $v_{j+1}$. Let $\eta_k$ equal one on $e_k$ and zero on every other edge. The vertex function
\[
 f(v_j):=\begin{cases}0,&j\le0,\\1,&j\ge1,\end{cases}
 \qquad \delta f=\eta_0,
\]
has infinite support, although its coboundary is supported on one edge. No finitely supported vertex function $b$ can have this coboundary, because
\[
 \sum_{j\in\mathbb Z}(\delta b)(e_j)
 =\sum_{j\in\mathbb Z}\bigl(b(v_j)+b(v_{j+1})\bigr)=0,
 \qquad \sum_{j\in\mathbb Z}\eta_0(e_j)=1.
\]
Each vertex coefficient occurs twice in the first sum, while $\eta_k-\eta_0$ is the coboundary of the characteristic function of the vertices between the two edges. All one-edge cochains therefore represent the same nonzero class, yielding $H_c^1(\R,\Ftwo)\simeq\Ftwo$ for this graph with two ends and no cycles.

On a half-line the vertices and edges are instead indexed by $j\ge0$. For each $k\ge0$, the finitely supported function
\[
 b_k(v_j):=\begin{cases}1,&0\le j\le k,\\0,&j>k,\end{cases}
 \qquad \delta b_k=\eta_k,
\]
shows that every one-edge cochain, and hence every finitely supported cochain, is a coboundary. Thus $H_c^1([0,\infty),\Ftwo)=0$. The finite endpoint allows $b_k$ to vanish beyond the marked edge, explaining the different role of the principal component below.

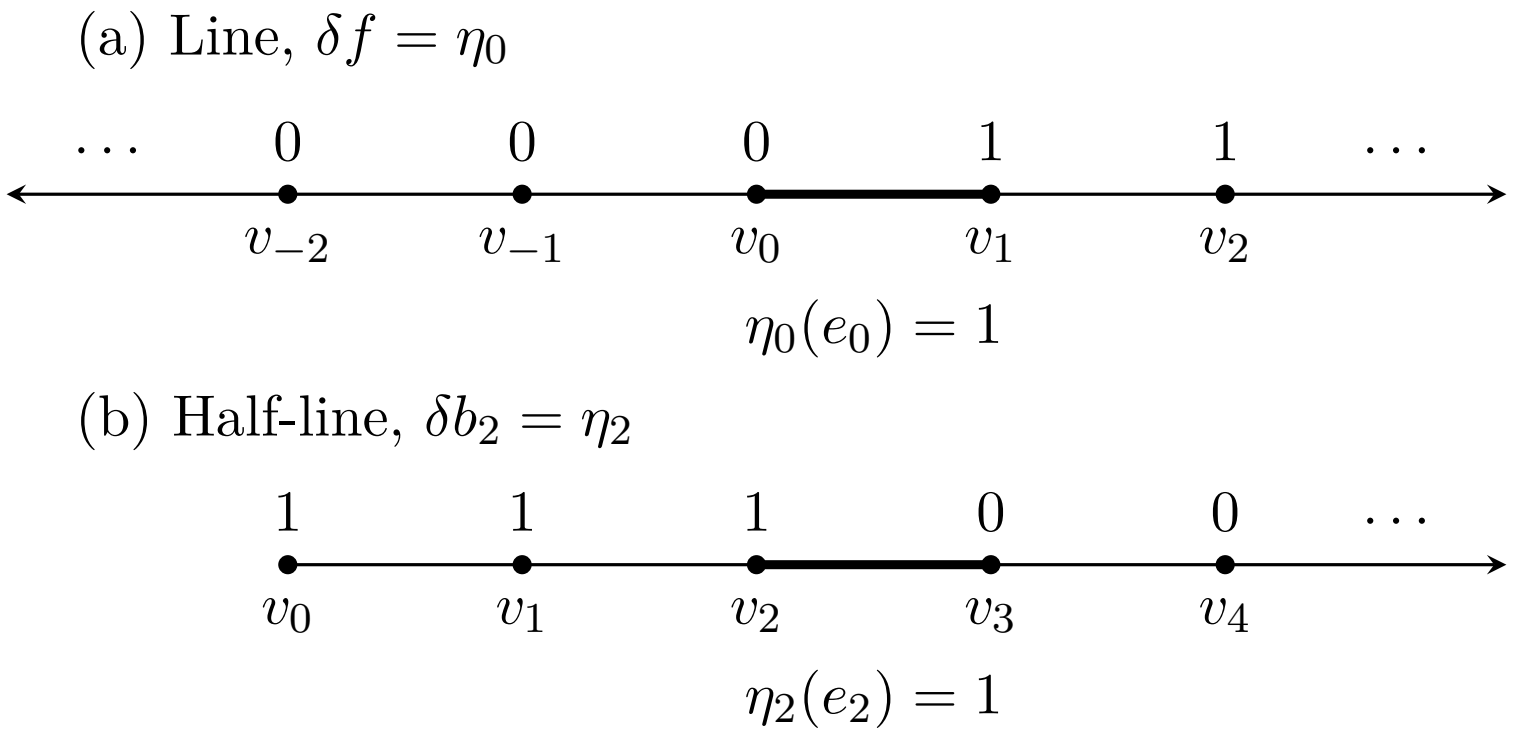
\begin{figure}[!htbp]
\centering
\begin{tikzpicture}[x=1.55cm,y=1cm,font=\small,>=stealth]
\begin{scope}
\node[anchor=west] at (-3,1.02) {(a) Line, $\delta f=\eta_0$};
\draw[<->,line width=.6pt] (-3.2,0)--(3.2,0);
\draw[line width=1.7pt] (0,0)--(1,0);
\foreach \j/\val in {-2/0,-1/0,0/0,1/1,2/1}{
 \fill (\j,0) circle (1.8pt);
 \node[above] at (\j,.08) {$\val$};
 \node[below] at (\j,-.08) {$v_{\j}$};
}
\node[below] at (.5,-.55) {$\eta_0(e_0)=1$};
\node at (-2.75,.28) {$\cdots$};\node at (2.75,.28) {$\cdots$};
\end{scope}
\begin{scope}[yshift=-2.45cm]
\node[anchor=west] at (-3,.95) {(b) Half-line, $\delta b_2=\eta_2$};
\draw[->,line width=.6pt] (-2,0)--(3.2,0);
\draw[line width=1.7pt] (0,0)--(1,0);
\foreach \x/\j/\val in {-2/0/1,-1/1/1,0/2/1,1/3/0,2/4/0}{
 \fill (\x,0) circle (1.8pt);
 \node[above] at (\x,.08) {$\val$};
 \node[below] at (\x,-.08) {$v_{\j}$};
}
\node[below] at (.5,-.55) {$\eta_2(e_2)=1$};
\node at (2.75,.28) {$\cdots$};
\end{scope}
\end{tikzpicture}
\caption{Vertex values are written above the dots, and vertex names below them. The bold edge has coboundary value one, while every other edge has value zero. On the line, $f$ remains one to the right and has infinite support. On the half-line, $b_2$ is supported on just the first three vertices. This difference gives a nonzero class in (a) and the zero class in (b).}\label{Fig9}
\end{figure}

\begin{figure}[!htbp]
\centering
\begin{tikzpicture}[x=1.5cm,y=1cm,font=\small,>=stealth]
\draw[<->,line width=.6pt] (-3.1,0)--(3.1,0);
\draw[line width=1.6pt] (0,0)--(1,0);
\draw[dashed,line width=1.2pt] (0,0) .. controls (-1,1.9) and (1,1.9) .. (0,0);
\foreach \j/\val in {-2/0,-1/0,0/0,1/1,2/1}{
 \fill (\j,0) circle (1.8pt);
 \node[below] at (\j,-.08) {$v_{\j}$};
 \node[above] at (\j,.08) {$\val$};
}
\node[right] at (.45,1.15) {$\eta_\ell(\ell)=1$};
\node[below] at (.5,-.58) {$\eta_0(e_0)=1$};
\node at (-2.65,.25) {$\cdots$};\node at (2.65,.25) {$\cdots$};
\end{tikzpicture}
\caption{A line with a loop. The solid bold edge supports $\eta_0$ and the dashed loop supports $\eta_\ell$. The values at the vertices define the same step function $f$ as on the line, with $\delta f=\eta_0$ and zero coboundary on the loop. The two cochains represent independent classes.}\label{Fig10}
\end{figure}
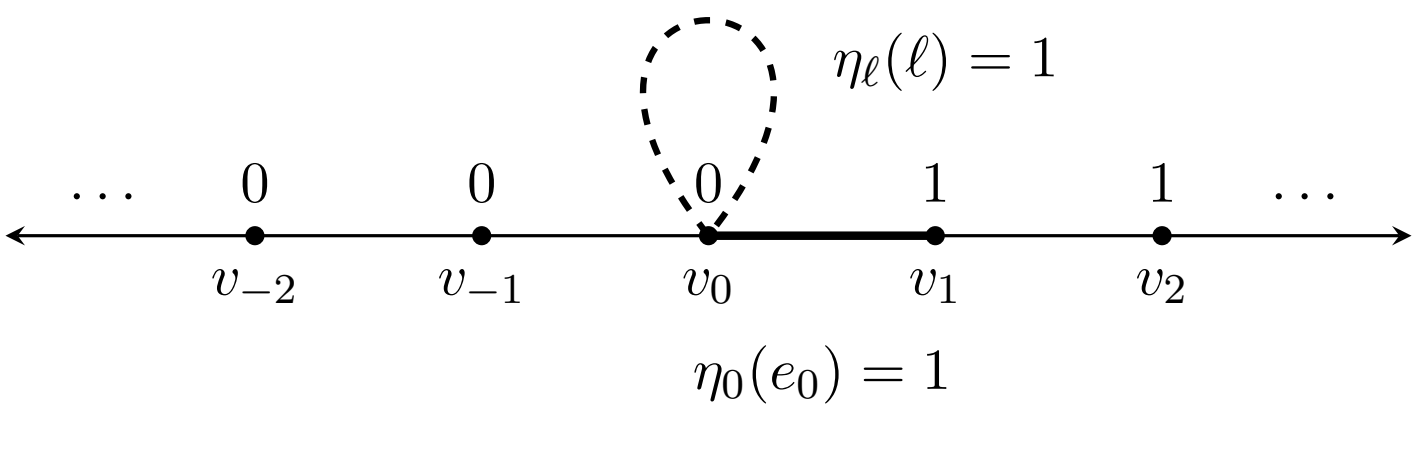

The next example separates the contribution of a cycle from that of the ends. Attach one loop $\ell$ at $v_0$ of the preceding line, and denote the resulting graph by $H$. It has two ends and one independent cycle. Let $\eta_\ell$ be one on the loop and zero on the line, and retain $\eta_0$, supported on the edge from $v_0$ to $v_1$. Figure~\ref{Fig10} marks the support of each cochain.

For every vertex function $b$, the two ends of the loop coincide, so
\[
 (\delta b)(\ell)=b(v_0)+b(v_0)=0.
\]
The cochain $\eta_\ell$ detects the loop and represents a nonzero class even when arbitrary vertex functions are allowed. The class of $\eta_0=\delta f$, by contrast, depends on the finite-support requirement. To prove their independence, suppose that $a\eta_0+d\eta_\ell=\delta b$ with $b$ finitely supported. Evaluation on the loop gives $d=0$, and summation over the line edges gives $a=0$ because each vertex value occurs twice. As in the preceding example, every cochain supported on the line is equivalent to a multiple of $\eta_0$, so
\[
 H_c^1(H,\Ftwo)=\operatorname{span}\{[\eta_0],[\eta_\ell]\}
 \simeq\Ftwo^2.
\]
One class comes from the two directions towards infinity and the other from the loop. Compactly supported cohomology can therefore detect both features, whereas the classes constructed from $J(u)$ below specifically record the labels at infinity.

To apply this construction to the interval labels, choose a nonexceptional level $c<\lambda_*$. Subdividing at the finitely many solutions with $\lambda=c$ makes each cut point a vertex, with every open edge lying on one side of the cut. For a vertex $v=(\lambda,u)$ below the cut, write $J(v):=J(u)$. For each nonempty $J\subset\{1,\ldots,n\}$, define
\[
 f_{J,c}:V\to\Ftwo,\qquad f_{J,c}(v):=
 \begin{cases}
 1,&\pr(v)<c\text{ and }J(v)=J,\\
 0,&\text{otherwise}.
 \end{cases}
\]
The function marks the vertices below the cut carrying label $J$, assigning zero to all other vertices, including those on the cut. At vertices sufficiently far along a ray labeled by $I$, its value is therefore one if $I=J$ and zero otherwise.

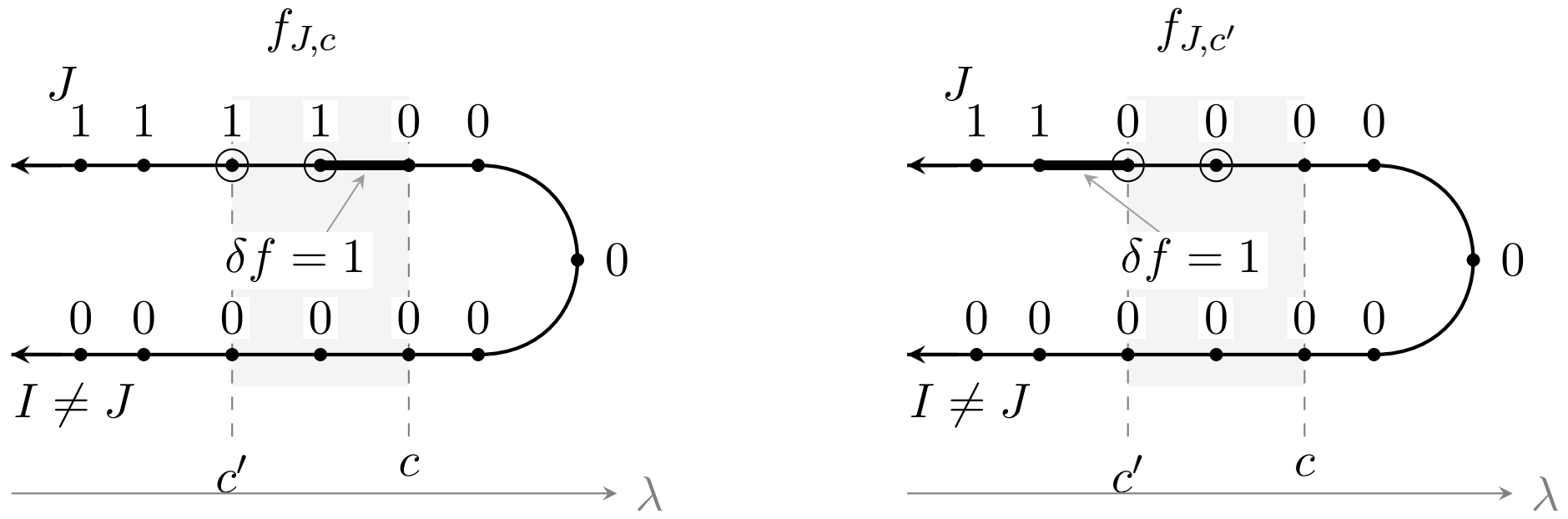
\begin{figure}[!htbp]
\centering
\begin{tikzpicture}[x=1cm,y=1cm,>=stealth,font=\small]
\foreach \panel/\shift in {0/0,1/7.1}{
\begin{scope}[xshift=\shift cm]
\node at (2.35,2.85)
 {\ifnum\panel=0 $f_{J,c}$\else $f_{J,c'}$\fi};
\fill[gray!9] (1.8,.05) rectangle (3.2,2.35);
\draw[dashed,gray] (1.8,-.35)--(1.8,2.35);
\draw[dashed,gray] (3.2,-.35)--(3.2,2.35);
\node[below] at (1.8,-.35) {$c'$};
\node[below] at (3.2,-.35) {$c$};
\draw[->,gray] (.05,-.8)--(4.85,-.8) node[right] {$\lambda$};
\draw[thick] (.05,1.8)--(3.75,1.8)
 ..controls (4.8,1.8) and (4.8,.3)..(3.75,.3)--(.05,.3);
\draw[thick,->] (.45,1.8)--(.05,1.8);
\draw[thick,->] (.45,.3)--(.05,.3);
\node at (.45,2.45) {$J$};
\node at (.55,-.1) {$I\ne J$};
\ifnum\panel=0
 \draw[line width=2.2pt] (2.5,1.8)--(3.2,1.8);
 \draw[->,gray!75] (2.55,1.27)--(2.85,1.73);
 \foreach \x/\val in {.6/1,1.1/1,1.8/1,2.5/1,3.2/0,3.75/0}{
  \fill (\x,1.8) circle (1.5pt);
  \node[above,fill=white,inner sep=1.2pt] at (\x,1.98) {$\val$};
 }
\else
 \draw[line width=2.2pt] (1.1,1.8)--(1.8,1.8);
 \draw[->,gray!75] (2.05,1.27)--(1.45,1.73);
 \foreach \x/\val in {.6/1,1.1/1,1.8/0,2.5/0,3.2/0,3.75/0}{
  \fill (\x,1.8) circle (1.5pt);
  \node[above,fill=white,inner sep=1.2pt] at (\x,1.98) {$\val$};
 }
\fi
\node[fill=white,inner sep=1.5pt] at (2.3,1.05) {$\delta f=1$};
\foreach \x in {.6,1.1,1.8,2.5,3.2,3.75}{
 \fill (\x,.3) circle (1.5pt);
 \node[above,fill=white,inner sep=1.2pt] at (\x,.42) {$0$};
}
\fill (4.5375,1.05) circle (1.5pt);
\node[right] at (4.61,1.05) {$0$};
\draw (1.8,1.8) circle (3.6pt);
\draw (2.5,1.8) circle (3.6pt);
\end{scope}
}
\end{tikzpicture}
\caption{Vertex values on an abstract component with ends labeled by $J$ and $I\ne J$. The bold edge in each panel marks the support of $\delta f$, while the value $1$ persists along the upper tail. Moving the cut from $c$ to $c'<c$ changes only the two circled vertex values. Their characteristic function $b=f_{J,c}-f_{J,c'}$ has finite support, so $\delta f_{J,c}-\delta f_{J,c'}=\delta b$ and both coboundaries represent $\alpha_J$.}\label{Fig11}
\end{figure}

Let $E_c$ be the set of edges incident to the cut vertices. It is finite by Lemma~\ref{Lem2.3} and local finiteness. Away from $E_c$, both endpoints of an edge below the cut have the same index set, while both endpoints of an edge above it have value zero. Thus the changes in $f_{J,c}$ occur only next to the cut,
\[
 \supp(\delta f_{J,c})\subset E_c,
 \qquad \delta f_{J,c}\in C_c^1(G),
\]
and give the cohomology class
\begin{equation}\label{Eq4.6}
 \alpha_J=[\delta f_{J,c}]\in H_c^1(G).
\end{equation}
Figure~\ref{Fig11} illustrates this construction on a component with two differently labeled ends, using the same subdivision for two choices of the cut.

Although $\delta f_{J,c}$ has finite support, Theorem~\ref{Th3.2} gives a ray labeled by $J$ with infinitely many vertices where $f_{J,c}=1$. Hence $f_{J,c}$ does not belong to $C_c^0(G)$. The class $\alpha_J$ vanishes precisely when $\delta f_{J,c}$ can also be obtained from a finitely supported vertex function, equivalently,
\[
 \alpha_J=0\quad\Longleftrightarrow\quad
 \delta(f_{J,c}-b)=0\text{ for some }b\in C_c^0(G).
\]
By~\eqref{Eq4.4}, the equation $\delta q=0$ says that $q$ is constant on each connected component. Vanishing of $\alpha_J$ would therefore allow us to make $f_{J,c}$ constant on each component by changing finitely many vertex values. If one component contains both a ray labeled by $J$ and a ray with another label, their different eventual values make this impossible.

For two nonexceptional cuts $c,c'<\lambda_*$, add the finitely many points on both cuts to the vertex set. Moving the cut changes the vertex values only in the intervening strip,
\[
 \supp(f_{J,c}-f_{J,c'})
 \subset V\cap\pr^{-1}([\min\{c,c'\},\max\{c,c'\}]).
\]
The right-hand set is finite because the part of $G$ in this strip is compact and meets finitely many vertices. Hence $f_{J,c}-f_{J,c'}\in C_c^0(G)$ and
\[
 [\delta f_{J,c}]-[\delta f_{J,c'}]
 =[\delta(f_{J,c}-f_{J,c'})]=0.
\]
The class $\alpha_J$ is therefore independent of the cut. It is preserved by subdivision as well, since the refined vertex function restricts to $f_{J,c}$ on the original vertices and the cancellation formula gives $S(\delta\widetilde f_{J,c})=\delta f_{J,c}$. The threshold itself does not affect the eventual labels. Indeed, fix $0<r<R$ and enlarge each $I_i^+$ to an interval $[a_i,b_i]$ by placing its endpoints inside the adjacent negative intervals, using $-1$ or $1$ when the positive interval reaches the boundary. By \ref{E4} and the Dirichlet data, $u(a_i),u(b_i)<r$ for all sufficiently negative $\lambda$, uniformly over the positive solutions. If $\norm u_{L^\infty(I_i^+)}\le R$ and $-\lambda>\norm{a^+}_\infty R^{p-1}$, then
\[
 u''=(-\lambda-au^{p-1})u\ge0\qquad\text{on }[a_i,b_i].
\]
On the positive interval this follows from the assumed bound, and on the added parts it follows from $a<0$. Convexity gives $\norm u_{L^\infty(I_i^+)}<r$. Thus, sufficiently far to the left, every interval maximum is either below $r$ or above $R$. Applying this to any two positive thresholds shows that their vertex functions agree below a common cut. Together with cut and subdivision invariance, this proves that $\alpha_J$ is independent of all these choices.

The labels satisfy one universal relation. Exactly one nonempty index set occurs at each vertex below the cut, where the sum of the $f_{J,c}$ is therefore one, while at and above the cut the sum is zero. It follows that
\[
 b_c:=1-\sum_{J\ne\varnothing}f_{J,c}
 =\mathbf1_{V\cap G_{\ge c}}\in C_c^0(G),
\]
where $\mathbf1$ denotes the characteristic function. Compactness of $G_{\ge c}$ gives finite support. The constant vertex function $1$ satisfies $\delta1=0$, which yields
\begin{equation}\label{Eq4.7}
 \sum_{J\ne\varnothing}\alpha_J=[-\delta b_c]=0.
\end{equation}
Let $N:=2^n-1$ and write each vector $t\in\Ftwo^N$ as $t=(t_J)_{J\ne\varnothing}$, with one coordinate $t_J\in\{0,1\}$ for each nonempty subset $J\subset\{1,\ldots,n\}$. To describe the linear relations among the classes $\alpha_J$, define
\[
 A:\Ftwo^N\longrightarrow H_c^1(G),\qquad
 A(t):=\sum_Jt_J\alpha_J.
\]
The image of $A$ is generated by the classes $\alpha_J$, and its kernel consists of their linear relations. We bound the number of relations by the number of noncompact components, then obtain a complementary bound from the even valences.
\begin{proposition}\label{Prop4.3}
Let $K$ be the number of noncompact connected components of $G$. If $K<\infty$, then
\begin{equation}\label{Eq4.8}
 \dim H_c^1(G,\Ftwo)\ge2^n-1-K,\qquad
 \dim H_c^1(G,\Ftwo)\ge K-1.
\end{equation}
If $K=\infty$, then $\dim H_c^1(G,\Ftwo)=\infty$.
\end{proposition}
\begin{proof}
Within this proof, write $h:=\dim H_c^1(G,\Ftwo)$. Assume first that $K<\infty$, and write $C_1,\ldots,C_K$ for the noncompact components. A vector $t=(t_J)_{J\ne\varnothing}\in\Null{A}$ is a relation among the classes. By the definition of the quotient $H_c^1(G)$, there is a finitely supported vertex function $b_t$ such that
\[
 \sum_{J\ne\varnothing}t_J\delta f_{J,c}=\delta b_t.
\]
Consequently the vertex function
\[
 q_t:=\sum_{J\ne\varnothing}t_Jf_{J,c}-b_t
\]
satisfies $\delta q_t=0$. Formula~\eqref{Eq4.4} says that $q_t$ has one constant value on the vertices of each connected component. Let $\theta_\ell(t)\in\Ftwo$ be its value on $C_\ell$, and collect these constants into a vector,
\[
 \begin{aligned}
 \Theta&:\Null{A}\longrightarrow\Ftwo^K,\\
 \Theta(t)&:=\bigl(\theta_1(t),\ldots,\theta_K(t)\bigr),
 \qquad \Theta(t)_\ell=\theta_\ell(t)=q_t(v)
 \quad(v\in V\cap C_\ell).
 \end{aligned}
\]
We must check that the auxiliary choice of $b_t$ does not change this vector. If $b_t'$ is another choice, then
\[
 \delta(b_t'-b_t)=0,\qquad b_t'-b_t\in C_c^0(G).
\]
The difference is constant on each component and finitely supported. Its value on every noncompact component is therefore zero, which makes $\Theta(t)$ independent of the choice of $b_t$. For $t,s\in\Null{A}$ we may use $b_t+b_s$ to represent $t+s$, giving
\[
 q_{t+s}=q_t+q_s,\qquad\Theta(t+s)=\Theta(t)+\Theta(s).
\]
Thus $\Theta$ is linear.

To prove injectivity, we recover each coordinate $t_J$ from these component constants. Fix a nonempty $J$, and choose the edge ray $r_J$ with label $J$ constructed in Theorem~\ref{Th3.2}. It lies in some $C_{\ell(J)}$ and passes through an infinite sequence of distinct vertices $\{v_{J,m}\}_{m\ge0}$. Along this sequence the parameter tends to $-\infty$. Since $b_t$ is finitely supported, for all sufficiently large $m$ we have
\[
 b_t(v_{J,m})=0,\qquad \pr(v_{J,m})<c,\qquad
 f_{I,c}(v_{J,m})=
 \begin{cases}1,&I=J,\\0,&I\ne J.
 \end{cases}
\]
Evaluating $q_t$ at any such vertex leaves exactly one term of the sum,
\[
 \Theta(t)_{\ell(J)}
 =q_t(v_{J,m})
 =\sum_{I\ne\varnothing}t_If_{I,c}(v_{J,m})-b_t(v_{J,m})
 =t_J.
\]
Thus a zero vector $\Theta(t)$ forces $t_J=0$ for every $J$, which proves injectivity. If several labels have rays in the same component, their coefficients in a relation must all equal that component's constant. Consequently $\dim\Null{A}\le K$, and rank--nullity gives
\[
 h\ge\dim\Range{A}
 =N-\dim\Null{A}\ge N-K.
\]

For the other bound, choose a nonprincipal noncompact component $C$. All its vertices have even valence because it does not contain the unique leaf $\mathfrak{o}$. Define the linear functional that adds the edge values in this component,
\[
 \varepsilon_C:C_c^1(G)\longrightarrow\Ftwo,\qquad
 \varepsilon_C(\eta):=\sum_{e\subset C}\eta(e).
\]
The sum has only finitely many nonzero terms. For a finitely supported vertex function $b$, each value $b(v)$ appears once for each edge incidence at $v$, including two appearances for a loop. Hence
\[
 \begin{aligned}
 \varepsilon_C(\delta b)
 &=\sum_{e\subset C}\bigl(b(v_0(e))+b(v_1(e))\bigr)\\
 &=\sum_{v\in V\cap C}\deg_G(v)b(v)=0
 \qquad\text{in }\Ftwo.
 \end{aligned}
\]
All sums are finite by local finiteness. In particular, changing a cochain by a finite coboundary does not change its total value on $C$. This proves that
\[
 \overline\varepsilon_C:H_c^1(G)\longrightarrow\Ftwo,
 \qquad\overline\varepsilon_C([\eta]):=\varepsilon_C(\eta)
\]
is well defined. A cochain equal to one on a single edge of $C$ has value one under this functional, and therefore cannot be a finite coboundary.

Now choose distinct nonprincipal noncompact components $C_1',\ldots,C_m'$ and an edge $e_j$ in each $C_j'$. Let $\eta_j$ equal one on $e_j$ and zero elsewhere. The functionals distinguish these classes separately,
\[
 \overline\varepsilon_{C_i'}([\eta_j])=
 \begin{cases}1,&i=j,\\0,&i\ne j.
 \end{cases}
\]
Indeed, the marked edge belongs to $C_i'$ exactly when $i=j$. If $\sum_{j=1}^m a_j[\eta_j]=0$, applying $\overline\varepsilon_{C_i'}$ gives
\[
 0=\sum_{j=1}^m a_j\overline\varepsilon_{C_i'}([\eta_j])=a_i
 \qquad(1\le i\le m).
\]
The classes are therefore independent. Taking all $K-1$ nonprincipal noncompact components gives $h\ge K-1$. If there are infinitely many, the same construction works for every finite $m$, proving $h=\infty$. For $K=1$ the second bound is simply $h\ge0$. The principal component cannot be used in this parity argument, since a vertex function supported at $\mathfrak{o}$ has a coboundary with total edge value one.
\end{proof}

Taking the larger of the two bounds eliminates $K$. When $G$ is connected, evaluation on the labeled rays also determines the kernel of $A$.

\begin{corollary}\label{Cor4.4}
One has $\dim H_c^1(G)\ge2^{n-1}-1$. If $G$ is connected, the kernel of $A$ is $\operatorname{span}\{(1,\ldots,1)\}$, giving a canonical injection
\begin{equation}\label{Eq4.9}
 \overline A:\Ftwo^{\,2^n-1}/\operatorname{span}\{(1,\ldots,1)\}\longrightarrow H_c^1(G).
\end{equation}
In particular, the connected case satisfies $\dim H_c^1(G)\ge2^n-2$.
\end{corollary}
\begin{proof}
If $K<\infty$, Proposition~\ref{Prop4.3} gives
\[
 \dim H_c^1(G)\ge\max\{2^n-1-K,K-1\}
 \ge\frac{2^n-2}{2}=2^{n-1}-1.
\]
The infinite case follows from the same proposition. If $G$ is connected, the proof of Proposition~\ref{Prop4.3} shows that all coefficients of a relation $t\in\Null{A}$ equal the single component constant. In the notation of that proof, $t_J=\Theta(t)_1$ for every nonempty $J$, hence
\[
 \Null{A}\subset\operatorname{span}\{\mathbf1\},
 \qquad\mathbf1=(1,\ldots,1).
\]
Relation~\eqref{Eq4.7} gives the reverse inclusion, so the only coefficient vectors in the kernel are the constant ones. Quotienting by their span removes precisely this relation and gives the map in~\eqref{Eq4.9},
\[
 \overline A([t]):=\sum_Jt_J\alpha_J.
\]
Its kernel is zero and its image has dimension $(2^n-1)-1=2^n-2$, proving the claimed embedding and the bound in the connected case.
\end{proof}

This proves Theorem~\ref{Th1.2}. For instance, when $n=2$ and $G$ is connected, the three labels give
\[
 \alpha_{\{1\}}+\alpha_{\{2\}}+\alpha_{\{1,2\}}=0,
 \qquad \dim\operatorname{span}
 \{\alpha_{\{1\}},\alpha_{\{2\}},\alpha_{\{1,2\}}\}=2.
\]
Any two classes are independent and their sum is the third. Further classes may occur in $H_c^1(G)$, so the dimension of the full space need not be two.

\section{Generic components and the equality case}\label{Sec5}
Branching allows a single component to carry several ends. We now vary the weight in the open class $\A_n$ defined in Section~\ref{Sec1}. For a residual subset of this class, the solution set consists of smooth curves, and each nonprincipal noncompact component has exactly two ends. Counting these ends gives the lower bound on $K$. When equality holds, each nonempty interval label occurs at exactly one end.

\subsection{Transversality with respect to the weight}
At a positive solution, surjectivity of the derivative in the pair $(\lambda,u)$ allows the implicit function theorem to describe all nearby solutions by one real coordinate. Varying the weight will provide the directions needed to ensure this surjectivity for generic weights. We call the solution set \emph{totally regular} when
\begin{equation}\tag{R}\label{EqR}
 D_{(\lambda,u)}\F_a(\lambda,u)[\mu,v]=\Lin_{\lambda,u}v-\mu u
 \quad\text{maps }\R\times X\text{ onto }Y
\end{equation}
at every positive solution. Since $\lambda$ is included as a variable, this condition allows the linearization $\Lin_{\lambda,u}$ to be singular at turning points.

At a solution where $\Lin_{\lambda,u}$ is invertible, the derivative in~\eqref{EqR} is onto. In the singular case, let $\phi$ span its kernel. The Fredholm alternative gives
\[
 \Range{\Lin_{\lambda,u}}=\{g\in Y:(g,\phi)_2=0\},\qquad
 \Range{D_{(\lambda,u)}\F_a}
 =\Range{\Lin_{\lambda,u}}+\operatorname{span}\{u\}.
\]
The parameter direction contributes the vector $-u$, which fills the missing range direction precisely when $(u,\phi)_2\ne0$. Thus, in the singular case,~\eqref{EqR} is equivalent to
\begin{equation}\label{Eq5.1}
 \int_{-1}^1u\phi\,dx\ne0.
\end{equation}
\begin{proposition}\label{Prop5.1}
Condition~\eqref{EqR} holds for a residual subset $\mathcal R_n$ of $\A_n$.
\end{proposition}
\begin{proof}
Allowing the weight to vary supplies the additional direction needed to resolve the possible one-dimensional obstruction in the linearized equation. Put $E_r=\mathcal C^r([-1,1])$ and define
\[
 \F:\A_n\times\R\times\Pp\longrightarrow Y,\qquad
 \F(a,\lambda,u):=-u''-\lambda u-au^p.
\]
The set $\A_n$ is open in $E_r$, so every $b\in E_r$ is an admissible infinitesimal weight variation. At a zero $z=(a,\lambda,u)$,
\[
 D\F(z)[b,\mu,v]=\Lin_{\lambda,u}v-\mu u-bu^p.
\]
If $\Lin_{\lambda,u}$ is invertible, any prescribed $g\in Y$ is obtained with $b=\mu=0$ and $v=\Lin_{\lambda,u}^{-1}g$.

Suppose instead that $\Null{\Lin_{\lambda,u}}=\operatorname{span}\{\phi\}$, with $(\phi,\phi)_2=1$. The Fredholm alternative gives
\[
 \Range{\Lin_{\lambda,u}}=\{g\in Y:(g,\phi)_2=0\}.
\]
Only the coefficient in the direction of $\phi$ is missing from this range. Since $u>0$ in $(-1,1)$ and $\phi\not\equiv0$, there is an interior interval on which $u^p\phi$ has a fixed nonzero sign. Choose a nonnegative smooth $b_0$, not identically zero and supported in that interval. Then
\[
 d:=(b_0u^p,\phi)_2\ne0.
\]
To use $b_0$ to supply this coefficient, let $X_1=\{v\in X:(v,\phi)_2=0\}$. The restriction of $\Lin_{\lambda,u}$ to $X_1$ is an isomorphism onto its range, so for any $g\in Y$ we may set
\[
 \begin{aligned}
 b_g&:=-\frac{(g,\phi)_2}{d}b_0,\\
 v_g&:=(\Lin_{\lambda,u}|_{X_1})^{-1}
       \left(g-\frac{(g,\phi)_2}{d}b_0u^p\right).
 \end{aligned}
\]
The argument of the inverse is orthogonal to $\phi$. Substitution gives
\[
 \Lin_{\lambda,u}v_g-b_gu^p=g.
\]
Thus the map $R_z:Y\to E_r\times\R\times X$, $R_zg:=(b_g,0,v_g)$, is a bounded right inverse of $D\F(z)$. In the invertible case, use $R_zg=(0,0,\Lin_{\lambda,u}^{-1}g)$ instead.

This right inverse makes the full solution set $\mathcal Z:=\F^{-1}(0)$ an analytic Banach manifold. Indeed, every variation decomposes uniquely as
\[
 E_r\times\R\times X=\Null{D\F(z)}\oplus\Range{R_z},\qquad
 \xi=\bigl(\xi-R_zD\F(z)\xi\bigr)+R_zD\F(z)\xi.
\]
Both summands are closed, and $D\F(z)$ maps the second isomorphically onto $Y$. The implicit function theorem therefore solves
\[
 \F(z+k+w)=0,
 \qquad k\in\Null{D\F(z)},\quad w\in\Range{R_z},
\]
for $w$ as an analytic function of small $k$. More precisely, there are an open neighborhood $U_z$ of zero in the kernel and an analytic map
\[
 \psi_z:U_z\longrightarrow\Range{R_z},\qquad
 \psi_z(0)=0,\quad D\psi_z(0)=0,
\]
such that the nearby zeros are exactly $z+k+\psi_z(k)$. These are local coordinates in a Banach space, which is the meaning of a Banach manifold in this setting. Its tangent space is
\[
 T_z\mathcal Z=\Null{D\F(z)}.
\]

We now separate the individual weights by the projection
\[
 \wpr:\mathcal Z\longrightarrow\A_n,\qquad
 \wpr(a,\lambda,u):=a.
\]
For brevity, put $\mathcal D_z[\mu,v]:=\Lin_{\lambda,u}v-\mu u$. A tangent vector to $\mathcal Z$ must satisfy the differentiated equation, so
\[
 \begin{aligned}
 T_z\mathcal Z&=\{(b,\mu,v):\mathcal D_z[\mu,v]=bu^p\},\\
 D\wpr(z)[b,\mu,v]&=b,\\
 \Null{D\wpr(z)}&=\{(0,\mu,v):(\mu,v)\in\Null{\mathcal D_z}\},\\
 \Range{D\wpr(z)}&=\{b\in E_r:bu^p\in\Range{\mathcal D_z}\}.
 \end{aligned}
\]

The operator $\mathcal D_z$ is Fredholm of index one. Indeed, regarding $\Lin_{\lambda,u}v$ as an operator on $\R\times X$ leaves its range unchanged and adds one kernel dimension, while $-\mu u$ is a finite-rank perturbation. Its range has codimension at most one, because it contains $\Range{\Lin_{\lambda,u}}$.

If $\mathcal D_z$ is onto, then every $bu^p$ lies in its range. Hence $D\wpr(z)$ is onto and
\[
 \dim\Null{D\wpr(z)}=\dim\Null{\mathcal D_z}=1.
\]
If $\mathcal D_z$ is not onto, the test~\eqref{Eq5.1} gives $(u,\phi)_2=0$, and
\[
 \Range{\mathcal D_z}=\{g:(g,\phi)_2=0\},\qquad
 \Range{D\wpr(z)}=\{b:(bu^p,\phi)_2=0\}.
\]
The second range is the kernel of a nonzero continuous functional, since $b_0$ gives a nonzero value. It is therefore closed and has codimension one. In this case the index formula for $\mathcal D_z$ gives
\[
 \dim\Null{D\wpr(z)}=\dim\Null{\mathcal D_z}=2,
 \qquad \ind D\wpr(z)=2-1=1.
\]
Thus $\wpr$ is a Fredholm map of index one, meaning that all its derivatives are Fredholm operators with that index. The two cases also prove
\[
 D\wpr(z)\text{ is onto}
 \quad\Longleftrightarrow\quad
 D_{(\lambda,u)}\F_a(\lambda,u)\text{ is onto}.
\]
A weight $a$ is a \emph{regular value} of $\wpr$ when the left-hand condition holds at every point of $\wpr^{-1}(a)$, or equivalently when~\eqref{EqR} holds at all its positive solutions.

The Sard--Smale theorem~\cite[Th.~1.3]{Smale} states that the regular values of a $\mathcal C^q$ Fredholm map between Banach manifolds with countable bases form a residual set when $q>\max\{\ind\wpr,0\}$. Here analyticity allows $q=2>1$, and separability of $E_r$, $X$ and $Y$ gives $\mathcal Z$ a countable base. We therefore obtain a residual set $\mathcal R_n\subset\A_n$ of weights satisfying~\eqref{EqR}. It is dense because $\A_n$ is an open subset of a Banach space and hence a Baire space.
\end{proof}

\subsection{Components and their ends}\label{Subsec5.2}
Under total regularity, local analytic coordinates make the solution set a one-dimensional manifold, with $\mathfrak{o}$ as its only boundary point. The classification of these manifolds then determines the number of ends on each noncompact component.

\begin{proposition}\label{Prop5.2}
Under~\eqref{EqR}, the number $K$ of noncompact connected components of $G$ satisfies $K\ge2^{n-1}$, and every nonprincipal noncompact component contains a turning point.
\end{proposition}
\begin{proof}
Let $z_0=(\lambda_0,u_0)\in\Sp$ and put $\mathcal D_0=D_{(\lambda,u)}\F_a(z_0)$. Allowing the additional scalar variable $\lambda$ raises the Fredholm index from zero to one, as in Proposition~\ref{Prop5.1}. Condition~\eqref{EqR} makes the range all of $Y$, so
\[
 \dim\Null{\mathcal D_0}
 =\ind\mathcal D_0+\operatorname{codim}_Y\Range{\mathcal D_0}=1.
\]
Choose a nonzero vector $\zeta_0$ in this kernel and a closed complementary space $W_0$ such that
\[
 \R\times X=\operatorname{span}\{\zeta_0\}\oplus W_0.
\]
The restriction $\mathcal D_0|_{W_0}$ is injective by this splitting and surjective because $\mathcal D_0$ is onto. It is therefore a bounded isomorphism from $W_0$ to $Y$.

Apply the analytic implicit function theorem to the map
\[
 \mathcal H:\R\times W_0\longrightarrow Y,\qquad
 \mathcal H(t,w):=\F_a(z_0+t\zeta_0+w).
\]
At $(0,0)$ its derivative in $w$ is $\mathcal D_0|_{W_0}$. Hence the nearby zeros are uniquely parametrized by
\[
 \begin{aligned}
 w&:(-\varepsilon,\varepsilon)\longrightarrow W_0,
       &w(0)&=0,\quad w'(0)=0,\\
 \gamma&:(-\varepsilon,\varepsilon)\longrightarrow\R\times X,
       &\gamma(t)&:=z_0+t\zeta_0+w(t).
 \end{aligned}
\]
The projection onto $\operatorname{span}\{\zeta_0\}$ recovers $t$ from $\gamma(t)-z_0$, so this parametrization is injective and has a continuous local inverse. Moreover, $\gamma'(0)=\zeta_0\ne0$. Each positive solution thus has an interval as a local coordinate neighborhood, with tangent space $\Null{\mathcal D_0}$. This gives a one-dimensional analytic manifold, and the openness of $\Pp$ ensures that no boundary appears at a positive solution.

There are no isolated positive solutions, and Lemma~\ref{Lem2.2} adds exactly one point to their closure, namely $\mathfrak{o}$. Its unique positive half-branch supplies a half-interval coordinate there. Consequently $G=\Sb$ is a one-dimensional manifold whose only boundary point is $\mathfrak{o}$. The classification of connected one-dimensional manifolds~\cite[App., Th., p.~55]{Milnor} gives
\[
 \begin{array}{c|c|c}
 \text{component}&\text{boundary points}&\text{topological type}\\ \hline
 \text{principal}&\mathfrak{o}&[0,\infty)\\
 \text{nonprincipal, noncompact}&\text{none}&\R\\
 \text{compact}&\text{none}&\text{circle}.
 \end{array}
\]
A compact interval is excluded because it has two boundary points.

These topological types determine the ends. A component is closed in $G$, so a compact subset of $G$ intersects it in a compact set. Under a homeomorphism from $[0,\infty)$ or $\R$, this intersection corresponds to a bounded closed set. Every exterior tail therefore defines a proper ray in $G$, along which $\lambda$ tends to $-\infty$ by~\eqref{Eq2.6}. The principal component has one such end and every other noncompact component has two. Thus, if $K<\infty$,
\[
 e_-(G)=1+2(K-1)=2K-1.
\]
Combining this identity with the end bound yields $2K-1\ge2^n-1$, hence $K\ge2^{n-1}$. The same inequality holds immediately when $K=\infty$.

Finally, let $\gamma:\R\to C$ be a homeomorphism onto a nonprincipal noncompact component and write $\ell(t):=\pr(\gamma(t))$. Both exterior directions are negative ends, so
\[
 \lim_{t\to-\infty}\ell(t)=\lim_{t\to+\infty}\ell(t)=-\infty.
\]
Choose $R>0$ so that $\ell(t)<\ell(0)-1$ for $|t|\ge R$. The maximum of $\ell$ on $[-R,R]$ is then a global maximum and is attained in its interior. Around a maximizing solution, use the regular analytic coordinate $(\lambda(s),u(s))$ constructed above, with the maximum at $s=0$. Differentiating the equation gives
\[
 \dot\lambda(0)=0,\qquad
 \Lin_{\lambda(0),u(0)}\dot u(0)=\dot\lambda(0)u(0)=0.
\]
The tangent vector $(\dot\lambda(0),\dot u(0))$ is nonzero, hence $\dot u(0)\ne0$ and the linearization is singular. The analytic function $\lambda(s)$ cannot be constant on a neighborhood, because that would contradict Lemma~\ref{Lem2.3}. At a maximum its first nonzero Taylor term has even order and negative coefficient, so
\[
 \lambda(s)=\lambda(0)-a_0s^{2m}+O(s^{2m+1}),
 \qquad a_0>0,\quad m\ge1.
\]
It follows that $\lambda(s)<\lambda(0)$ for all sufficiently small nonzero $s$. The maximum is therefore a strict turning point, including the possible higher-order cases.
\end{proof}

When the lower bound on $K$ is attained, the count also determines how often each interval label occurs.

\begin{proposition}\label{Prop5.3}
Under~\eqref{EqR} and $K=2^{n-1}$, the labeling by $J(u)$ gives a bijection between negative ends and nonempty subsets of $\{1,\ldots,n\}$.
\end{proposition}
\begin{proof}
By the classification in the proof of Proposition~\ref{Prop5.2}, the principal half-line contributes one negative end and each of the other $K-1$ lines contributes two. Since $K=2^{n-1}$, this gives $e_-(G)=1+2(K-1)=2^n-1$. On a proper ray $r(t)=(\lambda(t),u(t))$ representing one of these ends, choose $T$ large enough that $\lambda(t)<\lambda_*$ for $t\ge T$. The tail is connected, and the interval norms never equal $\rho$ in that parameter range. Hence $J(u(t))$ is constant for $t\ge T$. Changing $T$ or choosing another ray representing the same end leaves this constant unchanged, since their late tails lie in the same connected component below the cut $\lambda_*$. Thus every end has exactly one well-defined nonempty label.

Let $\mathcal E_-(G)$ denote the finite set of negative ends. The preceding assignment defines
\[
 \begin{aligned}
 \mathcal J&:\mathcal E_-(G)\longrightarrow
 \{J\subset\{1,\ldots,n\}:J\ne\varnothing\},\\
 \mathcal J([r])&:=J(u(t))\quad\text{for sufficiently large }t.
 \end{aligned}
\]
Every nonempty label occurs on an end by Theorem~\ref{Th3.2}, so $\mathcal J$ is surjective. Since its domain and codomain both have $2^n-1$ elements, it is a bijection.
\end{proof}

This completes the proof of Theorem~\ref{Th1.3}.

\section{Morse index and secondary bifurcation}\label{Sec6}
We seek points on the even solution family that are limits of nonsymmetric positive solutions. Every odd eigenvalue of the linearization is positive near the trivial endpoint, whereas a solution at a sufficiently negative parameter has a negative odd direction whenever it is large on a reflected pair of positive intervals. Estimating these directions through the Morse index will produce eigenvalue crossings and hence nonsymmetric solution branches.

\subsection{The Morse index and negative directions}
For a positive solution $(\lambda,u)$, the relevant eigenvalue problem is
\begin{equation}\label{Eq6.1}
 \begin{cases}
 -\phi''-\lambda\phi-pa(x)u^{p-1}\phi=\mu\phi
       &\text{in }(-1,1),\\
 \phi(-1)=\phi(1)=0,
 \end{cases}
 \qquad \phi\not\equiv0.
\end{equation}
Here $\mu$ is the spectral variable and $\lambda$ is the parameter of the nonlinear problem. We use the self-adjoint realization of $\Lin_{\lambda,u}$ in $L^2(-1,1)$ with domain $H^2(-1,1)\cap H_0^1(-1,1)$. The continuous potential and the bounded interval give a discrete spectrum, with simple eigenvalues because the Dirichlet conditions are separated. The eigenvalue sequence $\{\mu_j(\lambda,u)\}_{j\ge1}$ therefore satisfies
\[
 \mu_1(\lambda,u)<\mu_2(\lambda,u)<\cdots,
 \qquad \lim_{j\to\infty}\mu_j(\lambda,u)=+\infty.
\]
The \emph{Morse index} of the solution is
\begin{equation}\label{Eq6.2}
 m(\lambda,u):=\#\{j\ge1:\mu_j(\lambda,u)<0\}.
\end{equation}
Only finitely many eigenvalues are negative, each with multiplicity one. A zero eigenvalue makes the solution \emph{degenerate}, with $\Null{\Lin_{\lambda,u}}\ne\{0\}$, and contributes zero to the Morse index.

The same information can be expressed through the quadratic form
\begin{equation}\label{Eq6.3}
 Q_{\lambda,u}(v):=\int_{-1}^1
 \bigl(|v'|^2-\lambda v^2-pa(x)u^{p-1}v^2\bigr)\,dx,
 \qquad v\in H_0^1(-1,1).
\end{equation}
For an eigenfunction $\phi_j$ normalized by $\norm{\phi_j}_{L^2}=1$, integration by parts gives $Q_{\lambda,u}(\phi_j)=\mu_j$. More generally, orthonormal eigenfunctions satisfy
\[
 Q_{\lambda,u}\left(\sum_{j=1}^{\ell}t_j\phi_j\right)
 =\sum_{j=1}^{\ell}\mu_jt_j^2.
\]
Thus the span of the negative eigenfunctions is a space on which $Q_{\lambda,u}$ is strictly negative away from zero. Any subspace of larger dimension contains a nonzero vector orthogonal to all these eigenfunctions, on which the spectral expansion gives $Q_{\lambda,u}\ge0$. Consequently,
\begin{equation}\label{Eq6.4}
 m(\lambda,u)=\max\left\{\dim W\;\middle|\;
 \begin{gathered}
 W\subset H_0^1(-1,1)\text{ is a linear subspace},\\
 Q_{\lambda,u}(v)<0\text{ for every }v\in W\setminus\{0\}
 \end{gathered}\right\}.
\end{equation}
To prove $m(\lambda,u)\ge k$, it therefore suffices to find $k$ independent variations whose entire span is negative for $Q_{\lambda,u}$. The associated min--max formula is
\begin{equation}\label{Eq6.5}
 \mu_k(\lambda,u)=
 \min_{\substack{W\subset H_0^1(-1,1)\\\dim W=k}}
 \ \max_{v\in W\setminus\{0\}}
 \frac{Q_{\lambda,u}(v)}{\norm v_{L^2}^2}.
\end{equation}
These variational characterizations allow us to estimate the index by constructing suitable test functions.

The form $Q_{\lambda,u}$ is the second variation at $u$ of the energy
\[
 \mathscr E_\lambda(w)=\int_{-1}^1
 \left(\frac{|w'|^2}{2}-\frac{\lambda w^2}{2}
       -\frac{a(x)w^{p+1}}{p+1}\right)\,dx,
 \qquad w\in H_0^1(-1,1).
\]
Thus a negative direction decreases the energy to second order, and the Morse index counts the largest dimension of a subspace on which this decrease holds.

We now construct negative variations localized on the positive weight intervals. For $s>0$ and $v\in H_0^1(-1,1)$, put $E_s(v)=\int_{-1}^1(|v'|^2+sv^2)\,dx$. The following identity separates a negative energy contribution from the error caused by the cutoff.

\begin{lemma}\label{Lem6.1}
For $\chi\in \mathcal C^1([-1,1])$ and a positive solution with $\lambda=-s<0$,
\begin{equation}\label{Eq6.6}
 Q_{\lambda,u}(\chi u)=-(p-1)E_s(\chi u)
       +p\int_{-1}^1|\chi'|^2u^2\,dx.
\end{equation}
\end{lemma}
\begin{proof}
Both $\chi u$ and $\chi^2u$ vanish at the endpoints. Testing $-u''+su=au^p$ with $\chi^2u$ and expanding $| (\chi u)' |^2$ gives
\[
 \int_{-1}^1a\chi^2u^{p+1}\,dx
 =E_s(\chi u)-\int_{-1}^1|\chi'|^2u^2\,dx.
\]
Substituting this identity into $Q_{\lambda,u}(\chi u)=E_s(\chi u)-p\int_{-1}^1a\chi^2u^{p+1}\,dx$ proves the formula. Its positive error term is supported where the cutoff varies.
\end{proof}

Choose a cutoff equal to one on each positive interval and varying only where $a<0$. The uniform decay in~\ref{E4} makes the error in~\eqref{Eq6.6} small, while the maximum above $\rho$ on an interval belonging to $J(u)$ gives a positive lower bound for the localized energy. This comparison yields independent negative directions.

\begin{proposition}\label{Prop6.2}
There exists $\lambda_M<\lambda_*$ such that every positive solution with $\lambda<\lambda_M$ satisfies
\begin{equation}\label{Eq6.7}
 m(\lambda,u)\ge|J(u)|.
\end{equation}
\end{proposition}
\begin{proof}
For each positive interval choose $\chi_i\in\mathcal C^1([-1,1])$ with $0\le\chi_i\le1$, equal to one on $\overline{I_i^+}$, and tapering to zero only inside the neighboring negative intervals. The negative gap between two consecutive positive intervals has positive length. We can therefore place the two transition regions in disjoint closed subintervals of that gap, with a region where both cutoffs vanish between them. This makes the supports of the $\chi_i$ pairwise disjoint. At an exterior endpoint a cutoff may stay one up to the Dirichlet boundary, since multiplication by $u$ still gives a function vanishing there. The cutoffs are fixed independently of the solution and of $\lambda$.

For $i\in J(u)$ put $v_i=\chi_i u\in H_0^1(-1,1)$. Since $\chi_i=1$ on $I_i^+$, the definition of the label gives $\norm{v_i}_\infty>\rho$. Every $v\in H_0^1(-1,1)$ satisfies
\[
 |v(x)|^2
 =\left|\int_{-1}^xv'(y)\,dy\right|^2
 \le (x+1)\int_{-1}^x|v'(y)|^2\,dy
 \le2\int_{-1}^1|v'|^2\,dy.
\]
Thus, for $s=-\lambda>0$,
\[
 E_s(v_i)\ge\int_{-1}^1|v_i'|^2\,dx
 \ge\frac12\norm{v_i}_\infty^2>\frac{\rho^2}{2}.
\]

For $D_i=\supp\chi_i'$, the construction gives a compact subset of the negative region, where the uniform decay estimate~\ref{E4} yields
\[
 \int_{-1}^1|\chi_i'|^2u^2\,dx
 \le\left(\int_{-1}^1|\chi_i'|^2\,dx\right)
       \sup_{x\in D_i}u(x)^2.
\]
If $D_i$ is empty this error is zero. Otherwise its right-hand side tends to zero as $\lambda$ tends to $-\infty$, uniformly over the positive solutions. There are only finitely many cutoffs, so a single $\lambda_M<\lambda_*$ can be chosen such that
\[
 p\int_{-1}^1|\chi_i'|^2u^2\,dx
 <\frac{(p-1)\rho^2}{4}
 \qquad(1\le i\le n,\ \lambda<\lambda_M).
\]
Lemma~\ref{Lem6.1} now yields, for every $i\in J(u)$,
\[
 Q_{\lambda,u}(v_i)
 <-\frac{(p-1)\rho^2}{2}
   +\frac{(p-1)\rho^2}{4}<0.
\]
The nonzero functions $v_i$ have disjoint supports and are therefore independent. Since both $v_iv_j$ and $v_i'v_j'$ vanish almost everywhere for $i\ne j$, their quadratic form is the sum
\[
 Q_{\lambda,u}\left(\sum_{i\in J(u)}c_iv_i\right)
 =\sum_{i\in J(u)}c_i^2Q_{\lambda,u}(v_i)<0
 \quad\text{whenever }(c_i)_{i\in J(u)}\ne0.
\]
Their span is a negative subspace of dimension $|J(u)|$. The characterization~\eqref{Eq6.4} therefore proves the claimed uniform lower bound for the Morse index.
\end{proof}

\subsection{The even component and odd index}\label{Subsec6.2}
Assume from now on that $a$ is even. Reflection preserves the equation,
\[
 \Reflection u(x)=u(-x),\qquad
 \F_a(\lambda,\Reflection u)=\Reflection\F_a(\lambda,u).
\]
We separate even and odd variations by the decomposition
\[
 \begin{aligned}
 X_{\mathrm{ev}}&=\{v\in X:\Reflection v=v\},&
 v_{\mathrm{ev}}&=\tfrac12(v+\Reflection v),\\
 X_{\mathrm{odd}}&=\{v\in X:\Reflection v=-v\},&
 v_{\mathrm{odd}}&=\tfrac12(v-\Reflection v),\\
 X&=X_{\mathrm{ev}}\oplus X_{\mathrm{odd}},&
 v&=v_{\mathrm{ev}}+v_{\mathrm{odd}}.
 \end{aligned}
\]
We use the same parity subscripts for $Y$, $L^2(-1,1)$ and $H_0^1(-1,1)$. At an even solution, the potential $\lambda+pa(x)u^{p-1}$ is even, so the linearization has the block form
\[
 \Lin_{\lambda,u}\Reflection=\Reflection\Lin_{\lambda,u},\qquad
 \Lin_{\lambda,u}=
 \begin{pmatrix}\Lin_{\lambda,u}^{+}&0\\0&\Lin_{\lambda,u}^{-}\end{pmatrix}
 :X_{\mathrm{ev}}\oplus X_{\mathrm{odd}}
 \longrightarrow Y_{\mathrm{ev}}\oplus Y_{\mathrm{odd}}.
\]
The superscripts indicate the restrictions to even and odd functions. Their mixed quadratic term vanishes because its integrand is odd,
\[
 \int_{-1}^1\bigl(v_{\mathrm{ev}}'v_{\mathrm{odd}}'
       -(\lambda+pau^{p-1})v_{\mathrm{ev}}v_{\mathrm{odd}}\bigr)\,dx=0.
\]
Consequently,
\[
 Q_{\lambda,u}(v)=Q_{\lambda,u}(v_{\mathrm{ev}})
                  +Q_{\lambda,u}(v_{\mathrm{odd}}).
\]

Both spectra can be computed on $(0,1)$ using the same differential equation and the boundary condition at zero determined by parity,
\[
 \begin{gathered}
 -\phi''-\lambda\phi-pa(x)u^{p-1}\phi=\mu\phi\quad\text{in }(0,1),\\[3pt]
 \begin{array}{ll}
 \phi'(0)=\phi(1)=0&\text{for the even spectrum},\\
 \phi(0)=\phi(1)=0&\text{for the odd spectrum}.
 \end{array}
 \end{gathered}
\]
Even or odd reflection gives the corresponding eigenfunction on $(-1,1)$. Denote the increasing eigenvalue sequences by $\{\mu_j^+(\lambda,u)\}_{j\ge1}$ and $\{\mu_j^-(\lambda,u)\}_{j\ge1}$, respectively. The corresponding Morse indices are
\begin{equation}\label{Eq6.8}
 m_\pm(\lambda,u):=\#\{j\ge1:\mu_j^\pm(\lambda,u)<0\},
 \qquad m=m_++m_-.
\end{equation}
In the variational form of~\eqref{Eq6.4},
\[
 m_\pm=\max\{\dim W:W\subset H_0^1(-1,1),\
       \Reflection|_W=\pm I,\ Q_{\lambda,u}(v)<0\ (v\in W\setminus\{0\})\}.
\]
The odd index counts negative directions that break the symmetry of an even solution, since for $u$ even and $v\ne0$ odd,
\[
 \Reflection(u+tv)=u-tv\ne u+tv\qquad(t\ne0).
\]
The interval label $J(u)$ allows us to estimate these two indices separately.

Reflection reverses the positive intervals and their indices,
\[
 -I_i^+=(-\tau_i,-\sigma_i)=I_{n+1-i}^+,\qquad
 \Reflection J:=\{n+1-i:i\in J\}.
\]
A label is \emph{invariant under reflection} when $\Reflection J=J$. Every even solution with $\lambda<\lambda_*$ has an invariant label, because
\[
 \norm u_{L^\infty(I_i^+)}=\norm u_{L^\infty(I_{n+1-i}^+)}
 \quad\Longrightarrow\quad
 \bigl(i\in J(u)\ \Longleftrightarrow\ n+1-i\in J(u)\bigr).
\]
The label records which interval maxima exceed the threshold. Establishing symmetry of a solution requires comparison of its full profile with its reflection.

\begin{figure}[!htbp]
\centering
\newcommand{\threeintervals}[3]{%
 \draw[black!40] (-.15,0)--(4.55,0);
 \draw[densely dashed,black!50] (2.2,-.45)--(2.2,.3);
 \foreach \xx/\ii/\shade in {.1/1/#1,1.7/2/#2,3.3/3/#3}{%
  \draw[fill=\shade,line width=.65pt] (\xx,-.18) rectangle (\xx+1,.18);
  \node at (\xx+.5,0) {$I_{\ii}^+$};
 }%
 \node[below] at (2.2,-.25) {$0$};
}
\begin{tikzpicture}[x=1cm,y=1cm,font=\small,>=stealth]
\node[anchor=west] at (-.2,0) {(a) Reflection with a central positive interval};
\begin{scope}[shift={(0,-1.05)}]
 \node at (2.2,.53) {$J=\{1,3\}$};
 \threeintervals{black!30}{white}{black!30}
\end{scope}
\draw[->] (4.95,-1.05)--(6.45,-1.05) node[midway,above] {$\Reflection$};
\begin{scope}[shift={(7,-1.05)}]
 \node at (2.2,.53) {$\Reflection J=J$};
 \threeintervals{black!30}{white}{black!30}
\end{scope}
\begin{scope}[shift={(0,-2.55)}]
 \node at (2.2,.53) {$J=\{1\}$};
 \threeintervals{black!30}{white}{white}
\end{scope}
\draw[->] (4.95,-2.55)--(6.45,-2.55) node[midway,above] {$\Reflection$};
\begin{scope}[shift={(7,-2.55)}]
 \node at (2.2,.53) {$\Reflection J=\{3\}\ne J$};
 \threeintervals{white}{white}{black!30}
\end{scope}

\node[anchor=west] at (-.2,-3.8) {(b) Two positive intervals and $a(0)<0$};
\draw[<->] (3.05,-4.55) .. controls (4.2,-4.15) and (7.2,-4.15) .. (8.35,-4.55);
\node[fill=white,inner sep=2pt] at (5.7,-4.23) {$\Reflection$};
\draw[black!40] (.4,-4.95)--(11,-4.95);
\draw[densely dashed,black!50] (5.7,-4.6)--(5.7,-5.38);
\draw[fill=black!30,line width=.65pt] (1.8,-5.13) rectangle (4.3,-4.77);
\draw[fill=black!30,line width=.65pt] (7.1,-5.13) rectangle (9.6,-4.77);
\node at (3.05,-4.95) {$I_1^+$};
\node at (8.35,-4.95) {$I_2^+$};
\node[below] at (5.7,-5.05) {$0$};
\node at (5.7,-5.77) {$J=\{1,2\}=\Reflection J$};
\end{tikzpicture}
\caption{Filled intervals are those selected by $J$. In (a), reflection exchanges $I_1^+$ and $I_3^+$ and preserves the central interval $I_2^+$. The first selection is invariant, with $q=1$ and $\epsilon=0$, while the second is not. In (b), zero lies in a negative gap. There is no central positive interval, and a nonempty invariant label must select both intervals.}\label{Fig12}
\end{figure}
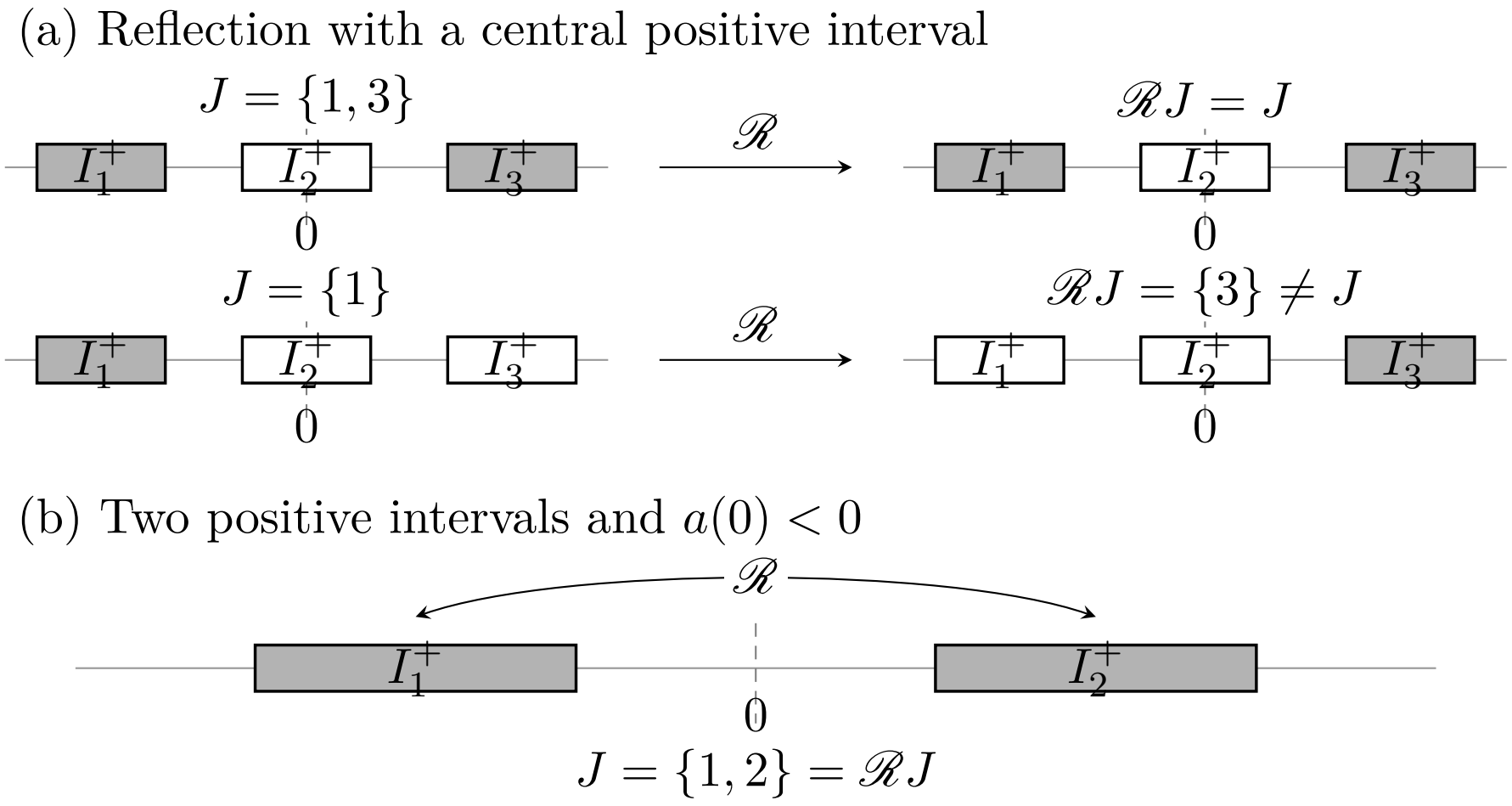

A \emph{central positive interval} is one fixed by reflection. It exists precisely when $n$ is odd, in which case its index is $i_0=(n+1)/2$ and $0\in I_{i_0}^+$. For an invariant label define
\[
 \begin{aligned}
 q(J)&:=\#\{i\in J:i<n+1-i\},\\
 \epsilon(J)&:=
 \begin{cases}1,&n\text{ is odd and }(n+1)/2\in J,\\
 0,&\text{otherwise}.
 \end{cases}
 \end{aligned}
\]
Thus $q(J)$ counts reflected pairs, while $\epsilon(J)$ records the central index, and
\[
 |J|=2q(J)+\epsilon(J).
\]
Figure~\ref{Fig12} compares an invariant label with a noninvariant one. A reflected pair will supply one even and one odd negative direction, while the central interval supplies an even direction.

\begin{proposition}\label{Prop6.3}
After decreasing the threshold $\lambda_M$ of Proposition~\ref{Prop6.2}, every even positive solution with $\lambda<\lambda_M$ satisfies
\begin{equation}\label{Eq6.9}
 m_-(\lambda,u)\ge q(J(u)),\qquad
 m_+(\lambda,u)\ge q(J(u))+\epsilon(J(u)).
\end{equation}
If $a(0)<0$, then $m_-(\lambda,u)\ge1$ at every sufficiently negative parameter.
\end{proposition}
\begin{proof}
Choose the cutoffs of Proposition~\ref{Prop6.2} symmetrically, so that
\[
 \chi_{n+1-i}=\Reflection\chi_i,\qquad
 v_i:=\chi_i u,\qquad v_{n+1-i}=\Reflection v_i.
\]
Their supports are pairwise disjoint, and $Q_{\lambda,u}(v_i)<0$ for $i\in J(u)$. Write
\[
 J_{\mathrm L}:=\{i\in J(u):i<n+1-i\},\qquad
 \#J_{\mathrm L}=q(J(u)).
\]
For each $i\in J_{\mathrm L}$ set
\[
 w_i^+:=\frac{v_i+v_{n+1-i}}{\sqrt2},\qquad
 w_i^-:=\frac{v_i-v_{n+1-i}}{\sqrt2},\qquad
 \Reflection w_i^\pm=\pm w_i^\pm.
\]
Disjoint supports and reflection invariance give
\[
 Q_{\lambda,u}(w_i^\pm)
 =\tfrac12\bigl(Q_{\lambda,u}(v_i)+Q_{\lambda,u}(v_{n+1-i})\bigr)
 =Q_{\lambda,u}(v_i)<0.
\]
Different pairs also have disjoint supports. Hence, for either sign and any nonzero coefficient vector $(c_i)_{i\in J_{\mathrm L}}$,
\[
 Q_{\lambda,u}\Bigl(\sum_{i\in J_{\mathrm L}}c_iw_i^\pm\Bigr)
 =\sum_{i\in J_{\mathrm L}}c_i^2Q_{\lambda,u}(v_i)<0.
\]
Their spans are therefore negative parity subspaces of dimension $q(J(u))$. If $\epsilon(J(u))=1$, choose the central cutoff even. The additional function $v_{i_0}$ then satisfies
\[
 \Reflection v_{i_0}=v_{i_0},\qquad Q_{\lambda,u}(v_{i_0})<0,
 \qquad\supp v_{i_0}\cap\supp w_i^\pm=\varnothing.
\]
It increases the dimension of the negative even subspace by one, proving both bounds through~\eqref{Eq6.4}.

If $a(0)<0$, every positive interval has a distinct reflected partner. The label of an even solution at sufficiently negative $\lambda$ is nonempty and invariant under reflection, so it contains at least one such pair. The corresponding odd test function gives the required negative direction.
\end{proof}

Recall that $\Cs$ is the connected component containing $\mathfrak{o}$ in
\[
 \{(\lambda,u)\in\Sp:\Reflection u=u\}\cup\{\mathfrak{o}\}.
\]
It is the principal component within the even solution set, which may be smaller than the principal component of the full graph $G$. To use the preceding index bounds, we need $\Cs$ to reach arbitrarily negative parameters and to join a small even solution to one in that region by a path of positive even solutions. Both properties follow from the mixed problem on $(0,1)$.

\begin{lemma}\label{Lem6.4}
The component $\Cs$ contains a proper ray with $\lambda\to-\infty$, and $\Cs\setminus\{\mathfrak{o}\}$ is connected.
\end{lemma}
\begin{proof}
Even solutions correspond to the mixed problem
\[
 \begin{cases}
 -w''=\lambda w+a(x)w^p&\text{in }(0,1),\\
 w'(0)=w(1)=0,\qquad w>0&\text{in }[0,1).
 \end{cases}
\]
More precisely, for $X_{\mathrm{mix}}=\{w\in\mathcal C^2([0,1]):w'(0)=w(1)=0\}$, restriction and even extension give inverse bounded linear maps
\[
 \begin{aligned}
 \operatorname{res}&:X_{\mathrm{ev}}\longrightarrow X_{\mathrm{mix}},&\operatorname{res}(u)&:=u|_{[0,1]},\\
 \operatorname{ext}&:X_{\mathrm{mix}}\longrightarrow X_{\mathrm{ev}},&\operatorname{ext}(w)(x)&:=w(|x|).
 \end{aligned}
\]
The condition $w'(0)=0$ ensures that the extension is $\mathcal C^2$ at zero. The mixed linearization is
\[
 \Lin^{\mathrm{mix}}_{\lambda,w}:X_{\mathrm{mix}}\to\mathcal C([0,1]),\qquad
 \Lin^{\mathrm{mix}}_{\lambda,w}v:=-v''-\lambda v-paw^{p-1}v.
\]
It is Fredholm of index zero, being a compact perturbation of the isomorphism $-\partial_x^2:X_{\mathrm{mix}}\to\mathcal C([0,1])$. Cauchy uniqueness gives
\[
 v\in\Null{\Lin^{\mathrm{mix}}_{\lambda,w}},\quad v(0)=0
 \quad\Longrightarrow\quad v(0)=v'(0)=0
 \quad\Longrightarrow\quad v\equiv0,
\]
so its kernel has dimension at most one.

By reflection, Lemmas~\ref{Lem2.2} and~\ref{Lem2.3} give compactness in bounded parameter strips and finite parameter fibers for this mixed problem. The reduction and link argument of Proposition~\ref{Prop2.4} therefore give a locally finite graph of nonisolated even solutions with
\[
 \deg_{\mathrm{ev}}(z)\in2\N
 \qquad\text{at every nonisolated positive even solution }z.
\]
Here the valence counts only even half-arcs. The local branch at $\mathfrak{o}$ is even as well. Its parameter is $s=(u,\varphi_1)_2$, where $\varphi_1$ is the normalized positive first Dirichlet eigenfunction. Reflection leaves this parameter unchanged,
\[
 (\Reflection u,\varphi_1)_2=(u,\Reflection\varphi_1)_2
 =(u,\varphi_1)_2=s,
\]
and the uniqueness in Lemma~\ref{Lem2.2} gives $\Reflection u(s)=u(s)$. Thus $\Cs$ is a connected graph with
\[
 \deg_{\Cs}(\mathfrak{o})=1,\qquad
 \deg_{\Cs}(z)\in2\N\quad(z\ne\mathfrak{o}).
\]
If it were compact, its finite edge set would satisfy
\[
 2|E(\Cs)|=\sum_{v\in V(\Cs)}\deg_{\Cs}(v)
 \equiv1\pmod2,
\]
which is impossible. It is therefore noncompact and contains a simple infinite edge ray $r:[0,\infty)\to\Cs$, by the finite-path construction in Theorem~\ref{Th3.2}. Local finiteness makes $r$ proper. For every $c\in\R$,
\[
 \Cs\cap G_{\ge c}\text{ is compact}
 \quad\Longrightarrow\quad
 r^{-1}(\Cs\cap G_{\ge c})\subset[0,T_c]
\]
for some $T_c<\infty$. Hence $\pr(r(t))<c$ for $t>T_c$, proving $\lim_{t\to\infty}\pr(r(t))=-\infty$.

Finally, any two points of $\Cs\setminus\{\mathfrak{o}\}$ are joined by a finite simple path in $\Cs$. A leaf cannot lie in the interior of such a path, since an interior point has two distinct incident path segments. The joining path therefore avoids $\mathfrak{o}$, proving path connectedness of the punctured component.
\end{proof}

\subsection{Odd eigenvalue crossings and symmetry breaking}\label{Subsec6.3}
We now pass from negative odd directions to nonsymmetric solutions. The argument follows an odd eigenvalue along the even family and reduces the equation near a sign-changing zero to a scalar equation in $\lambda$ and an odd amplitude $t$.

For two even solutions $z=(\lambda,u)$ and $\widetilde z=(\widetilde\lambda,\widetilde u)$ at the same weight, put
\[
 d(z,\widetilde z):=|\lambda-\widetilde\lambda|
       +p\norm{a(u^{p-1}-\widetilde u^{p-1})}_\infty.
\]
Their quadratic forms satisfy
\[
 |Q_z(v)-Q_{\widetilde z}(v)|
 \le d(z,\widetilde z)\norm v_{L^2}^2.
\]
Applying~\eqref{Eq6.5} in each parity subspace and in both directions gives
\begin{equation}\label{Eq6.10}
 |\mu_j^\pm(z)-\mu_j^\pm(\widetilde z)|\le d(z,\widetilde z)
 \qquad(j\ge1).
\end{equation}
In particular, if $g:=\inf_j|\mu_j^-(z)|>0$, then
\[
 d(z,\widetilde z)<g/2
 \quad\Longrightarrow\quad
 \operatorname{sgn}\mu_j^-(\widetilde z)=\operatorname{sgn}\mu_j^-(z)
 \quad\text{for every }j.
\]
Thus the odd index can change only at a zero odd eigenvalue.

At such a zero, the full Sturm--Liouville kernel is spanned by an odd eigenfunction. Simplicity and the parity decomposition imply
\[
 \Null{\Lin_{\lambda_0,u_0}}
 =\Null{\Lin^-_{\lambda_0,u_0}}\oplus\Null{\Lin^+_{\lambda_0,u_0}},
 \qquad\dim\Null{\Lin^-_{\lambda_0,u_0}}=1,
 \qquad\Null{\Lin^+_{\lambda_0,u_0}}=\{0\}.
\]
The even restriction has Fredholm index zero and is therefore invertible. The implicit function theorem makes the nearby even solutions a single analytic graph $(\lambda,u_{\mathrm s}(\lambda))$. A zero of
\[
 \mu(\lambda):=\mu_j^-(\lambda,u_{\mathrm s}(\lambda))
\]
is called \emph{sign-changing} if, for some $\varepsilon>0$,
\[
 \mu(\lambda_0-s)\mu(\lambda_0+s)<0
 \qquad(0<s<\varepsilon).
\]
The proof below also establishes analyticity of $\mu$. For a nonidentically zero analytic eigenvalue, this condition means that its first nonzero Taylor term has odd order. It does not require $\mu'(\lambda_0)\ne0$.

Near the trivial endpoint, the odd eigenvalues are positive,
\begin{equation}\label{Eq6.11}
 \lim_{(\lambda,u)\to\mathfrak{o}}\mu_j^-(\lambda,u)
 =(j\pi)^2-\lambda_1>0,\qquad j\ge1.
\end{equation}
This follows from the odd Dirichlet problem on $(0,1)$ and~\eqref{Eq6.10}. If $m_-$ later reaches $k$, the first $k$ odd eigenvalues have become negative. The resulting sign changes on $\Cs$ connect the odd-index bound with the spectral bifurcation mechanism of~\cite[Th.~2.1]{PW}.

\begin{lemma}\label{Lem6.5}
At an even positive solution, if $\mu_j^-(\lambda_0,u_0)=0$ changes sign along the local even graph, then $(\lambda_0,u_0)$ is a symmetry-breaking point with at least two nonsymmetric half-arcs, exchanged by reflection, in addition to the two even half-arcs.
\end{lemma}
\begin{proof}
We apply the Lyapunov--Schmidt reduction of Section~\ref{Subsec2.2} and~\cite[Th.~6.1]{LGB}, now centered along the local even branch rather than at a fixed function. Reflection will impose an odd symmetry on the reduced scalar equation.

Choose an odd function $\phi\in\Null{\Lin_{\lambda_0,u_0}}$ normalized by $(\phi,\phi)_2=1$. The even restriction is invertible, so the implicit function theorem gives
\[
 u_{\mathrm s}:I\longrightarrow X_{\mathrm{ev}},\qquad
 \F_a(\lambda,u_{\mathrm s}(\lambda))=0,\qquad
 u_{\mathrm s}(\lambda_0)=u_0,
\]
where $I$ is an open interval containing $\lambda_0$. We use
\[
 \Lin_{\mathrm s}(\lambda):=\Lin_{\lambda,u_{\mathrm s}(\lambda)}:X\to Y,
 \qquad \Lin_0:=\Lin_{\mathrm s}(\lambda_0).
\]
Thus $\Lin_{\mathrm s}(\lambda)$ is the linearization along the even branch, and $\Lin_0$ is its value at the bifurcation point.

As in~\eqref{Eq2.2}, separate the kernel direction from its orthogonal complement,
\[
 \begin{aligned}
 X_1&:=\{w\in X:(w,\phi)_2=0\},&
 X&=\operatorname{span}\{\phi\}\oplus X_1,\\
 Y_0&:=\{g\in Y:(g,\phi)_2=0\},&
 Y&=\operatorname{span}\{\phi\}\oplus Y_0,\\
 Q&:Y\longrightarrow Y_0,&
 Qg&:=g-(g,\phi)_2\phi.
 \end{aligned}
\]
The Fredholm alternative gives $\Range{\Lin_0}=Y_0$. Removing its kernel makes the restriction
\[
 \Lin_0|_{X_1}:X_1\longrightarrow Y_0
\]
an isomorphism. For a function near the even branch, write
\[
 u=u_{\mathrm s}(\lambda)+t\phi+w,
 \qquad t\in\R,\quad w\in X_1.
\]
The equation $\F_a(\lambda,u)=0$ is equivalent to the two conditions
\[
 Q\F_a(\lambda,u)=0,
 \qquad (\F_a(\lambda,u),\phi)_2=0.
\]
We first solve the equation with values in $Y_0$. Define
\[
 \F_{\mathrm r}:I\times\R\times X_1\longrightarrow Y_0,
 \qquad
 \F_{\mathrm r}(\lambda,t,w):=
 Q\F_a(\lambda,u_{\mathrm s}(\lambda)+t\phi+w).
\]
Its derivative in $w$ at $(\lambda_0,0,0)$ is $\Lin_0|_{X_1}$. The analytic implicit function theorem gives a unique small correction
\[
 \Yred:\mathcal U\longrightarrow X_1,
 \qquad \mathcal U:=I\times(-\delta,\delta),
 \qquad \F_{\mathrm r}(\lambda,t,\Yred(\lambda,t))=0,
\]
after shrinking $I$ and choosing $\delta>0$. Since $u_{\mathrm s}(\lambda)$ already solves the equation, uniqueness implies $\Yred(\lambda,0)=0$. Differentiation at $(\lambda_0,0)$ also gives
\[
 \Lin_0|_{X_1}\,\partial_t\Yred(\lambda_0,0)
 =-Q\Lin_0\phi=0,
 \qquad \partial_t\Yred(\lambda_0,0)=0.
\]
The candidate solution and the remaining scalar equation are
\[
 \begin{aligned}
 U&:\mathcal U\longrightarrow X,&
 U(\lambda,t)&:=u_{\mathrm s}(\lambda)+t\phi+\Yred(\lambda,t),\\
 \Red&:\mathcal U\longrightarrow\R,&
 \Red(\lambda,t)&:=(\F_a(\lambda,U(\lambda,t)),\phi)_2.
 \end{aligned}
\]
Because $Q\F_a(\lambda,U(\lambda,t))=0$, the value of the full equation has no component in $Y_0$. Its only remaining component is the multiple of $\phi$ given by
\begin{equation}\label{Eq6.12}
 \F_a(\lambda,U(\lambda,t))=\Red(\lambda,t)\phi,
 \qquad (U(\lambda,t)-u_{\mathrm s}(\lambda),\phi)_2=t.
\end{equation}
Consequently $U(\lambda,t)$ solves the original equation exactly when $\Red(\lambda,t)=0$. Since $t$ is its coefficient in the odd kernel direction, each solution with $t\ne0$ is nonsymmetric.

Reflection preserves $X_1$, commutes with $Q$ and sends $\phi$ to $-\phi$. Applying reflection to the range equation gives
\[
 Q\F_a\bigl(\lambda,u_{\mathrm s}(\lambda)-t\phi
                         +\Reflection\Yred(\lambda,t)\bigr)=0.
\]
Uniqueness of the correction at $(\lambda,-t)$ therefore implies
\[
 \Yred(\lambda,-t)=\Reflection\Yred(\lambda,t),\qquad
 U(\lambda,-t)=\Reflection U(\lambda,t).
\]
Taking the inner product of the reflected equation with $\phi$ yields
\[
 \Red(\lambda,-t)
 =(\Reflection\F_a(\lambda,U(\lambda,t)),\phi)_2
 =(\F_a(\lambda,U(\lambda,t)),-\phi)_2
 =-\Red(\lambda,t).
\]
Hence the Taylor series contains only odd powers of $t$. Factoring out $t$ and writing $z=t^2$ gives an analytic function
\begin{equation}\label{Eq6.13}
 \begin{aligned}
 h&:I\times(-\delta^2,\delta^2)\longrightarrow\R,&
 h(\lambda,z)&:=\sum_{m\ge0}a_m(\lambda)z^m,\\
 \Red(\lambda,t)&=t\,h(\lambda,t^2),&
 h(\lambda,0)&=\partial_t\Red(\lambda,0),
 \end{aligned}
\end{equation}
where $\Red(\lambda,t)=\sum_{m\ge0}a_m(\lambda)t^{2m+1}$ and the domain may be reduced. The even branch is given by $t=0$, while the nonsymmetric solutions correspond to zeros of $h(\lambda,t^2)$ with $t\ne0$.

To relate $h(\lambda,0)$ to the crossing odd eigenvalue, we construct an eigenpair along $u_{\mathrm s}$ through the map
\[
 \begin{aligned}
 \mathscr H&:I\times X_{\mathrm{odd}}\times\R
                \longrightarrow Y_{\mathrm{odd}}\times\R,\\
 \mathscr H(\lambda,e,\nu)&:=
       (\Lin_{\mathrm s}(\lambda)e-\nu e,(e,\phi)_2-1).
 \end{aligned}
\]
The second coordinate fixes the normalization of the eigenfunction. At $(\lambda_0,\phi,0)$, the derivative in $(e,\nu)$ sends $(w,\eta)$ to
\[
 (\Lin_0w-\eta\phi,(w,\phi)_2).
\]
For a prescribed $(g,r)\in Y_{\mathrm{odd}}\times\R$, taking the inner product of the first equation with $\phi$ forces $\eta=-(g,\phi)_2$. Writing $w=r\phi+w_\perp$ then leaves $\Lin_0w_\perp=Qg$, so the unique solution is
\[
 w=r\phi+(\Lin_0|_{X_1})^{-1}Qg,
 \qquad \eta=-(g,\phi)_2.
\]
The inverse preserves odd functions because $\Lin_0$ and $Q$ commute with reflection. Thus this derivative is an isomorphism on the stated odd spaces. The implicit function theorem gives analytic maps $e:I\to X_{\mathrm{odd}}$ and $\mu:I\to\R$ with
\[
 \Lin_{\mathrm s}(\lambda)e(\lambda)=\mu(\lambda)e(\lambda),
 \quad(e(\lambda),\phi)_2=1,
 \quad e(\lambda_0)=\phi,\quad\mu(\lambda_0)=0.
\]
Simplicity and continuity identify $\mu(\lambda)$ with $\mu_j^-(\lambda,u_{\mathrm s}(\lambda))$.

Now differentiate the first identity in~\eqref{Eq6.12} with respect to $t$ at $t=0$. Since $U(\lambda,0)=u_{\mathrm s}(\lambda)$, the chain rule gives
\[
 \begin{aligned}
 v(\lambda)&:=\partial_tU(\lambda,0)
             =\phi+\partial_t\Yred(\lambda,0),\\
 \Lin_{\mathrm s}(\lambda)v(\lambda)
       &=\partial_t\Red(\lambda,0)\phi=h(\lambda,0)\phi,
 \qquad v(\lambda_0)=\phi.
 \end{aligned}
\]
Taking the inner product with $e(\lambda)$ and using its normalization and self-adjointness gives
\[
 \begin{aligned}
 h(\lambda,0)
 &=h(\lambda,0)(\phi,e(\lambda))_2\\
 &=(\Lin_{\mathrm s}(\lambda)v(\lambda),e(\lambda))_2\\
 &=(v(\lambda),\Lin_{\mathrm s}(\lambda)e(\lambda))_2
   =\mu(\lambda)(v(\lambda),e(\lambda))_2.
 \end{aligned}
\]
The continuous factor $\kappa(\lambda):=(v(\lambda),e(\lambda))_2$ equals one at $\lambda_0$ and remains positive on a smaller interval, where
\begin{equation}\label{Eq6.14}
 h(\lambda,0)=\kappa(\lambda)\mu(\lambda),\qquad
 \kappa(\lambda_0)=1,\qquad\kappa(\lambda)>0.
\end{equation}
The assumed eigenvalue crossing is consequently a sign change of $h(\lambda,0)$.

Take $\lambda^-<\lambda_0<\lambda^+$ arbitrarily close to $\lambda_0$ with $h(\lambda^-,0)$ and $h(\lambda^+,0)$ of opposite signs. For all sufficiently small $t>0$, the values $h(\lambda^-,t^2)$ and $h(\lambda^+,t^2)$ retain those signs. The intermediate value theorem gives a zero $h(\lambda,t^2)=0$ between them. This yields nonsymmetric solutions with $t>0$ arbitrarily close to $(\lambda_0,u_0)$, all positive by openness of $\Pp$.

The reduced zero set has finitely many incident half-arcs by Proposition~\ref{Prop2.4}. One of them therefore contains points with $t>0$ accumulating at the crossing. In an analytic parametrization of that half-arc the coordinate $t$ has a fixed sign near the endpoint, so the half-arc lies in $t>0$. Its reflection lies in $t<0$, giving the two nonsymmetric half-arcs in addition to the even ones on $t=0$.
\end{proof}

We apply the crossing lemma to the first $k$ odd eigenvalues along a path in the even component.

\begin{proposition}\label{Prop6.6}
If $\Cs$ contains a solution with $m_-(\lambda,u)\ge k$ for an integer $k\ge1$, then it contains at least $k$ distinct symmetry-breaking points.
\end{proposition}
\begin{proof}
Let $E=\Cs\setminus\{\mathfrak{o}\}$. By Lemma~\ref{Lem6.4}, it is path connected. Choose $z_1\in E$ with $m_-(z_1)\ge k$ and $z_0\in E$ sufficiently near $\mathfrak{o}$. Formula~\eqref{Eq6.11} gives
\[
 0<\mu_1^-(z_0)<\mu_2^-(z_0)<\cdots,
 \qquad \mu_1^-(z_1)<\cdots<\mu_k^-(z_1)<0.
\]
Join them by a simple path $\gamma:[0,1]\to E$ with $\gamma(0)=z_0$ and $\gamma(1)=z_1$.

We first exclude intervals of zero eigenvalues. For each $j$ define
\[
 W_j:=\{z\in E:\mu_j^-|_U\equiv0
               \text{ for some relative neighborhood }U\text{ of }z\}.
\]
The set is open by definition. To prove closedness, let $\{z_\ell\}_{\ell\ge1}\subset W_j$ converge to $z\in E$, so continuity gives $\mu_j^-(z)=0$. Invertibility of the even linearization provides a local analytic graph $(\lambda,u_{\mathrm s}(\lambda))$ through $z$. The eigenpair construction in the preceding proof applies to every simple odd eigenvalue, giving the analytic function
\[
 \vartheta_j:I\longrightarrow\R,\qquad
 \vartheta_j(\lambda):=\mu_j^-(\lambda,u_{\mathrm s}(\lambda)).
\]
For large $\ell$, the neighborhood witnessing $z_\ell\in W_j$ contains a nonempty open subinterval of this graph. The identity theorem yields
\[
 \vartheta_j=0\text{ on that subinterval}
 \quad\Longrightarrow\quad\vartheta_j\equiv0\text{ on }I
 \quad\Longrightarrow\quad z\in W_j.
\]
Thus $W_j$ is also closed. Since $E$ is connected and $\mu_j^-(z_0)>0$,
\[
 W_j=\varnothing.
\]
Every zero of $\mu_j^-$ is consequently isolated in its local even graph. The set
\[
 Z_j:=\{z\in E:\mu_j^-(z)=0\}
\]
is closed and discrete, so $Z_j\cap\gamma([0,1])$ is finite.

For $1\le j\le k$, the function $\mu_j^-\circ\gamma$ starts positive, ends negative and has only finitely many zeros. At one of them, $s_j\in(0,1)$, it therefore changes sign,
\[
 \mu_j^-(\gamma(s_j-\varepsilon))
 \mu_j^-(\gamma(s_j+\varepsilon))<0
 \qquad\text{for all sufficiently small }\varepsilon>0.
\]
Near $\gamma(s_j)$, the even set is a graph over $\lambda$. The simple path is injective in this graph, so $\pr\circ\gamma$ is strictly monotone on a small interval about $s_j$. Thus $\mu_j^-$ changes sign across $\pr(\gamma(s_j))$ along the even branch. Lemma~\ref{Lem6.5} gives symmetry breaking at $\gamma(s_j)$.

The crossing points are distinct. Otherwise, for $i\ne j$ and $z=\gamma(s_i)=\gamma(s_j)$ we would have
\[
 \mu_i^-(z)=\mu_j^-(z)=0,
\]
contradicting the strict ordering of the odd eigenvalues. We have therefore obtained at least $k$ distinct points in the solution space, possibly sharing parameter values.
\end{proof}

The preceding proposition and Lemma~\ref{Lem6.5} establish the counting and half-arc assertions of Theorem~\ref{Th1.4}. If $a(0)<0$, Lemma~\ref{Lem6.4} gives even solutions on $\Cs$ at arbitrarily negative parameters, where Proposition~\ref{Prop6.3} yields $m_-\ge1$ and hence symmetry breaking.

More generally, the same argument yields at least $k$ bifurcation points whenever $\Cs$ contains a solution with $q(J(u))\ge k$ and $\lambda<\lambda_M$.

\section{Examples and further questions}\label{Sec7}
Two even polynomial weights illustrate how the positive intervals determine labels and constrain symmetry breaking. We then compare the topological bounds with abstract graphs attaining them.

\subsection{Polynomial weights}
Consider the even polynomial
\[
 a_2(x)=\left(x^2-\frac1{16}\right)\left(\frac9{16}-x^2\right).
\]
Its positive intervals $(-3/4,-1/4)$ and $(1/4,3/4)$ have simple endpoints and a central negative gap, as in Figure~\ref{Fig12}(b). For every integer $p\ge2$, the end, branching and cohomological bounds give at least three negative ends, $K+B_{\mathrm{nc}}\ge2$ and $\dim H_c^1(G)\ge1$. If the graph is connected, it has at least one branching point and $\dim H_c^1(G)\ge2$.

Here $\{1,2\}$ is the only nonempty index set invariant under reflection. The odd index is zero near $\mathfrak{o}$ and at least one on the sufficiently negative part of $\Cs$, so Theorem~\ref{Th1.4} gives a secondary bifurcation point with at least two nonsymmetric positive half-arcs exchanged by $\Reflection$.

The singleton labels $\{1\}$ and $\{2\}$ also give two distinct negative ends consisting of nonsymmetric solutions, since reflection exchanges the labels. Thus the end count supplies unbounded nonsymmetric families. Whether these ends belong to the branches issuing from a particular symmetry-breaking point remains open.

\begin{figure}[!htbp]
\centering
\includegraphics{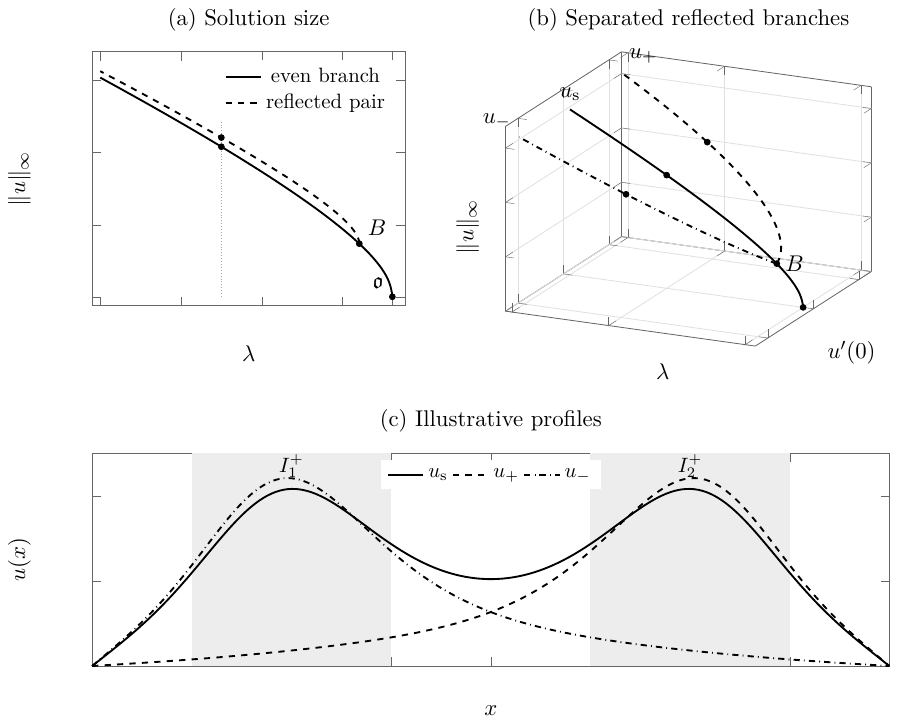}
\caption{Qualitative sketch of one possible symmetry-breaking configuration at $B$. The reflected branches overlap in (a) and separate in (b) through their opposite central slopes. Panel (c) shows the three marked profiles, with shaded positive weight intervals.}\label{Fig13}
\end{figure}

Figure~\ref{Fig13} illustrates the reflection geometry of a possible symmetry-breaking configuration. Since $a_2$ is even, Cauchy uniqueness at $x=0$ gives
\[
 u'(0)=0\quad\Longleftrightarrow\quad\Reflection u=u,
 \qquad (\Reflection u)'(0)=-u'(0).
\]
Reflected solutions therefore have equal uniform norms and opposite central slopes. Their projections onto the $(\lambda,\norm u_\infty)$ plane coincide, while the coordinate $u'(0)$ separates them. Reflection fixes the plane $u'(0)=0$ containing the even branch and exchanges the nonsymmetric pair across it.

For a residual set of unrestricted perturbations near $a_2$ in $\A_2$, Theorem~\ref{Th1.3} gives at least two noncompact components. If there are exactly two, the three ends carry the three distinct labels $\{1\}$, $\{2\}$ and $\{1,2\}$.

For comparison, the weight
\[
 a_3(x)=\left(\frac{25}{36}-x^2\right)
 \left(x^2-\frac14\right)\left(x^2-\frac1{36}\right)
\]
has three positive intervals $(-5/6,-1/2)$, $(-1/6,1/6)$ and $(1/2,5/6)$, with simple zeros and negative gaps. It gives at least seven ends, $K+B_{\mathrm{nc}}\ge4$ and $\dim H_c^1(G)\ge3$. The four labels $\{1\}$, $\{3\}$, $\{1,2\}$ and $\{2,3\}$ are not invariant under reflection. Theorem~\ref{Th3.2} therefore gives four distinct negative ends represented by rays of nonsymmetric solutions. For a solution on $\Cs$, the two outer intervals enter the label together because reflection exchanges them. At sufficiently negative parameters, such a pair gives $m_-\ge1$ by~\eqref{Eq6.9} and hence symmetry breaking by Theorem~\ref{Th1.4}. Thus, in the absence of symmetry breaking, the only possible nonempty label there is $\{2\}$.

\subsection{Sharpness and realization}
The component and cohomological bounds are sharp among abstract graphs satisfying the valence and end constraints of $G$. They are attained by
\[
 G_{\min}=[0,\infty)\ \sqcup\ \bigsqcup_{j=1}^{2^{n-1}-1}\R.
\]
This graph has $2^n-1$ ends, $\dim H_c^1(G_{\min})=2^{n-1}-1$ and $b_1(G_{\min})=0$. Every end becomes negative upon setting $\lambda(t)=\lambda_1-t$ on the half-line and $\lambda(t)=c-t^2$ on each line, for a fixed $c\in\R$.

The connected bounds are attained by a tree with one vertex of valence $2^n$, one finite leaf and $2^n-1$ exterior rays. It has $B_{\mathrm{nc}}=2^{n-1}-1$ and $\dim H_c^1=2^n-2$. Optimality for the differential equation requires weights realizing these models.

It remains to determine whether the lower bound of $2^{n-1}$ noncompact components is attained by some weight in $\A_n$ for every $n$, which labels can occur at the ends of the principal component for even weights, and whether the nonsymmetric ends can be connected to the local branches issuing from $\Cs$.

\section*{Acknowledgements}

I am grateful to Juli\'an L\'opez-G\'omez for taking the time to read
the manuscript and for his thoughtful comments.


\begin{thebibliography}{99}
\bibitem{AG}
H. Amann and J. L\'opez-G\'omez,
\emph{A priori bounds and multiple solutions for superlinear indefinite elliptic problems},
J. Differential Equations \textbf{146} (1998), 336--374.

\bibitem{BHJY}
D. Bartolucci, Y. Hu, A. Jevnikar and W. Yang,
\emph{Generic properties of the Rabinowitz unbounded continuum},
Adv. Nonlinear Stud. \textbf{23} (2023), Article 20220062.

\bibitem{BM}
E. Bierstone and P. D. Milman,
\emph{Semianalytic and subanalytic sets},
Publ. Math. Inst. Hautes Ètudes Sci. \textbf{67} (1988), 5--42.

\bibitem{CR}
M. G. Crandall and P. H. Rabinowitz,
\emph{Bifurcation from simple eigenvalues},
J. Funct. Anal. \textbf{8} (1971), 321--340.

\bibitem{CLT2022}
P. Cubillos, J. L\'opez-G\'omez and A. Tellini,
\emph{Multiplicity of nodal solutions in classical non-degenerate logistic equations},
Electron. Res. Arch. \textbf{30} (2022), 898--928.

\bibitem{CLT2023}
P. Cubillos, J. L\'opez-G\'omez and A. Tellini,
\emph{Global structure of the set of 1-node solutions in a class of degenerate diffusive logistic equations},
Commun. Nonlinear Sci. Numer. Simul. \textbf{125} (2023), Article 107389.

\bibitem{CLT2024}
P. Cubillos, J. L\'opez-G\'omez and A. Tellini,
\emph{High multiplicity of positive solutions in a superlinear problem of Moore--Nehari type},
Commun. Nonlinear Sci. Numer. Simul. \textbf{136} (2024), Article 108118.

\bibitem{Dancer}
E. N. Dancer,
\emph{Global structure of the solutions of non-linear real analytic eigenvalue problems},
Proc. Lond. Math. Soc. (3) \textbf{27} (1973), 747--765.

\bibitem{Diestel}
R. Diestel,
\emph{Locally finite graphs with ends: a topological approach, I. Basic theory},
Discrete Math. \textbf{311} (2011), 1423--1447.

\bibitem{FLS}
G. Feltrin, J. L\'opez-G\'omez and J. C. Sampedro,
\emph{A proof of a conjecture on the multiplicity of positive solutions of an indefinite superlinear problem},
J. Differential Equations \textbf{463} (2026), Article 114162.

\bibitem{FT}
G. Feltrin and C. Troestler,
\emph{Uniqueness, non-degeneracy, and exact multiplicity of positive solutions for superlinear elliptic problems},
J. Lond. Math. Soc. \textbf{113} (2026), Article e70476.

\bibitem{FZ}
G. Feltrin and F. Zanolin,
\emph{Multiple positive solutions for a superlinear problem: a topological approach},
J. Differential Equations \textbf{259} (2015), 925--963.

\bibitem{FL}
M. Fencl and J. L\'opez-G\'omez,
\emph{Global bifurcation diagrams of positive solutions for a class of 1D superlinear indefinite problems},
Nonlinearity \textbf{35} (2022), 1213--1248.

\bibitem{FuMcCrory}
J. H. G. Fu and C. McCrory,
\emph{Stiefel--Whitney classes and the conormal cycle of a singular variety},
Trans. Amer. Math. Soc. \textbf{349} (1997), 809--835.

\bibitem{LG2001}
J. L\'opez-G\'omez,
\emph{Spectral Theory and Nonlinear Functional Analysis},
CRC Press, Boca Raton, 2001.

\bibitem{LGMR2017}
J. L\'opez-G\'omez, M. Molina-Meyer and P. H. Rabinowitz,
\emph{Global bifurcation diagrams of one node solutions in a class of degenerate boundary value problems},
Discrete Contin. Dyn. Syst. Ser. B \textbf{22} (2017), 923--946.

\bibitem{LMZ2025}
J. L\'opez-G\'omez, E. Mu\~noz-Hern\'andez and F. Zanolin,
\emph{Global multiplicity results in a Moore--Nehari type problem with a spectral parameter},
J. Differential Equations \textbf{447} (2025), Article 113628.

\bibitem{LR15}
J. L\'opez-G\'omez and P. H. Rabinowitz,
\emph{Nodal solutions for a class of degenerate boundary value problems},
Adv. Nonlinear Stud. \textbf{15} (2015), 253--288.

\bibitem{LR16}
J. L\'opez-G\'omez and P. H. Rabinowitz,
\emph{The effects of spatial heterogeneities on some multiplicity results},
Discrete Contin. Dyn. Syst. \textbf{36} (2016), 941--952.

\bibitem{LR20}
J. L\'opez-G\'omez and P. H. Rabinowitz,
\emph{The structure of the set of 1-node solutions of a class of degenerate BVP's},
J. Differential Equations \textbf{268} (2020), 4691--4732.

\bibitem{LGB}
J. L\'opez-G\'omez and J. C. Sampedro,
\emph{Bifurcation theory for Fredholm operators},
J. Differential Equations \textbf{404} (2024), 182--250.

\bibitem{LGS}
J. L\'opez-G\'omez and J. C. Sampedro,
\emph{Global Leray--Schauder continuation for Fredholm operators},
Preprint, arXiv:2512.08412 (2025).

\bibitem{Mandel}
R. Mandel,
\emph{Global secondary bifurcation, symmetry breaking and period-doubling},
Topol. Methods Nonlinear Anal. \textbf{53} (2019), 779--800.

\bibitem{Milnor}
J. W. Milnor,
\emph{Topology from the Differentiable Viewpoint},
University Press of Virginia, Charlottesville, 1965.

\bibitem{MW}
J. F. Manning and O. H. Wang,
\emph{Cohomology and the Bowditch boundary},
Michigan Math. J. \textbf{69} (2020), 633--669.

\bibitem{PW}
J. Pejsachowicz and N. Waterstraat,
\emph{Bifurcation of critical points for continuous families of $\mathcal C^2$ functionals of Fredholm type},
J. Fixed Point Theory Appl. \textbf{13} (2013), 537--560.

\bibitem{Rabinowitz}
P. H. Rabinowitz,
\emph{Some global results for nonlinear eigenvalue problems},
J. Funct. Anal. \textbf{7} (1971), 487--513.

\bibitem{Rynne}
B. P. Rynne,
\emph{The structure of Rabinowitz' global bifurcating continua for generic quasilinear elliptic equations},
Nonlinear Anal. \textbf{32} (1998), 167--181.

\bibitem{Smale}
S. Smale,
\emph{An infinite dimensional version of Sard's theorem},
Amer. J. Math. \textbf{87} (1965), 861--866.

\bibitem{Sullivan}
D. Sullivan,
\emph{Combinatorial invariants of analytic spaces},
in \emph{Proceedings of Liverpool Singularities--Symposium I},
C. T. C. Wall (ed.), Lecture Notes in Math., vol. 192,
Springer, Berlin, 1971, pp. 165--168.

\bibitem{Tanaka}
S. Tanaka,
\emph{Morse index and symmetry-breaking for positive solutions of one-dimensional H\'enon type equations},
J. Differential Equations \textbf{255} (2013), 1709--1733.

\bibitem{Toland}
J. F. Toland,
\emph{Path-connectedness in global bifurcation theory},
Electron. Res. Arch. \textbf{29} (2021), 4199--4213.
\end{thebibliography}
\end{document}